\documentclass[11pt]{amsart}
\usepackage[margin=1.25in,bottom=1.25in]{geometry}                % See geometry.pdf to learn the layout options. There are lots.
\usepackage{graphicx}

\date{\today}        

\usepackage[colorlinks,citecolor=blue,linkcolor=blue]{hyperref}
\usepackage{mathabx}
\usepackage{amssymb}
\usepackage{mathdots}
\usepackage{caption}
\usepackage{subcaption}
\usepackage{tikz}
\usepackage{extarrows}
\usepackage[all]{xy}
\usepackage{MnSymbol}
\usepackage[justification=raggedright, singlelinecheck=false]{caption}
\usepackage{verbatim}
\usepackage{mathtools}
\usepackage{bbm}
\usepackage{amsthm}
\usepackage{nccmath}
\usepackage{amsfonts}
\usepackage[color=yellow]{todonotes}
\usepackage{pgf}
\usepackage[pscoord]{eso-pic}% The zero point of the coordinate system is the lower left corner of the page (the default).
\usepackage[nameinlink]{cleveref}

\newcommand{\placetextbox}[3]{
  \setbox0=\hbox{#3}% Put <stuff> in a box
  \AddToShipoutPictureFG*{% Add <stuff> to current page foreground
    \put(\LenToUnit{#1\paperwidth},\LenToUnit{#2\paperheight}){\vtop{{\null}\makebox[0pt][c]{#3}}}%
  }
}

\newcommand{\tx}{\text{}}
\newcommand{\p}{\text{}\\[10 pt]}
\newcommand{\pp}{\text{}\\[0 pt]}
\newcommand{\ppp}[1]{\text{}\\[#1 pt]}

\newcommand{\h}{\mathbb{H}}

\newcommand{\R}{\mathbb{R}}

\newcommand{\C}{\mathbb{C}}

\newcommand{\Z}{\mathbb{Z}}

\newcommand{\del}{\partial}

\newcommand{\qt}[2]{\small\text{\raisebox{.2ex}{$#1/_{\displaystyle #2}$\normalsize}}}   %quotient
\newcommand{\bg}{\begin}

\newcommand{\mr}{\mathrm}

\newcommand{\mc}{\mathcal}

\newtheorem{theorem}{Theorem}[section]

\theoremstyle{definition}
\newtheorem{defin}[theorem]{Definition}

\newtheorem{rem}[theorem]{Remark}

\begin{document}

%\rfoot{\thepage}
%\maketitle
\tx

\title{Experimental statistics of veering triangulations}

\author{William Worden} 
\address{Department of Mathematics, Temple University\\ Philadelphia, PA 19122}
\email{william.worden@temple.edu}

\begin{abstract}
Certain fibered hyperbolic 3-manifolds admit a \emph{layered veering triangulation}, which can be constructed algorithmically given the stable lamination of the monodromy. These triangulations were introduced by Agol in 2011, and have been further studied by several others in the years since. We obtain experimental results which shed light on the combinatorial structure of veering triangulations, and its relation to certain topological invariants of the underlying manifold.
\end{abstract}

\let\thefootnote\relax\footnote{This research was supported in part by the National Science Foundation
 through major research instrumentation grant number CNS-09-58854.}

\maketitle
\section{Introduction}
\label{sec:intro}

\p
In 2011, Agol \cite{Ag} introduced the notion of a layered veering triangulation for certain fibered hyperbolic manifolds. In particular, given a surface $\Sigma$ (possibly with punctures), and a pseudo-Anosov homeomorphism $\varphi:\Sigma\to\Sigma$, Agol's construction gives a triangulation $\mc{T}$ for the mapping torus $M_{\varphi^\circ}$ with fiber $\Sigma^\circ$ and monodromy $\varphi^\circ=\varphi_{\mid_{\Sigma^\circ}}$, where $\Sigma^\circ$ is the surface resulting from puncturing $\Sigma$ at the singularities of its $\varphi$-invariant foliations. 

For $\Sigma$ a once-punctured torus or $4$-punctured sphere, the resulting triangulation is well understood: it is the monodromy (or Floyd--Hatcher) triangulation considered in \cite{FlHa} and \cite{Jo}. In this case $\mc{T}$ is a geometric triangulation, meaning that tetrahedra in $\mc{T}$ can be realized as ideal hyperbolic tetrahedra isometrically embedded in $M_{\varphi^\circ}$. In fact, these triangulations are geometrically canonical, i.e., dual to the Ford--Voronoi domain of the mapping torus (see \cite{Ak}, \cite{La}, and \cite{Gu:thesis}). Agol's definition of $\mc{T}$ was conceived as a generalization of these monodromy triangulations, and so a natural question is whether these generalized monodromy triangulations are also geometric. Hodgson--Issa--Segerman \cite{HoIsSe} answered this question in the negative, by producing the first examples of non-geometric layered veering triangulations. In work in progress Schleimer and Segerman produce infinitely many examples of non-geometric veering triangulations using a Dehn-surgery argument.

While it is somewhat unsatisfying that these triangulations are not always geometric, they nevertheless have some nice properties. Hodgson--Rubinstein--Segerman--Tillmann \cite{HoRuSeTi:veering} show that veering triangulations admit positive angle structures, and Futer--Gu\'{e}ritaud \cite{FuGu:veering} take this a step further, by giving explicit constructions of such angle structures. In 2016, Gu\'{e}ritaud \cite{Gu:cannon} showed that there is a correspondence between the combinatorial structure of the cusp cross-section of $\mc{T}$ (lifted to the universal cover), and the tessellation of the plane induced by the image of the Cannon--Thurston map. This remarkable result extends earlier work of Cannon--Dicks \cite{CaDi} and Dicks--Sakuma \cite{DiSa}, which was in the context of punctured torus bundles. More recently, Minsky--Taylor \cite{MiTa} showed that the distance in the arc and curve complex of a subsurface $Y$ of $\Sigma^\circ$, between the stable and unstable laminations of $\varphi^\circ$ (lifted to $Y$) is coarsely bounded by $|\mc{T}|$, the number of tetrahedra in $\mc{T}$. In fact, if we fix a fibered face of the Thurston norm ball, then $\mc{T}$ is uniquely determined by this face (see \cite{Ag}), so the above statement actually holds for $\Sigma^\circ$ any fiber in the face (with corresponding laminations).

Veering triangulations have two types of combinatorially distinct tetrahedra, called \emph{hinges} and \emph{non-hinges}. 
Our main interest is in geometric aspects of these two tetrahedra types, and how their relative abundance and arrangement in the triangulation impact the geometry of the manifold. In this paper we study this relation between geometry and combinatorics via experimental methods, i.e., we use high performance computing to pseudo-randomly generate---via random walks on $\mr{Mod}(\Sigma)$ with respect to a fixed set of generators---and analyze, a large number of these layered veering triangulations. Statistics of hyperbolic 3-manifolds have previously been studied by Rivin \cite{Ri}, and by Dunfield and various collaborators (see \cite{Du:knots}, \cite{DuTh:haken}, \cite{DuTh:tunnel}). Veering triangulations are a relatively recent discovery, however, and experimental study of their statistical properties has not previously been undertaken in any significant form. Our experiments are made possible by the recent development and implementation of fast algorithms for computing veering triangulations, in the form of the computer programs \texttt{flipper}, by Bell \cite{Be}; and \texttt{Veering}, by Issa \cite{Is}. We summarize below some observations resulting from these experiments.

\subsection{Experimental Results}
First, we study the \emph{shape parameters} for the tetrahedra in $\mc{T}$. The shape parameter is a complex number which identifies the isometry class of an ideal hyperbolic tetrahedron (see \Cref{sec:shapes} for a more precise definition). Each of the tetrahedra $\Delta\subset \mc{T}$ can be assigned a shape $z_\Delta$ so that the corresponding hyperbolic tetrahedron is isometrically embedded in $M_{\varphi^\circ}$, though possibly with orientation reversed if the triangulation is not geometric (in this case we will have $\mr{Im}(z_\Delta)<0$). The first natural question to consider is: how often do we have some $\Delta \subset \mc{T}$ with negative imaginary part (i.e., when are these triangulations non-geometric)? Our data strongly suggests that in some sense these non-geometric triangulations are the generic case,leading us to make the following conjecture:

\bg{conj}
Given a surface $\Sigma$ of complexity $\xi(\Sigma)\ge 2$, the probability $P_{k,\Sigma}$ that a simple random walk of length $k$ on a set of generators for $\mr{Mod}(\Sigma)$ yields a geometric veering triangulation, decays to 0 like $e^{-kc}$ for some constant $c=c(\Sigma)$.

\end{conj}

Contrasting this observation, in \Cref{sec:others} we give an infinite family of veering triangulations that we propose, based on strong supporting experimental evidence, may all be geometric.

The shape parameters mentioned in the preceding paragraph are elements of a subset $B \subset \C$---the intersection of two unit discs centered at $0$ and $1$ on the real line. In \Cref{sec:shapes} we study the distribution in $B$ of these shapes, for randomly generated monodromies $\varphi:\Sigma_{g,n}\to \Sigma_{g,n}$ for fixed genus $g$ and $n$ punctures. We find that for hinges and non-hinges, the distributions are quite different, with non-hinges tending to be ``flatter" (i.e., their shape parameters have small imaginary part). We also observe a marked dependence of the shape distributions of both hinges and non-hinges on the genus $g$, while in most cases the distributions seem to be far more subtly dependent on the number of punctures $n$.

In \Cref{sec:volume} we observe a coarse linear relation between the number of hinge tetrahedra in $\mc{T}$, and the volume of the underlying hyperbolic manifold. The slope of this linear relation appears to depend on both the genus and number of punctures of the surface $\Sigma$. It is generally believed that this slope should be bounded below, independent of $\Sigma$, and our data bears this out. Whether there is a universal upper bound is not clear from our data.

Finally, in \Cref{sec:systole}, we investigate what the arrangement of non-hinge tetrahedra in $\mc{T}$ can tell us about the geometry of $M_{\varphi_\circ}$. In particular, our data suggests that long chains of non-hinge tetrahedra, with respect to an ordering on the tetrahedra in $\mc{T}$ to be defined, give short geodesics. These experiments are motivated by the Length Bound Theorem of Brock--Canary--Minsky \cite{BCM}, and our results are consistent with the theorem's implications.

The paper is arranged as follows. First, we review some necessary background material in \Cref{sec:bg}, and discuss experimental methodology in \Cref{sec:method}. In \Cref{sec:nongeo}, we give evidence supporting a conjecture that non-geometricity of veering triangulations is in some sense generic, when the complexity of $\Sigma$ is greater than one. In \Cref{sec:shapes} we analyze the distribution of shape parameters of tetrahedra of veering triangulations, observing an apparent dependence on the genus of the fiber of the mapping torus. In \Cref{sec:volume} we give evidence of a course linear relation between volume $\mr{Vol}(\mc{T})$ of a mapping torus and the number of hinges $|\mathcal{T}|_H$ in the veering triangulation. In \Cref{sec:systole} we explore a possible relation between systole length and non-hinge tetrahedra in the veering triangulation. In \Cref{sec:others} we discuss a few other experimental observations of note. Finally, in \Cref{sec:supp}, we direct the reader to supplemental materials available on the author's website, including the full dataset on which all results depend, all scripts used for the experiments, and additional figures which are not included herein.

\subsection{Acknowledgements}
\label{ssec:ack}

We are very grateful to Mark Bell for many helpful conversations about veering triangulations and his program \texttt{flipper}, and for excellent technical support for \texttt{flipper}! We also thank Saul Schleimer and Sam Taylor for illuminating conversations about veering triangulations, which at times helped guide further inquiry. Most of all we would like to thank Dave Futer for frequent conversations which were invaluable to the success of this project, and for much helpful guidance during the writing process. Additionally, we thank the referee for many helpful comments and suggestions.

\begin{figure}[h!]
        \centering
        \includegraphics[scale=.85]{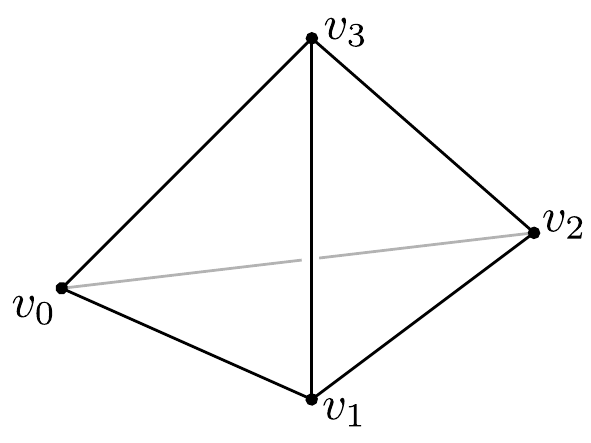}
        \caption{The standard 3-simplex}
        \label{fig:std_simplex}
\end{figure}

\section{Background}
\label{sec:bg}

Let $M$ be an orientable cusped hyperbolic $3$-manifold, i.e., $M$ is the interior of a compact orientable $3$-manifold $\overline{M}$ with at least one torus boundary component. Given a tessellation $\overline{\mathcal{T}}$ of $\overline{M}$ by truncated tetrahedra (i.e., tetrahedra with their corners chopped off), such that each triangular face of a truncated corner of a tetrahedron lies in $\del \overline{M}$, we get an \textbf{ideal triangulation} $\mathcal{T}$ of $M$ by taking the interior of $\overline{\mathcal{T}}$. The following more precise definition will also be useful:

\bg{defin}
An \textbf{ideal triangulation} $\mathcal{T}$ is a disjoint union $\bigsqcup_j \Delta_j$ of standard ideal 3-simplices, along with a surjection $i:\bigsqcup_j \Delta_j \to M$ which is an embedding on the interior  of each 3-simplex, and preserves the orientation given by the vertex ordering of the standard simplex (see \Cref{fig:std_simplex}). We will often identify $\mathcal{T}$ with the image $i(\bigsqcup_j \Delta_j)$. 

\end{defin}

The triangulation $\mathcal{T}$ consists of \textbf{ideal tetrahedra}: that is, tetrahedra with their vertices removed. Note that each tetrahedron $\Delta$ in $\mathcal{T}$ inherits an ordering $\{v_0,v_1,v_2,v_3\}$ of its vertices from the standard simplex.

 In this paper we will be interested in a certain type of ideal triangulation, called a \textbf{veering triangulation}, introduced by Agol \cite{Ag}. Before describing this triangulation, we make a few preliminary definitions.

\begin{figure}
        \begin{subfigure}[b]{0.35\textwidth}
                \includegraphics[scale=.75]{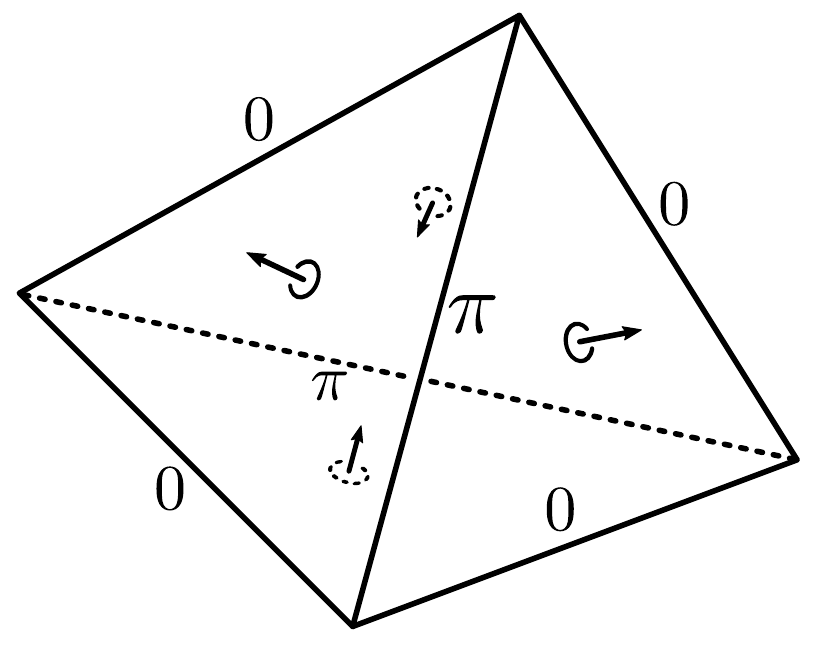}
                \caption{}
                \label{fig:taut_tet}
        \end{subfigure}
       \qquad\quad
        \begin{subfigure}[b]{0.5\textwidth}
                \includegraphics[scale=.65]{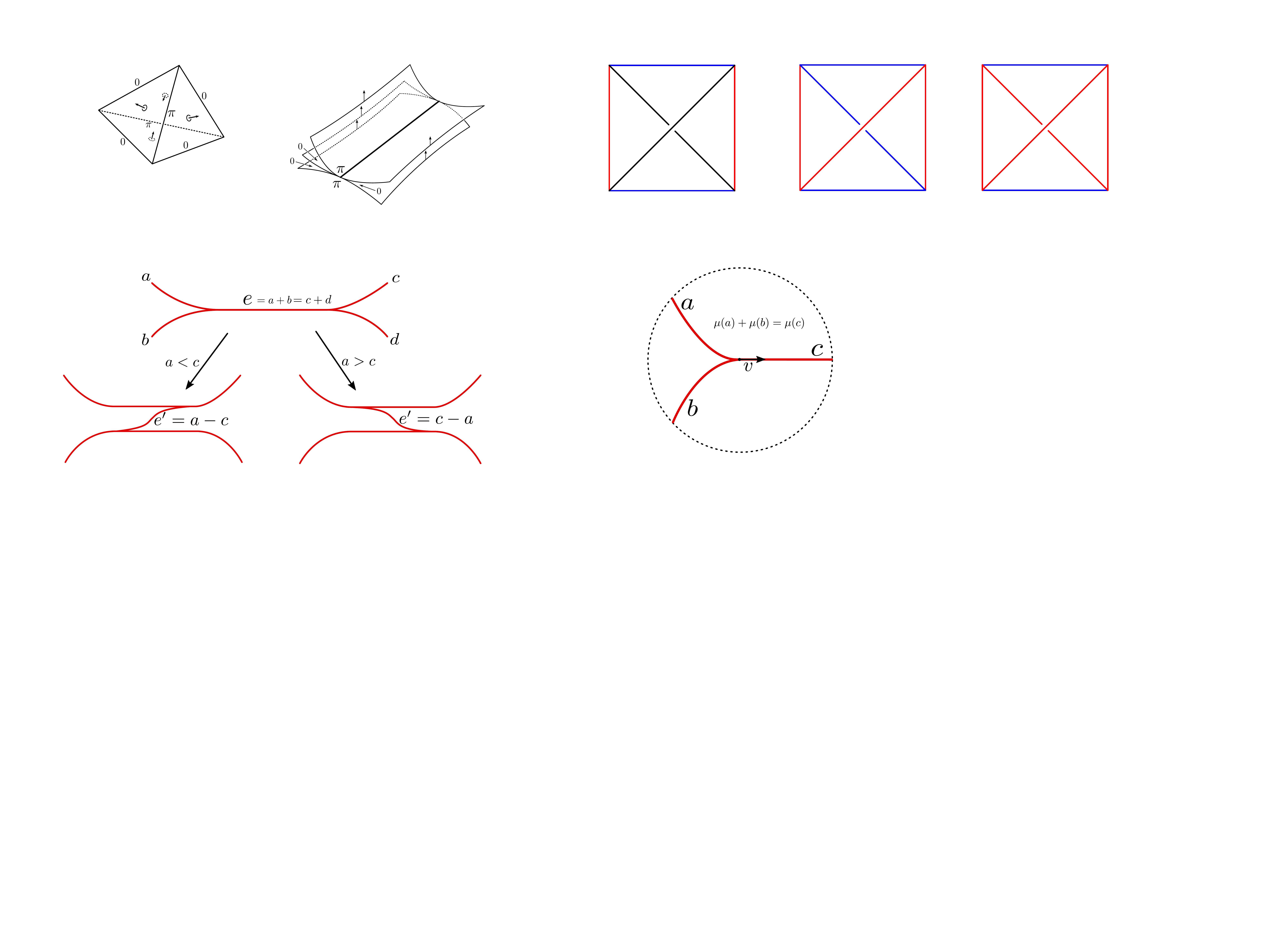}
                \caption{}
                \label{fig:taut_edge}
        \end{subfigure}       
        \caption{Left: A taut tetrahedron. Right: The angles around an edge of a taut tetrahedron sum to $2\pi$.}
        \label{fig:taut}
\end{figure}

\bg{defin}
A \textbf{taut ideal tetrahedron} $\Delta$ is an ideal tetrahedron such that each face of $\Delta$ has a coorientation, with two faces oriented into $\Delta$ and two oriented outward. In addition, we assign $\Delta$ dihedral angles as follows: an edge $e$ has dihedral angle $\pi$ if the two faces incident on $e$ have the same coorientation, and the dihedral angle of $e$ is $0$ otherwise (see \Cref{fig:taut_tet}). A \textbf{taut ideal triangulation} of a 3-manifold $M$ is an ideal triangulation $\mc{T}$, such that the faces of $\mc{T}$ can be cooriented so that every tetrahedron is taut, with the dihedral angles around each edge summing to $2\pi$ (see \Cref{fig:taut}).
\end{defin}

\begin{figure}
        \centering
        \begin{subfigure}[b]{0.22\textwidth}
                \includegraphics[scale=.37]{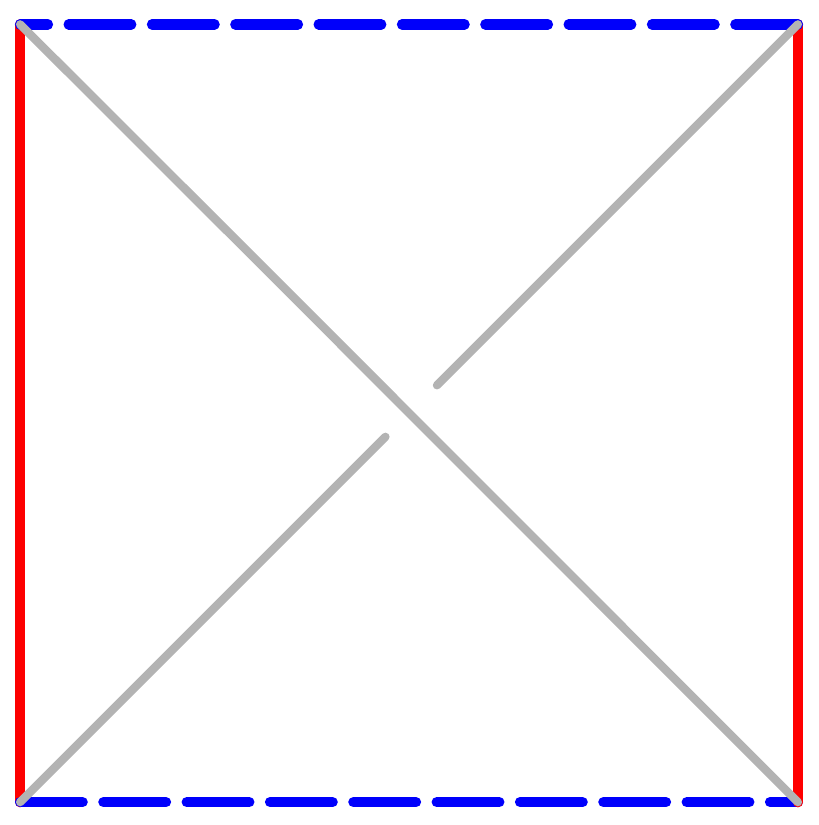}
                \caption{}
                \label{fig:veering_tet}
        \end{subfigure}
       \quad
        \begin{subfigure}[b]{0.22\textwidth}
                \includegraphics[scale=.37]{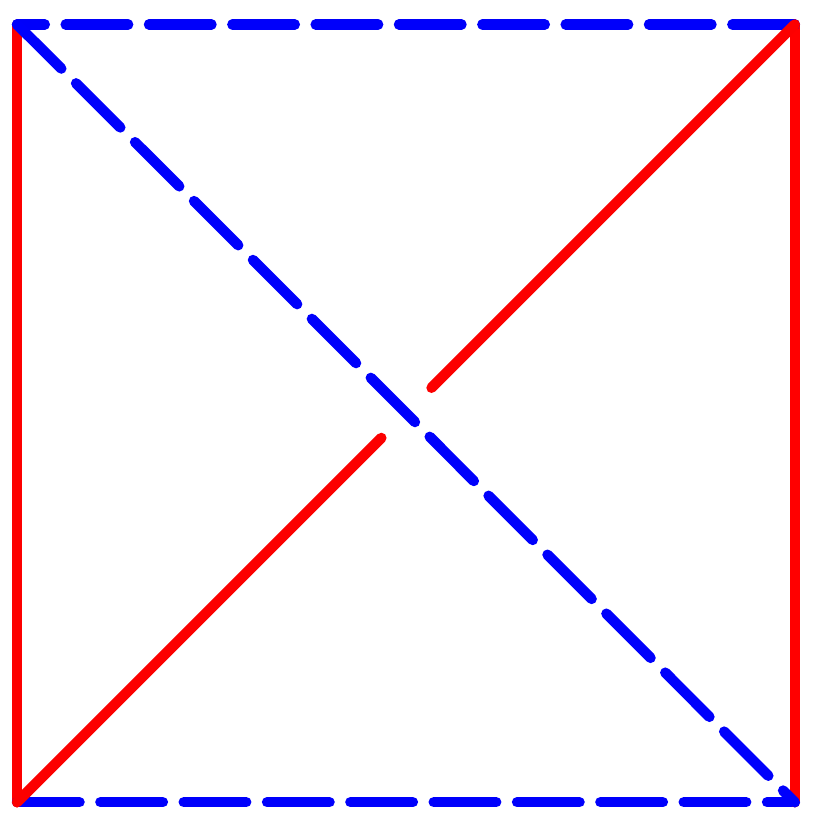}
                \caption{}
                \label{fig:hinge}
        \end{subfigure}    
         \quad
        \begin{subfigure}[b]{0.22\textwidth}
                \includegraphics[scale=.37]{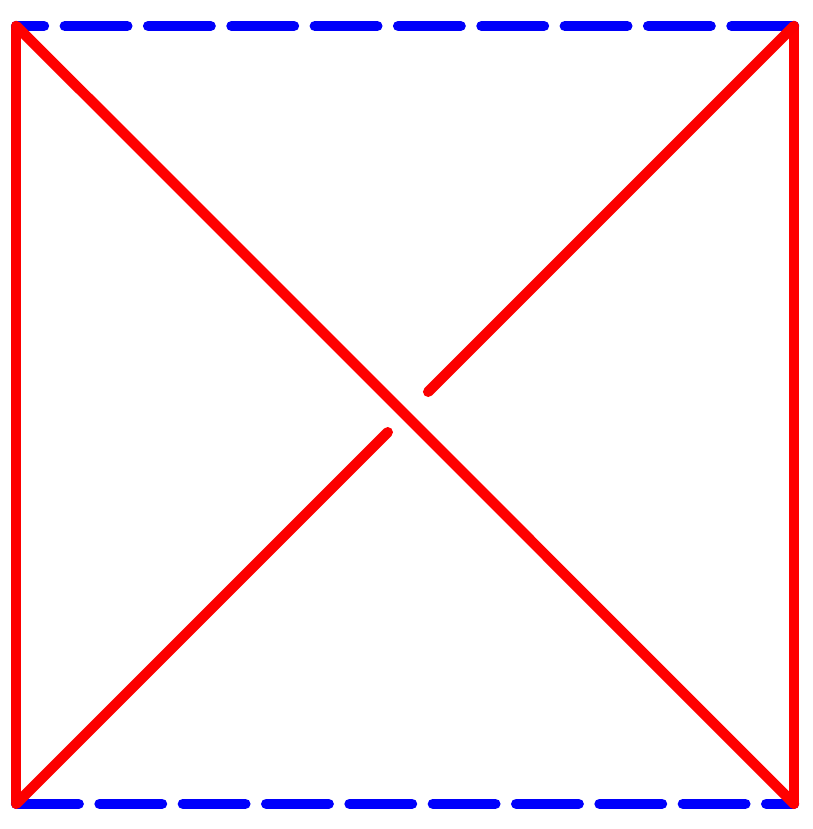}
                \caption{}
                \label{fig:non-hinge}
        \end{subfigure}
        \quad
        \begin{subfigure}[b]{0.22\textwidth}
                \includegraphics[scale=.5]{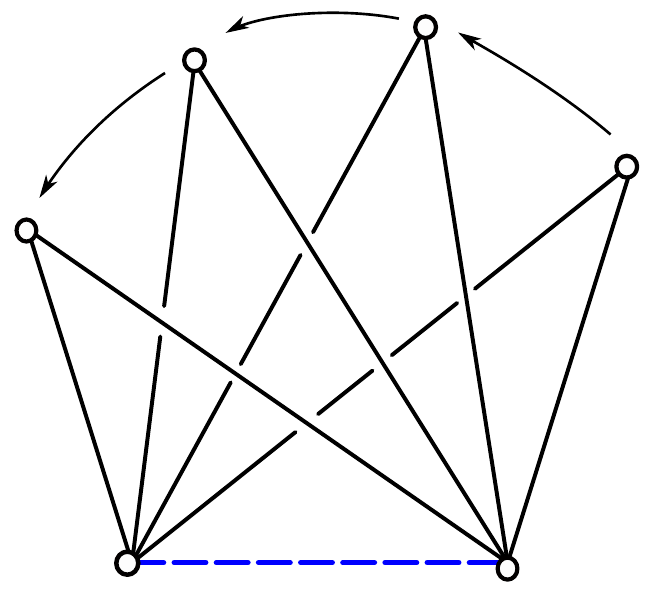}
                \caption{}
                \label{fig:left_veering}
        \end{subfigure}       
        \caption{Left: A veering tetrahedron, flattened on the page, with angle$-\pi$ edges on the diagonals. Edges with dihedral angle $\pi$ will inherit the color of the angle-$0$ edges that glue to it. Center-left: A hinge tetrahedron. Center-right: A non-hinge tetrahedron. Right: A left veering edge.}
        \label{fig:veer_tets}
\end{figure}

\bg{defin}\label{veering}
A \textbf{veering tetrahedron} $\Delta$ is a taut ideal tetrahedron with each edge colored either red (solid) or blue (dashed), such that $\Delta$ is combinatorially the same as the tetrahedron shown in \Cref{fig:veering_tet} (here we have drawn $\Delta$ flattened onto the page, with angle-$\pi$ edges on the diagonals, and angle-$0$ edges on the top, bottom, and sides). That is, at any corner of $\Delta$, the counter-clockwise ordering on the edges (as we look in from the cusp) is red angle-0, angle-$\pi$, blue angle-0. The angle-$\pi$ edges can be colored any combination of red and blue, and for this reason they are left uncolored in \Cref{fig:veering_tet}. A \textbf{veering triangulation} of a 3-manifold $M$ is an ideal triangulation $\mc{T}$ such that every tetrahedron is veering (note that the coloring for an angle-$\pi$ edge will be determined by the color of the angle-0 edges which are glued to it).
\end{defin}

In a veering triangulation $\mc{T}$, edges that are colored blue are called \textbf{left-veering}, and edges that are colored red are called \textbf{right-veering}. This terminology comes from Agol's original definition of the veering condition: for Agol, an edge $e$ is left- (resp. right-) veering if the triangles incident to it veer left (resp. right), with respect to the ordering induced by the transverse orientation, as shown in \Cref{fig:left_veering}. A taut triangulation is veering if every edge in the triangulation is either left- or right-veering. \Cref{veering} above is due to Hodgson--Rubinstein--Segerman--Tillmann \cite{HoRuSeTi:veering}, and is shown by the same authors to be equivalent to Agol's definition, in the case of taut triangulations.

Both Agol's definition and \Cref{veering} are useful perspectives when working with veering triangulations. There is a third definition, due to Gu\'{e}ritaud \cite{Gu:cannon}, in terms of the quadratic differential associated to a pseudo-Anosov $\varphi$, which is also quite useful and offers further insight into the structure of veering triangulations of mapping tori. This definition is somewhat more technical, so we will not include it here, but the interested reader may refer to \cite{Gu:cannon} and \cite{MiTa} for details.

The definition of a veering tetrahedron allows for three different color configurations for the edges, if we ignore face co-orientation: either both angle-$\pi$ edges are blue, or both are red, or one is red and the other is blue. This leads to the following definition:

\bg{defin}
If the two angle-$\pi$ edges of $\Delta$ are colored differently, then we say that $\Delta$ is a \textbf{hinge} tetrahedron. Otherwise, $\Delta$ is a \textbf{non-hinge} tetrahedron (see \Cref{fig:hinge} and \Cref{fig:non-hinge}).
\end{defin}

The veering triangulations constructed by Agol are built by layering taut tetrahedra onto the fiber (a punctured surface) of a mapping torus---we will describe this in more detail below. In this case the veering triangulations are called \textbf{layered}. In the same paper, Agol asked whether there might be veering triangulations which are \emph{not} layered. Indeed, it has been shown by Hodgson--Rubinstein--Segerman--Tillmann \cite{HoRuSeTi:veering} that, in general, veering triangulations may not be layered, and in fact there are triangulations that satisfy a weaker version of the veering condition, called veering angle-taut triangulations.

\begin{figure}
        \centering
        \begin{subfigure}[b]{0.3\textwidth}
                \includegraphics[scale=.35]{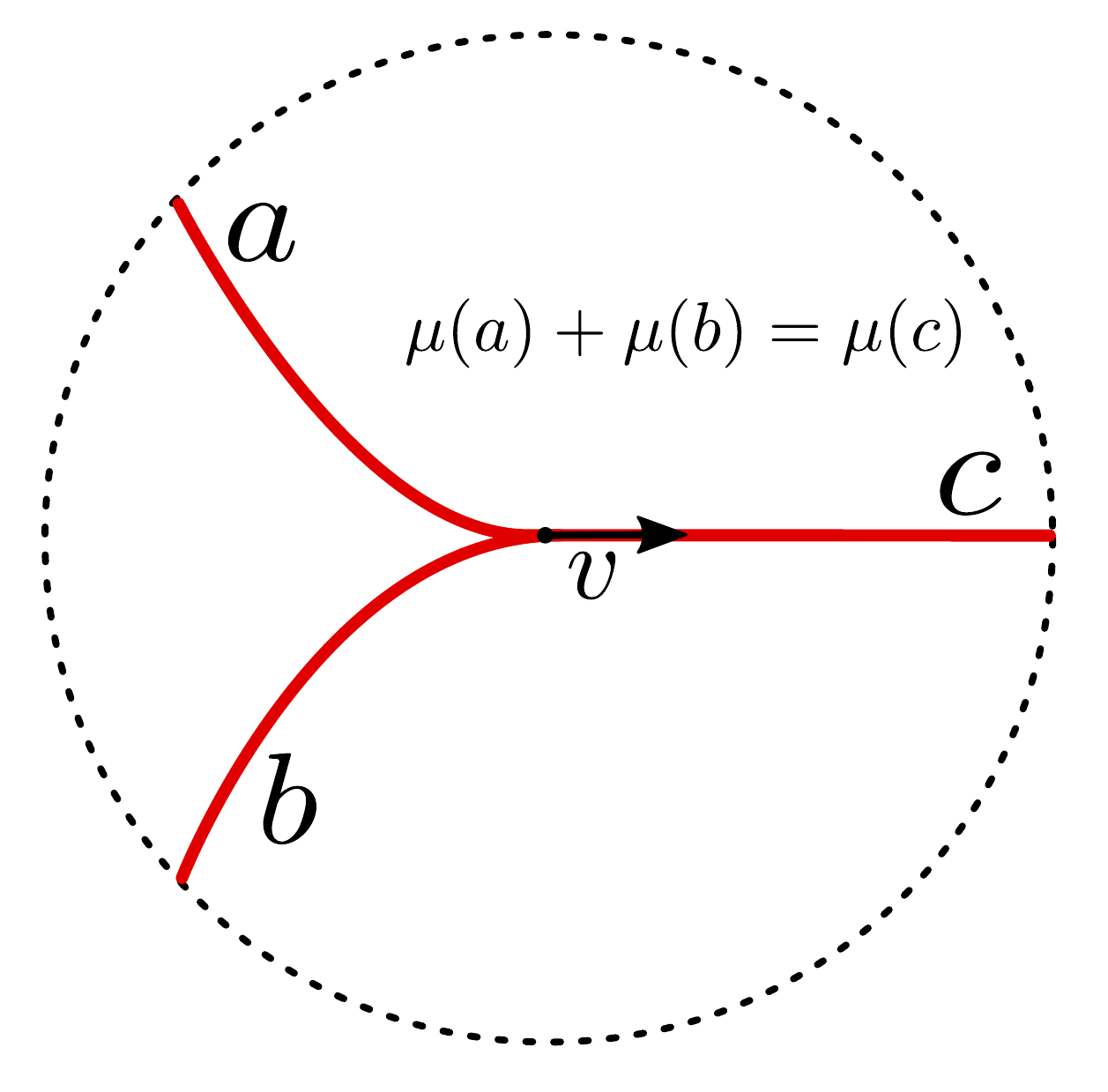}
                \caption{}
                \label{fig:branch}
        \end{subfigure}
       \qquad
        \begin{subfigure}[b]{0.63\textwidth}
                \includegraphics[scale=.35]{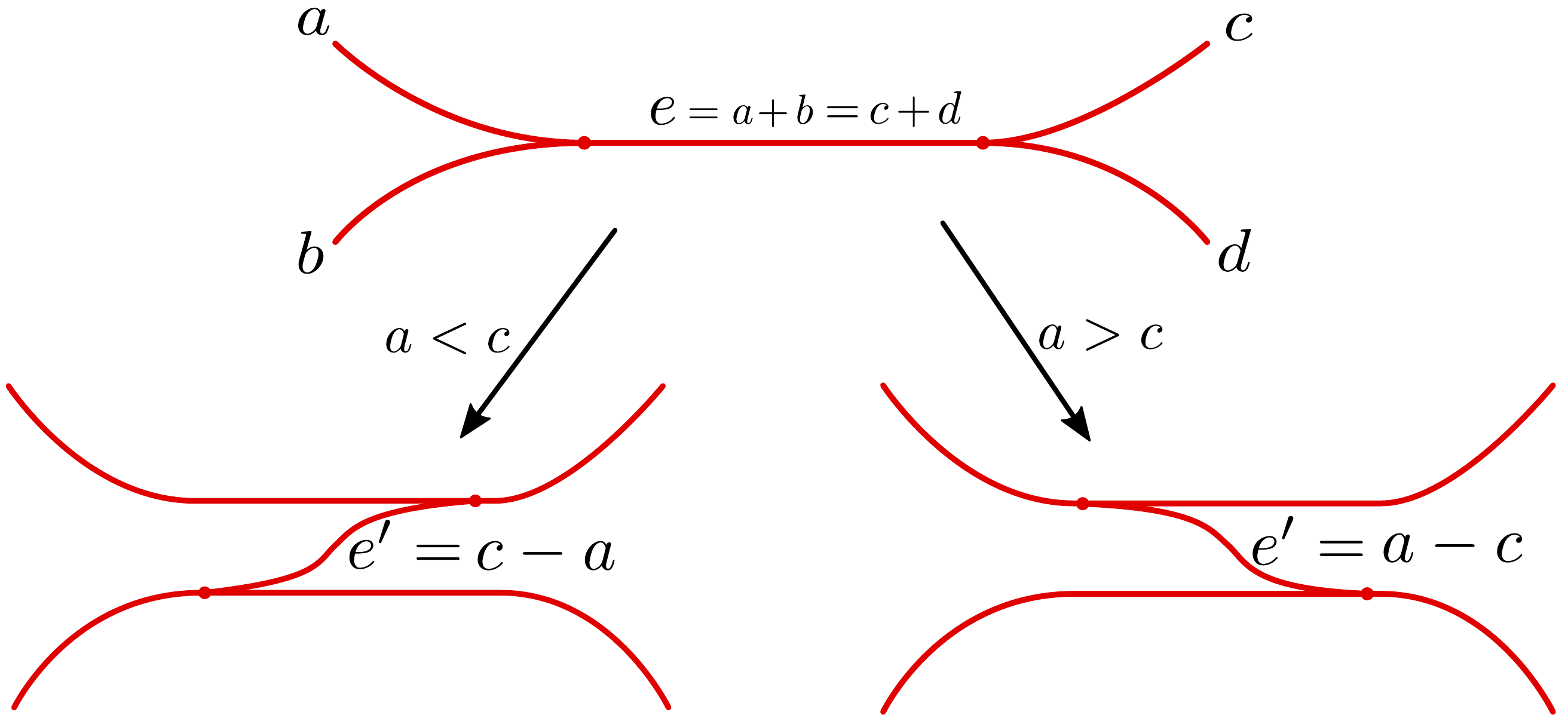}
                \caption{}
                \label{fig:splitting}
        \end{subfigure}       
        \caption{Left: A branch of a measured train track. Right: Splitting a measured train track (here we have suppressed notation, writing $a$ instead of $\mu(a)$, etc.)}
        \label{fig:train_track}
\end{figure}

We will now describe Agol's veering triangulation construction. Let $\Sigma$ be a surface, possibly with punctures, with Euler characteristic $\chi(\Sigma)<0$, and let $\varphi:\Sigma \to \Sigma$ be a pseudo-Anosov mapping class. Associated to the mapping class $\varphi$, we have the \emph{stable} and \emph{unstable} measured geodesic laminations $\mc{L}^+$ and $\mc{L}^-$, which satisfy $\varphi(\mathcal{L}^+)=\lambda_\varphi \mathcal{L}^+$ and $\varphi(\mathcal{L}^-)=\lambda_\varphi^{-1} \mathcal{L}^-$, for some $\lambda_\varphi \in (1,\infty)$. The number $\lambda_\varphi$ is called the \textbf{dilatation} of $\varphi$. Let $\Sigma^\circ$ be the surface resulting from puncturing each complementary region of $\mc{L}^+$ in $\Sigma$ (if it is not already punctured). The restriction $\varphi^\circ:=\varphi_{\mid_{\Sigma^\circ}}$ of $\varphi$ to $\Sigma^\circ$ is a pseudo-Anosov mapping class for our new surface $\Sigma^\circ$, having the same associated measured laminations $\mathcal{L}^+$ and $\mathcal{L}^-$ and the same dilatation $\lambda_\varphi$ as $\varphi$. Let $M_{\varphi^\circ}=\qt{\Sigma^\circ\times I}{\{(x,0)=(\varphi^\circ(x),1)\}}$ be the mapping torus with fiber $\Sigma^\circ$ and monodromy $\varphi^\circ$. 
To describe the veering triangulation of $M_{\varphi^\circ}$, we first need the following definition of Thurston \cite{Th}:

\bg{defin}
A \textbf{train track} $\tau$ on a surface $\Sigma$ is a trivalent graph embedded in $\Sigma$ so that each vertex has a well defined tangent direction,and such that $\chi(D(C))<0$ for each component $C$ of $\Sigma\setminus \tau$, where $D(C)$ is the double of $C$ across arcs of $\del C$ parallel to branches of $\tau$ (this means $D(C)$ has punctures at non-smooth points of $\del C$). Vertices of $\tau$ are called \textbf{switches}, and edges are called \textbf{branches}. A \textbf{measured train track} $(\tau,\mu)$ on $\Sigma$ is a train track $\tau$ along with a weight $\mu(e)\in\R$ associated to each branch $e$ of $\tau$, such that the weights satisfy the \textbf{switch condition}: if $a,b,$ and $c$ are branches of $\tau$ that meet at a switch $s$ as shown in \Cref{fig:branch}, then $\mu(a)+\mu(b)=\mu(c)$.
\end{defin}

A lamination $\mathcal{L}$ on $\Sigma$ is said to be \textbf{carried by} $\tau$ if there is a differentiable map $f:\Sigma \to \Sigma$ that is homotopic to the identity, non-singular on the tangent spaces of the leaves of $\mathcal{L}$, and satisfies $f(\mathcal{L})\subset \tau$. If $\mc{L}$ is a measured lamination, then $f$ induces a measure $\mu$ on $\tau$: given a branch $e\in \tau$, $\mu(e)$ is the $\mc{L}$-measure of the segment $f^{-1}(p)$ transverse to $\mc{L}$, where $p$ is a point in the interior of the branch $e$. We will say a measured lamination $\mc{L}$ is carried by the measured train track $(\tau,\mu)$ if $\mc{L}$ is carried by $\tau$ and $\mu$ is induced by the measure on $\mc{L}$. With these definitions Thurston provides a powerful tool for understanding laminations on surfaces---in particular, we can now replace a relatively complicated object, namely the measured lamination $\mc{L}$, with a comparatively simple combinatorial object.

Now, suppose we have a measured train track $(\tau,\mu)$ on $\Sigma^\circ$, such that the stable lamination $\mc{L}^+$ of $\varphi^\circ:\Sigma^\circ\to\Sigma^\circ$ is carried by $(\tau,\mu)$. We can obtain a new train track $(\tau',\mu')$ that also carries $\mu({L}^+)$ by performing the move described in \Cref{fig:splitting}, called a \textbf{splitting}. If $(\tau',\mu')$ is obtained from $(\tau,\mu)$ by splitting along all branches of $\tau$ with maximal weight, then we denote this \textbf{maximal splitting} by $(\tau,\mu)\rightharpoonup (\tau',\mu')$. If we have a sequence of maximal splittings (i.e., a \textbf{splitting sequence}), then we'll denote this by $(\tau_1,\mu_1)\rightharpoonup^n (\tau_n,\mu_n):=(\tau_1,\mu_1)\rightharpoonup (\tau_2,\mu_2)\rightharpoonup \dots \rightharpoonup (\tau_n,\mu_n)$. The following theorem of Agol will allow us to build a veering triangulation of the mapping torus $M_{\varphi^\circ}$ defined above, with fiber $\Sigma^\circ$ and monodromy $\varphi^\circ$: 

\begin{figure}
        \centering
        \begin{subfigure}[b]{0.6\textwidth}
                \includegraphics[scale=.35]{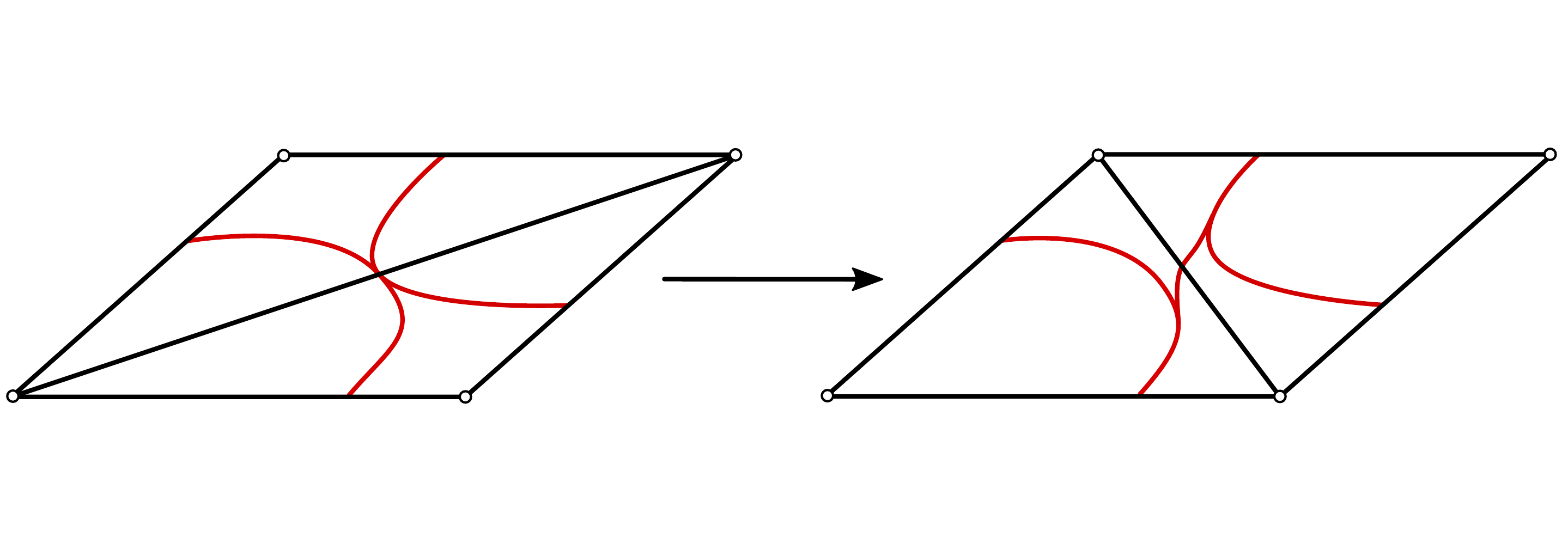}
                \caption{}
                \label{fig:dual_split}
        \end{subfigure}
       \qquad\quad
        \begin{subfigure}[b]{0.3\textwidth}
                \includegraphics[scale=.35]{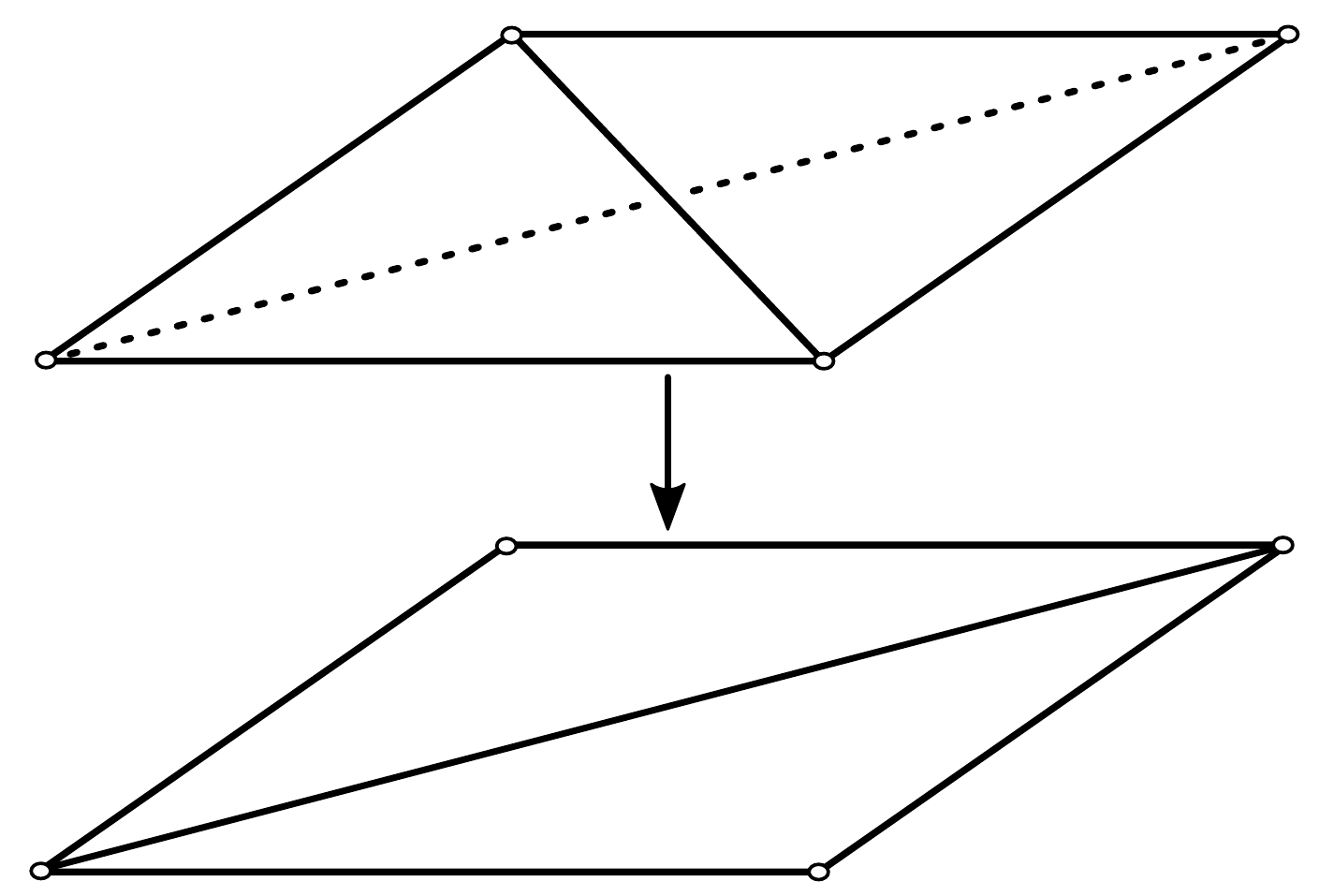}
                \caption{}
                \label{fig:whitehead}
        \end{subfigure}       
        \caption{Left: A splitting of the train track $\tau$ results in a diagonal exchange in the dual triangulation. Right: Each diagonal exchanged can be thought of as layering on a flat tetrahedron.}
        \label{fig:whitehead_split}
\end{figure}

\bg{thm}\cite[Theorem 3.5]{Ag}
\label{ag_thm}
If $\varphi:\Sigma\to\Sigma$ is a pseudo-Anosov map, with stable lamination $\mc{L}^+$ carried by $(\tau_1,\mu_1)$, then there exists $n,m$ such that 
$$
(\tau_1,\mu_1)\rightharpoonup^n (\tau_n,\mu_n)\rightharpoonup^m (\tau_{n+m},\mu_{n+m}),
$$
and $\tau_{n+m}=\varphi(\tau_n)$ and $\mu_{n+m}=\lambda_\varphi^{-1} \varphi(\mu_n)$.
Furthermore, the periodic portion of the sequence $(\tau_n,\mu_n)\rightharpoonup^m (\tau_{n+m},\mu_{n+m})$ is unique (i.e., it depends only on the conjugacy class of $\varphi$, and not on the initial choice of $(\tau_1, \mu_1)$).
\end{thm}

This says that if we start with any $(\tau,\mu)$ carrying $\mc{L}$ and repeatedly apply maximal splittings, the sequence will eventually be cyclic modulo the action of the monodromy. Recall that $\Sigma^\circ$ is the surface that results from puncturing the complementary regions of the stable lamination. Let $(\tau_i,\mu_i)$ be as in the theorem. If we consider $\tau_i$ as a train track on $\Sigma^\circ$, then the fact that the stable lamination $\mc{L}^+$ of $\varphi$ is carried by $\tau_i$ implies that every connected component of $\Sigma^\circ \setminus \tau_i$ will be a disc with a single puncture. If we take these punctures as the vertices of the dual graph of $\tau_i$, then this graph gives the edges of a triangulation $T_i$ of $\Sigma^\circ$, since $\tau_i$ is trivalent. When we apply the splitting $(\tau_i,\mu_i)\rightharpoonup (\tau_{i+1},\mu_{i+1})$, the dual triangulation is changed by a diagonal exchange, as shown in \Cref{fig:whitehead} (note that two maximal branches of a train track cannot share a vertex, so even if a maximal splitting involves multiple maximal weight branches, we can consider the branches individually). If we think of each diagonal exchange as layering onto $\Sigma^\circ$ a flat tetrahedron (see \Cref{fig:dual_split}), then the splitting sequence $(\tau_n,\mu_n)\rightharpoonup^m (\tau_{n+m},\mu_{n+m})$ corresponds to adding successive layers of tetrahedra to $\Sigma^\circ$. Since $\tau_{n+m}=\varphi(\tau_n)$ and $\mu_{n+m}=\lambda_\varphi^{-1}(\mu_n)$, we can glue $\Sigma^\circ \times \{0\}$ to $\Sigma^\circ\times \{1\}$ via $(x,0)\sim (\varphi^\circ(x),1)$, so that $\varphi^\circ(T_n)$ is glued to $T_{n+m}$. The fact that $\varphi^\circ$ is a pseudo-Anosov map is enough to guarantee that this process gives a triangulation $\mc{T}$ of the mapping torus $M_{\varphi^\circ}$, since no power of $\varphi^\circ$ can fix an edge of the triangulation $T_n$.

%Every edge of $\tau_n$ eventually splits at some step of this sequence (see \cite[Lemma 2.1]{Ag}), so we get a 3-dimensional triangulation $\mc{T}'$ of $\Sigma^\circ \times I$. Since $\tau_{n+m}=\varphi(\tau_n)$ and $\mu_{n+m}=\lambda_\varphi^{-1}(\mu_n)$, we can glue $\Sigma^\circ \times \{0\}$ to $\Sigma^\circ\times \{1\}$ via $(x,0)\sim (\varphi(x),1)$, so that $\varphi(T_n)$ is glued to $T_{n+m}$. Hence we get a triangulation $\mc{T}$ of the mapping torus $M_{\varphi^\circ}$.

By construction, this triangulation is layered. It is also quite clear that $\mc{T}$ is taut, since the layering gives a natural way to co-orient the faces and assign appropriate dihedral angles to the edges. Furthermore, Agol proves that $\mc{T}$ is in fact a veering triangulation \cite[Proposition 4.2]{Ag}, which by \Cref{ag_thm} above is an invariant of the fibration. All of the triangulations considered in this paper will be constructed as above, and therefore will be layered veering triangulations. We will often drop the adjective ``layered," and simply refer to our triangulations as veering triangulations. 

We close this section with another definition, which will be useful going forward:

\begin{defin}
The \textbf{complexity} of a surface $\Sigma$ is defined to be $\xi(\Sigma)=3g-3+n$, where $g$ is the genus and $n$ is the total number of punctures and boundary components of $\Sigma=\Sigma_{g,n}$.
\end{defin}

\section{Methodology}
\label{sec:method}

\begin{figure}
        \centering
        \begin{subfigure}[b]{0.35\textwidth}
                \includegraphics[scale=.25]{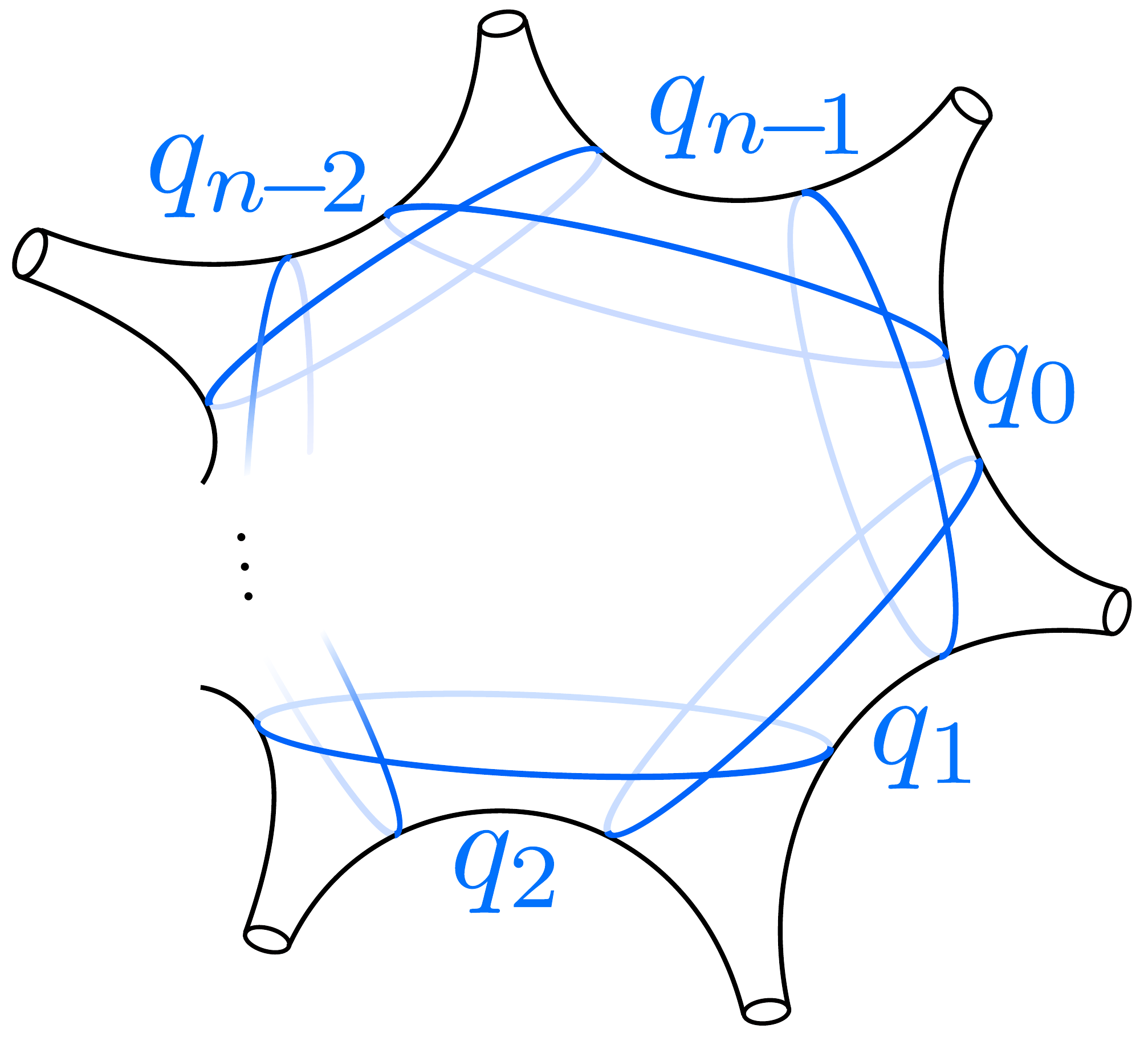}
        		\caption{}
        		\label{fig:S0n_gens}
       \end{subfigure}
       \qquad
        \begin{subfigure}[b]{0.55\textwidth}
                \includegraphics[scale=.35]{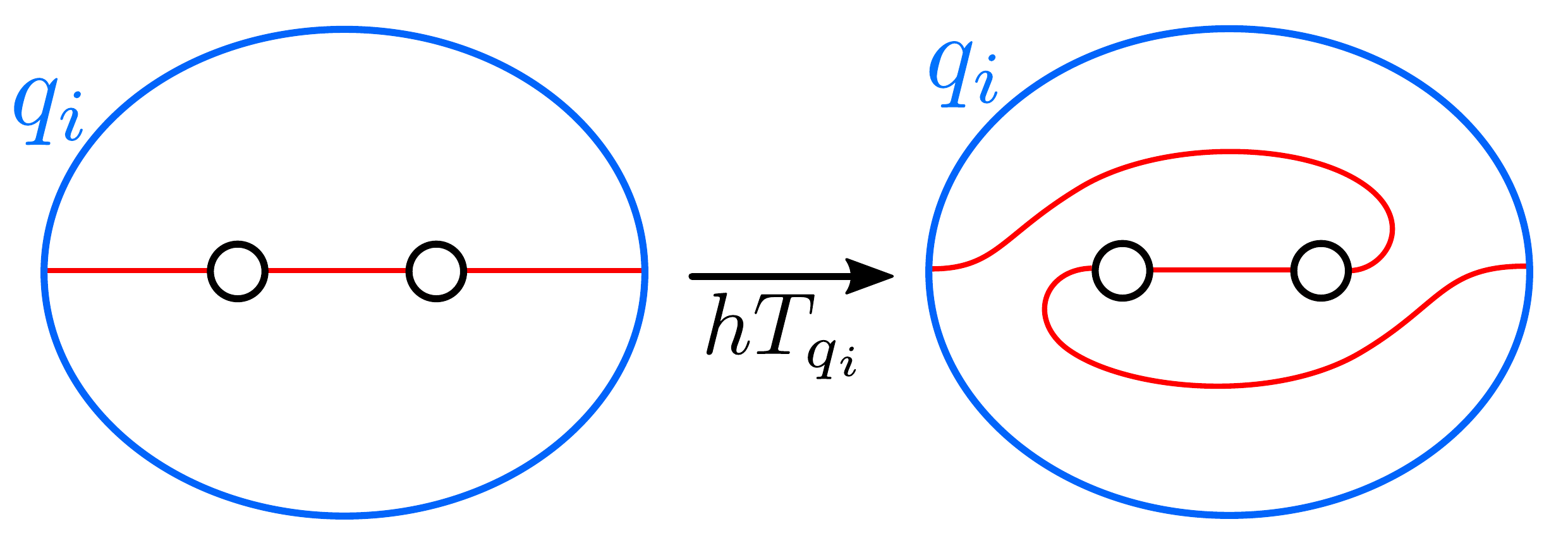}
                \caption{}
                \label{fig:half_twist}
        \end{subfigure}       
       \caption{Left: The mapping class group of the $n$-punctured sphere $\Sigma_{0,n}$, $n\ge 4$, is generated by half-twists about the $q_i$ curves shown. A half twist $hT_{q_i}$ about a curve $q_i$ is a homeomorphism that is the identity on the complement of the twice punctured disc bounded by $q_i$, and transposes the punctures in this disc, as shown in the right frame.}
        \label{fig:S0n}
\end{figure}

\begin{figure}
        \centering
        \includegraphics[scale=.38]{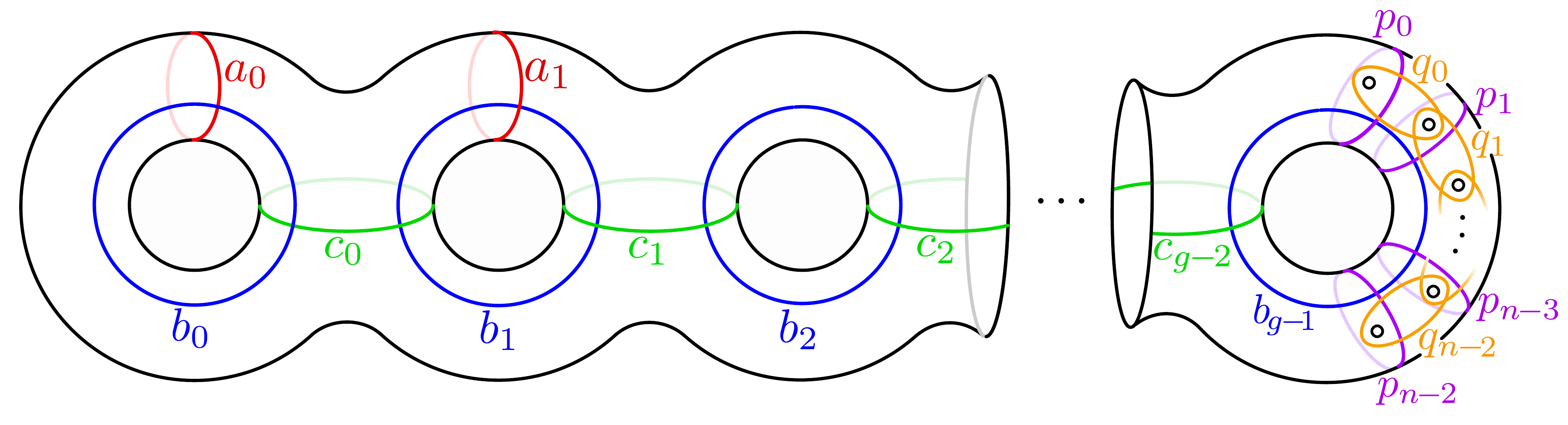}
        \caption{The mapping class group $\mr{Mod}(\Sigma_{g,n})$, for $g>0$, is generated by full Dehn-twists about the $a_i$, $b_i$, $c_i$, and $p_i$ curves, and half-twists about the $q_i$ curves. If the $q_i$ curves are omitted, than Dehn twists about the remaining curves generate the pure mapping class group $\mr{PMod}(\Sigma_{g,n})$.}
        \label{fig:surf_gens}
\end{figure}

To construct examples of veering triangulations following the above method, we use the computer program \texttt{flipper} \cite{Be}, written by Mark Bell. Let $\Sigma=\Sigma_{g,n}$ be a punctured surface of genus $g$, with $n>0$ punctures. We get a mapping class $\varphi\in\mr{Mod}(\Sigma)$ by pseudo-randomly selecting $l(\varphi)$ letters, one at a time, from a set of generators for $\mr{Mod}(\Sigma)$. In other words, $\varphi$ is the result of a random walk on the Cayley graph of $\mr{Mod}(\Sigma)$. For these random walks we will use the generating sets shown in Figures \Cref{fig:S0n} and \Cref{fig:surf_gens}, depending on whether $g=0$, or $g>0$, respectively. In some cases we will want to sample from only the pure mapping class group $\mr{PMod}(\Sigma)$, in which case the generators $q_i$ in \Cref{fig:surf_gens} are omitted. We also note that in \Cref{sec:systole} we will use a somewhat more restrictive notion of random word, to be described in that section.

 Given the mapping class $\varphi$, \texttt{flipper} constructs a veering triangulation $\mc{T}$ of the mapping torus $M_{\varphi^\circ}$, as described above (assuming that $\varphi$ was pseudo-Anosov). Recall that the fiber $\Sigma^\circ$ of $M_{\varphi^\circ}$ is the surface resulting from puncturing $\Sigma$ on the complementary regions of the stable lamination $\mc{L}^+$ (if such a region is not already punctured). If we Dehn-fill every cusp coming from a punctured singularity along its fiber slope (i.e., the slope parallel to the fibers), then the result is the mapping torus $M_\varphi$ with fiber $\Sigma$. It makes sense then to refer to $\Sigma$ as the \textbf{filled fiber}. In most cases this will not be the same as the fiber $\Sigma^\circ$, as follows from the Poincar\'{e}-Hopf theorem and a result of Gadre--Maher \cite{GaMa} which says that generic mapping classes have 1-pronged singularities at punctures and 3-pronged singularities elsewhere.

In the case when $\Sigma$ has $n$ punctures, the action of the mapping class group on the punctures of $\Sigma$ is by the symmetric group $S_n$. And when $g=0$, all of the generators of $\mr{Mod}(\Sigma)$ act as transpositions of adjacent punctures. Given a word $\varphi$ in these generators, let $\sigma_\varphi$ be the corresponding element of $S_n$. The conjugacy class of $\sigma_\varphi$, which corresponds to a partition of $n$, determines the number of cusps of $M_{\varphi}$, and how the meridian of the cusp projects to the fundamental group $\pi_1(S^1)$ of the base $S^1$ (i.e., how many times the cusp wraps around the mapping torus). For example, if $n=7$ and $\sigma_\varphi$ is in the conjugacy class corresponding to the partition $(4,2,1)$, then $M_{\varphi}$ will have $3$ cusps, with longitudes projecting to $4, 2$, and $1$ in $\pi_1(S^1)=\Z$. Clearly, $\sigma_\varphi$ will be in the alternating group $A_n\le S_n$ if and only if $l(\varphi)$ is even. Hence, if we want our random mapping classes to have a chance of hitting every conjugacy class, we will certainly need to ensure that both even and odd lengths $l(\varphi)$ occur with equal probability. Actually, we would like to have each conjugacy class $[\sigma]$ occurring with probability $|[\sigma]|/n!$, the proportion of elements of $S_n$ which are in that conjugacy class. In fact, experiments suggest that ensuring that even and odd length elements appear with equal frequency is sufficient to guarantee that the proportion of the sample that lies in a particular equivalence class $[\sigma]$ will be $\sim |[\sigma]|/n!$.

In practice, if we want (unreduced) words of length $L$, for each word we will choose the length to be a random choice between $L-3$ and $L+2$. So, for example, when we say the words are of length $l(\varphi)\sim 200$, they will be in the range $197\le l(\varphi)\le 202$. This choice of range is admittedly somewhat arbitrary---any range that ensures that even and odd length words occur with equal probability will, according to our experiments, result in the right distribution of conjugacy classes.

It follows from the construction of $\mc{T}$ that none of its edges are homotopic into a cusp of $M_{\varphi^\circ}$. Therefore, we can pull the edges of $\mc{T}$ tight to obtain a \textbf{geodesic triangulation}, i.e., we can homotope each edge to a geodesic representative. After pulling tight, each tetrahedron in $\mc{T}$ will be isometric to an ideal tetrahedron in $\h^3$. However, in some cases pulling a tetrahedron $\Delta$ tight will be a homotopy, but not an isotopy. That is, two edges of $\Delta$ may have to pass through each other to become geodesic (or one edge may pass through itself). In this case the orientation of the pulled tight $\Delta$ may be opposite to its original orientation, i.e., it will be negatively oriented in $M_{\varphi^\circ}$. Note that a geodesic triangulation with negatively oriented tetrahedra is not a triangulation in the usual sense: after pulling tight, the map $i:\bigsqcup_{\Delta\in\mc{T}} \Delta \to M_{\varphi^\circ}$ is no longer an embedding on the interiors of the tetrahedra.

Starting from the topological triangulation $\mc{T}$ that \texttt{flipper} outputs, we obtain the geodesic triangulation described above (possibly with negatively oriented tetrahedra) using the program \texttt{SnapPy} \cite{SnapPy}, by Culler, Dunfield, Goerner, and Weeks. \texttt{SnapPy} finds complex shape parameters $z\in \C$ for the tetrahedra in $\mc{T}$ (discussed further is \Cref{sec:shapes}), such that the tetrahedra glue up consistently. In particular, a tetrahedron with shape parameter $z$ will have edges with complex dihedral angles $\{z, \frac{1}{1-z},1-\frac{1}{z}\}$ (opposite edges have the same complex dihedral angle), and for each edge $e$ of $\mc{T}$, having incident complex dihedral angles $\{w_j\}_j$, we must have $\sum_j \log(w_j)=2\pi i$. This ensures that the total angle around each edge is $2\pi$, and that the metric is complete (see \cite{We} for a more thorough discussion of gluing equations). In general we may have $\mr{Im}(z)<0$ for some tetrahedron which thus will be negatively oriented in $M_{\varphi^\circ}$. If every tetrahedron in the geodesic realization of $\mc{T}$ is positively oriented, then we will say that the triangulation $\mc{T}$ is \textbf{geometric}.

In addition to the tetrahedra shapes, \texttt{SnapPy} also computes the volume $\mr{Vol}(\mc{T})$ of the triangulation. When \texttt{SnapPy} computes the volume, it does so by adding up the volume of all tetrahedra in the triangulation, taking the volume of a tetrahedron to be negative when it is negatively oriented. If the triangulation is geometric (i.e., no tetrahedra have negative volume), then in practice this is a good approximation of the volume of the manifold $M_{\varphi^\circ}$. In principle this volume calculation can be rigorously certified by \texttt{SnapPy}, using the HIKMOT \cite{hikmot} method or a descendant thereof, though we do not do this. If there are negatively oriented tetrahedra in $\mc{T}$, we ask \texttt{SnapPy} to retriangulate $M_{\varphi^\circ}$ to get a geometric triangulation, and if it succeeds it computes the volume using this new triangulation. In some cases this retriangulation will fail to produce a geometric triangulation. If this happens we compute the volume using the non-geometric veering triangulation. Out of $\sim$766,000 veering triangulations, we found that about 496,000 were non-geometric, and hence were retriangulated in an attempt to obtain a geometric triangulation from which to compute the volume. Of those that were retriangulated, about 51,000 failed to produce geometric triangulations. This might be concerning, but of the remaining 445,000, for which we \emph{did} find a geometric triangulation, the largest difference between the new volume (of the geometric triangulation) and the old volume (of the original veering triangulation) was on the order of $2.5\times 10^{-11}$, a negligible difference. Hence it is likely that the volumes of these 51,000 examples, for which we were unable to find geometric triangulations, are well approximated by the volume of their (non-geometric) veering triangulations.

For the experiments in \Cref{sec:systole}, we will also need \texttt{SnapPy} to compute the length of the systole of $M_{\varphi^\circ}$. For this it is necessary to compute the Dirichlet domain, which is computationally difficult and frequently fails for large cusped manifolds. One reason for this is that the Dirichlet domains of larger manifolds will typically have many very small faces, with vertices that are very close together. Since SnapPy uses numerical approximation, it becomes difficult to determine if two such vertices are indeed distinct, or if they should be considered the same vertex. For this reason, we are only able to compute the systole for relatively short words, of length $50\le l(\varphi)\le 100$. And despite this restriction, we still get about $12\%$ of these that either fail the Dirichlet domain computation with a runtime error, or terminate before completing due to compute cluster time limits. This situation seems unavoidable, since experience shows that computation time for the Dirichlet domain can vary wildly, even between two words of the same length.

\subsection{\texttt{flipper} Computation Time}
\begin{figure}
        \centering
        \begin{subfigure}[b]{0.49\textwidth}
                \includegraphics[scale=.36]{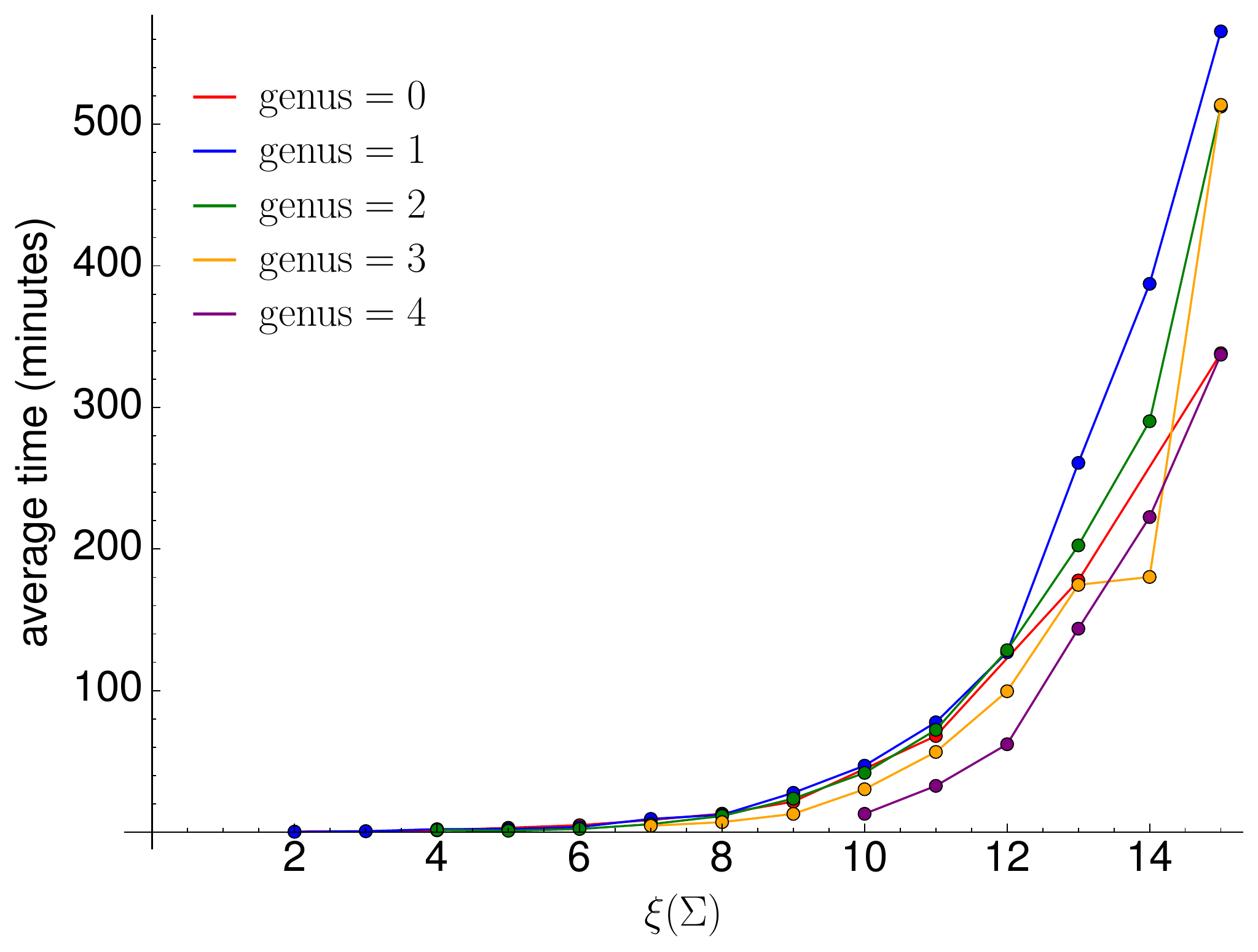}
        		\caption{}
        		\label{fig:cx_time}
       \end{subfigure}
        \begin{subfigure}[b]{0.5\textwidth}
                \includegraphics[scale=.36]{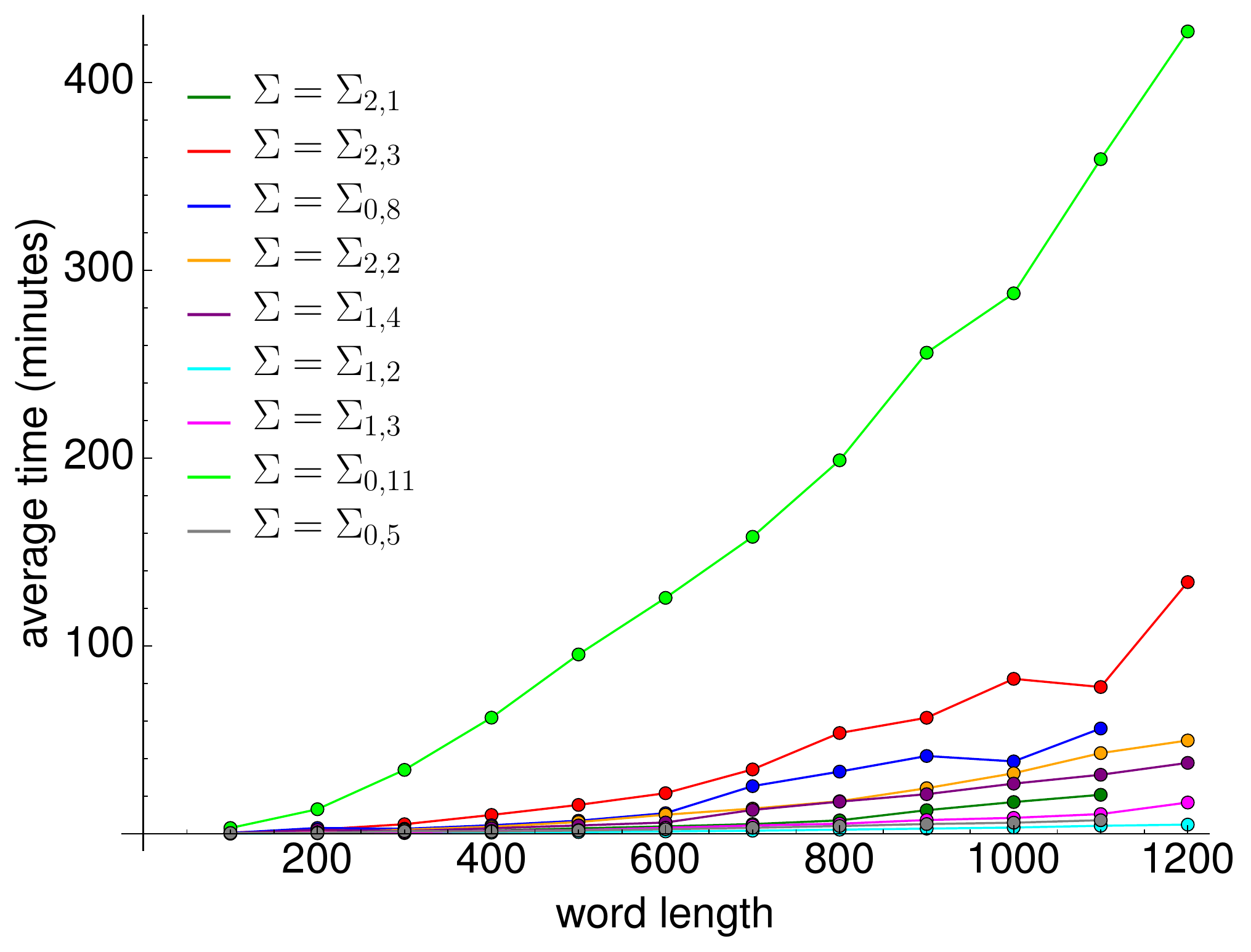}
                \caption{}
                \label{fig:wl_time}
        \end{subfigure}       
       \caption{Computation time as a function of complexity (left) and mapping class word length (right). Composing these with log gives plots (not shown) that suggest that computation time grows exponentially in complexity, and polynomially in word length.}
        \label{fig:time}
\end{figure}

In Figures \Cref{fig:cx_time} and \Cref{fig:wl_time} we plot computation time as a function of surface complexity $\xi(\Sigma_{g,n})=3g-3+n$ and mapping class word length, respectively (in the first plot mapping class length is fixed at $\sim 200$). The time shown is for both the \texttt{flipper} and \texttt{SnapPy} computations, plus some additional overhead computations. Compared to the \texttt{flipper} computations, however, the \texttt{SnapPy} and overhead calculations take a very small amount of time, so we should regard these plots as indications of approximate computation time for \texttt{flipper}. The experiments for which the \texttt{SnapPy} computations are very difficult, which were discussed in the preceding paragraph, are not included in these plots. We note that since our experiments were run on a multi-machine computing cluster, the hardware used for computation is not consistent across all experiment batches, and this may account for some of the noise we see in these plots.

\section{Genericity of Non-Geometric Veering Triangulations}
\label{sec:nongeo}

If $\Sigma$ is a once-punctured torus or a four-punctured sphere, and $\varphi:\Sigma\to \Sigma$ is a pseudo-Anosov mapping class, then the veering triangulation of $M_{\varphi^\circ}$ is geometric, as discussed in the introduction. For $\Sigma=\Sigma_{2,1}$, Hodgson--Issa--Segerman \cite{HoIsSe} have shown that of the $603$ distinct (up to conjugacy and inversion in $\mr{Mod}(\Sigma)$) pseudo-Anosov mapping classes contained in the ball of radius $7$ inside the Cayley graph of $\mr{Mod}(\Sigma)$ (using the generators in \Cref{fig:surf_gens}), $48$ are non-geometric. In general, when $\Sigma$ is not the 4-punctured sphere or the once-punctured torus, our experiments suggest that for very long mapping classes (in the generators given if Figures \Cref{fig:S0n_gens} and \Cref{fig:surf_gens}), geometric veering triangulations are rare.

\begin{figure}
        \centering
        \begin{subfigure}[b]{0.5\textwidth}
                \includegraphics[scale=.37]{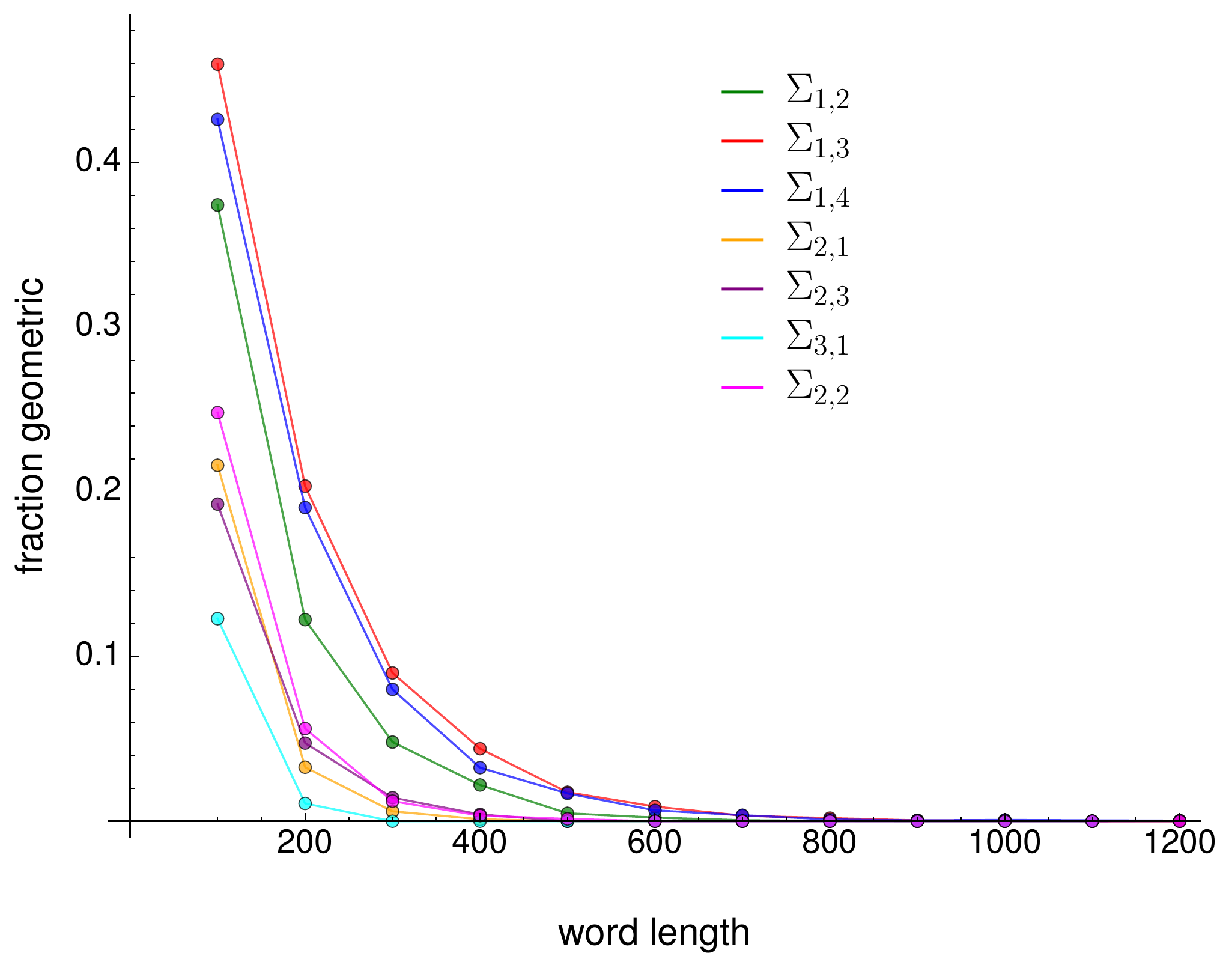}
        		\caption{}
        		\label{fig:fraction_geom_lin}
       \end{subfigure}
        \begin{subfigure}[b]{0.49\textwidth}
                \includegraphics[scale=.37]{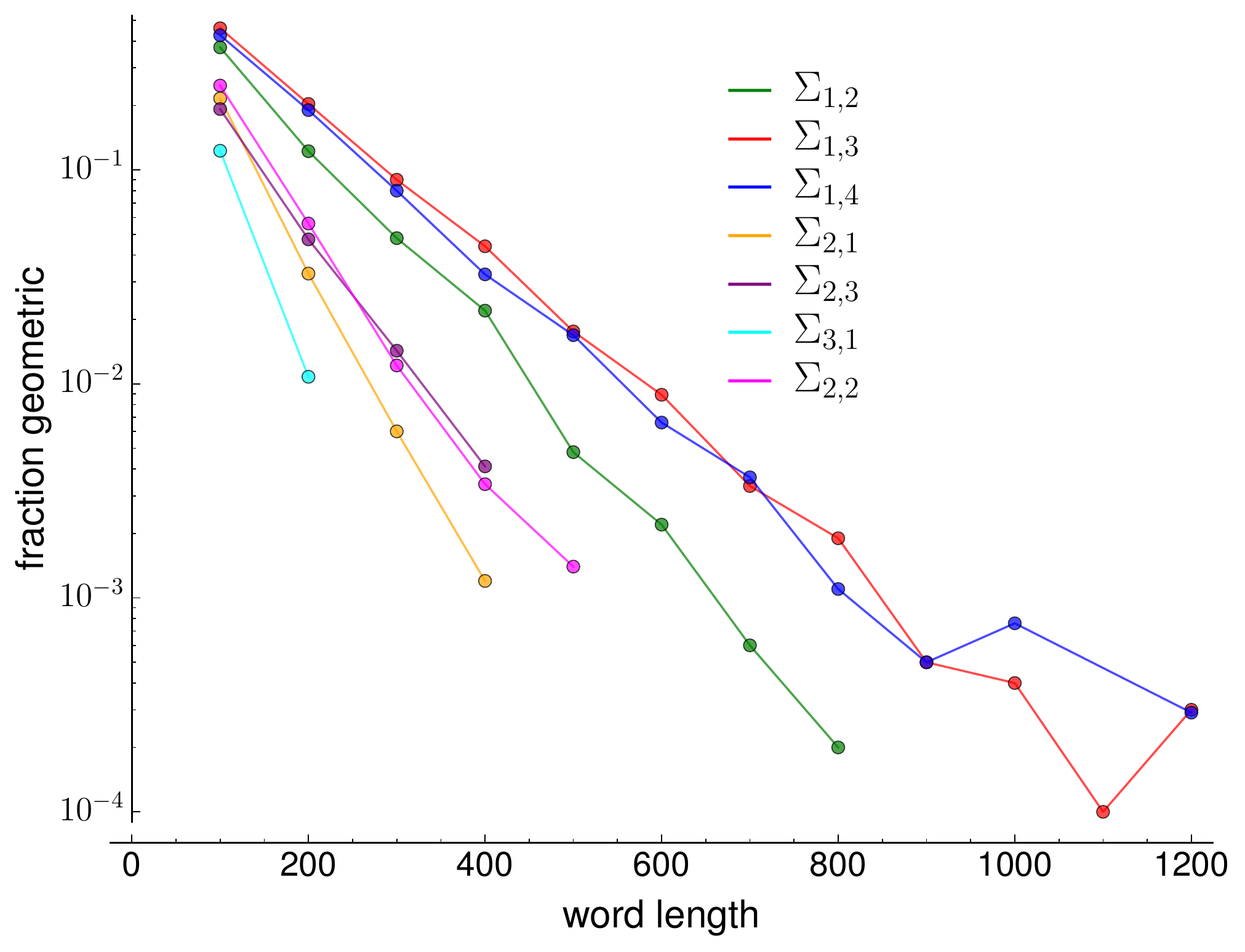}
                \caption{}
                \label{fig:fraction_geom_log}
        \end{subfigure}       
       \caption{Left: Each plot point represents between $\sim 2000$ and $\sim 20000$ triangulations with associated word length shown on the $x$-axis, and the proportion of these that are geometric is given by the $y$-coordinate. Right: The same data as on left, but with a log scale on the $y$-axis (plot point with $y=0$ are not shown as $\log(0)$ is undefined).}
        \label{fig:fraction_geom}
\end{figure}

\begin{figure}
        \centering
        \begin{subfigure}[b]{0.5\textwidth}
                \includegraphics[scale=.37]{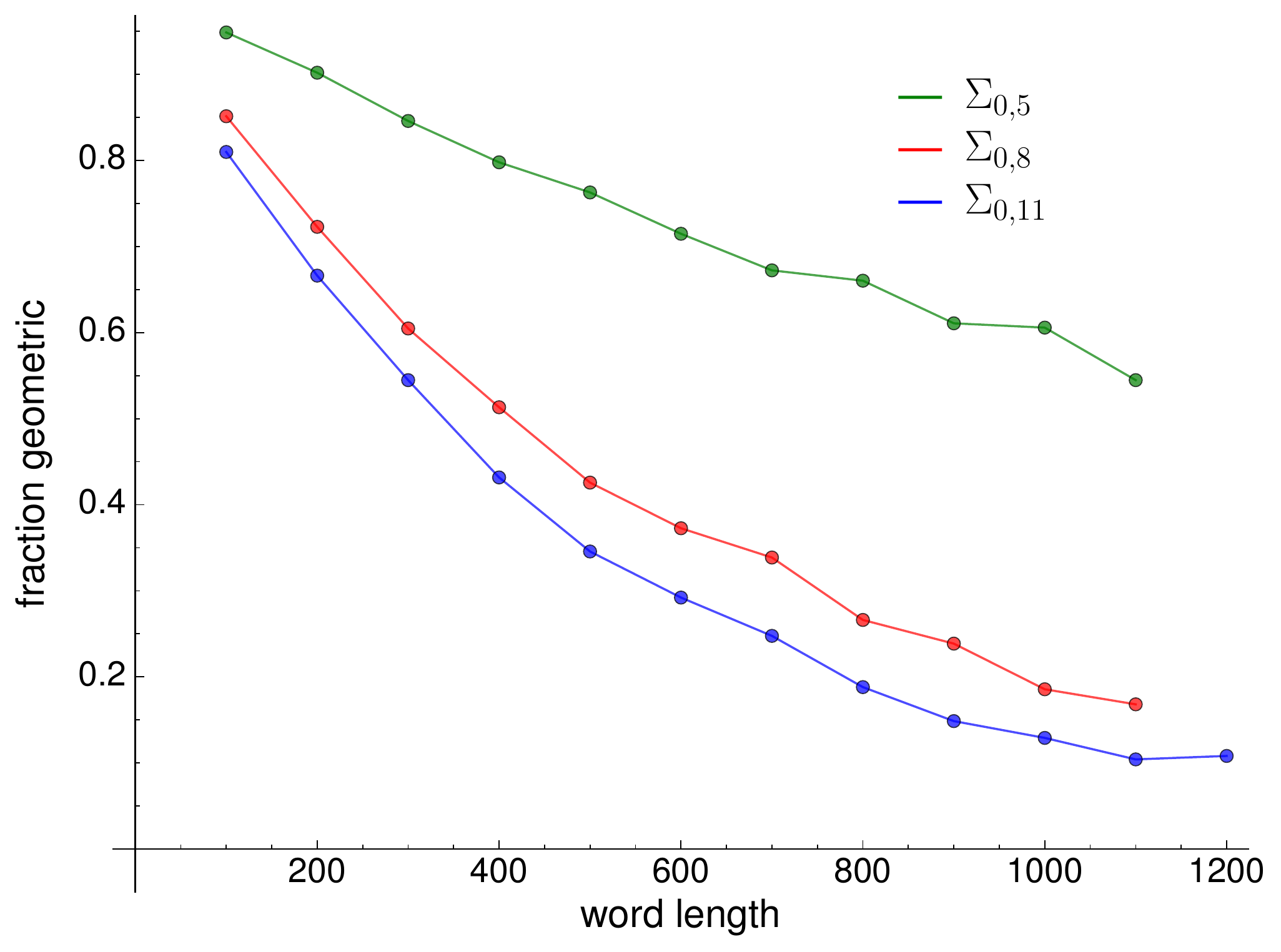}
        		\caption{}
        		\label{fig:fraction_geom_lin_0}
       \end{subfigure}
        \begin{subfigure}[b]{0.49\textwidth}
                \includegraphics[scale=.37]{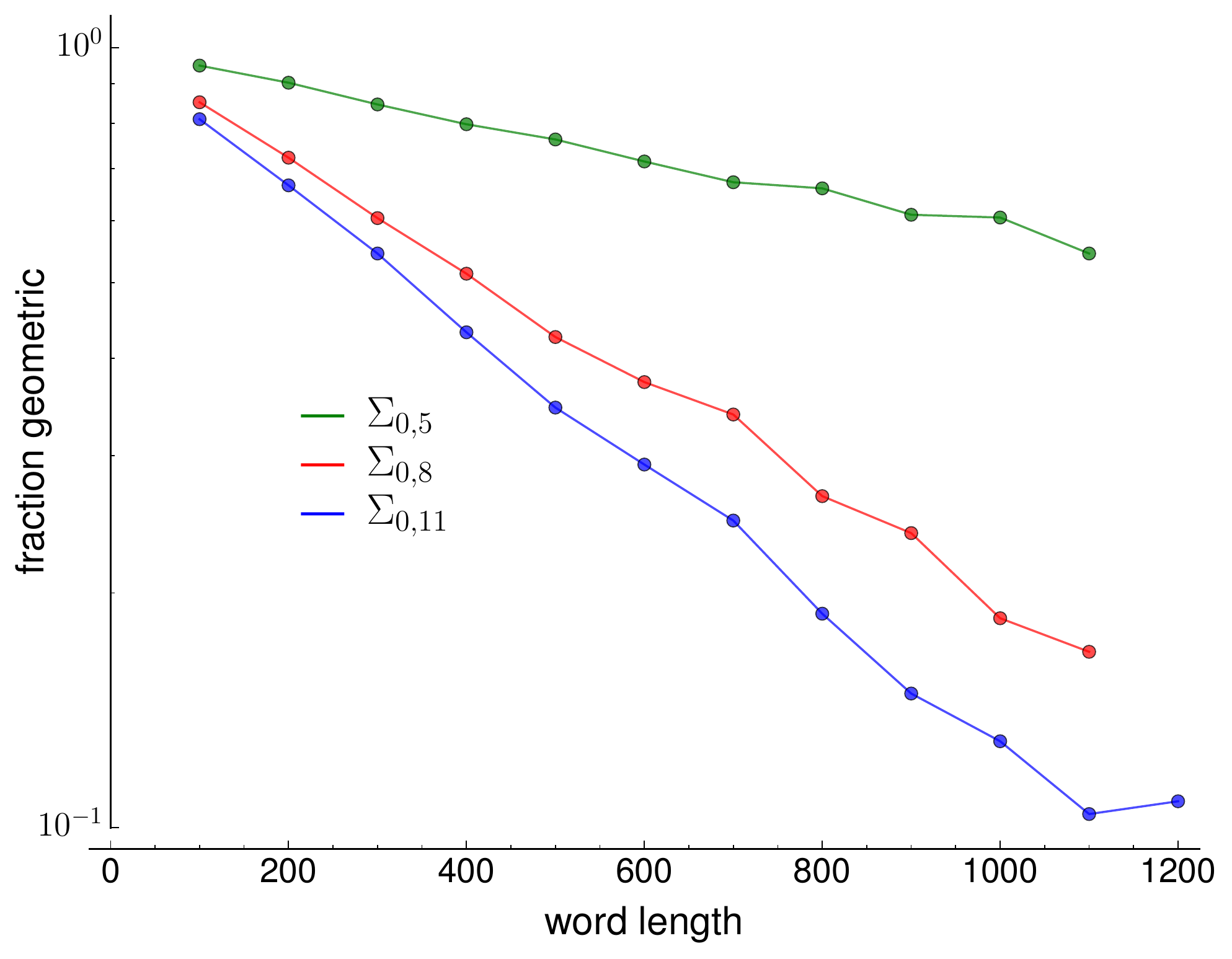}
                \caption{}
                \label{fig:fraction_geom_log_0}
        \end{subfigure}       
       \caption{Left: Each plot point represents $\sim 2000$ triangulations with associated word length shown on the $x$-axis, and the proportion of these that are geometric is given by the $y$-coordinate. Right: The same data as on left, but with a log scale on the $y$-axis.}
        \label{fig:fraction_geom_g0}
\end{figure}

\Cref{fig:fraction_geom_lin} shows, for several surfaces, the percentage of triangulations that are geometric, out of $\sim2000 - 20000$ examples for each word length shown. The plots in \Cref{fig:fraction_geom_log} reflect the same data as \Cref{fig:fraction_geom_lin} but with a log scale on the $y-$axis. The plots in \Cref{fig:fraction_geom_log} all end at their last non-zero value, since we cannot plot $y=0$ on a $y-$axis semi-log plot---For instance, for $\Sigma_{3,1}$ there are no geometric examples for word length $\ge 300$. One important consideration here is that for fractions less than $\sim 10^{-3}$, a sample size of $\sim 20000$ (for each word length) is likely not large enough, so we do not need to be concerned about the rather jagged tails of the plots for $\Sigma_{1,3}$ and $\Sigma_{1,4}$. For genus $0$ surfaces the decay is somewhat slower (see Figures \Cref{fig:fraction_geom_lin_0} and \Cref{fig:fraction_geom_log_0}), but nevertheless it still appears that for large enough words, geometric triangulations will be rare. These observations support the following conjecture:

\bg{conj}
Let $\Sigma=\Sigma_{g,n}$ be a surface of complexity $\xi(\Sigma)\ge 2$, and let $\mr{P}_{k,\Sigma}$ be the probability that a simple random walk of length $k$, on a set of generators for $\mr{Mod}(\Sigma)$, yields a pseudo-Anosov mapping class for which the layered veering triangulation of the mapping torus $M_{\varphi^\circ}$ is geometric. Then there exists a constant $c=c(\Sigma)$ such that 

$$
\lim_{k\to\infty}\frac{\mr{P}_{k,\Sigma}}{e^{-k\cdot c}}=0,
$$
i.e., $\mr{P}_{k,\Sigma}\ll e^{-kc}$ for $k\gg 0$.
\end{conj}

In addition to the above experimental evidence for this conjecture, consider the following partial heuristic. Suppose we have a word $\rho$ in the $\mr{Mod}(\Sigma)$ generators which has the following property: given any $\varphi\in \mr{Mod}(\Sigma)$ that is \textit{cyclically reduced}---i.e., all of its cyclic permutations are reduced, if $\varphi$ has $\rho$ as a subword, then the veering triangulation of $M_{\varphi^\circ}$ is non-geometric. In this case, the probability that a mapping class word of length $k$ contains $\rho$ as a subword approaches 1 exponentially as $k\to \infty$. This is easy to show---in short it is because the length of $\rho$ is fixed while the number of ways that $\rho$ can appear as a subword is growing without bound. Hence the existence of such a ``poison" word $\rho$ would imply that $P_{k,\Sigma}\to 0$ exponentially as $k\to\infty$. 

\begin{figure}
         \centering
         \includegraphics[scale=.75]{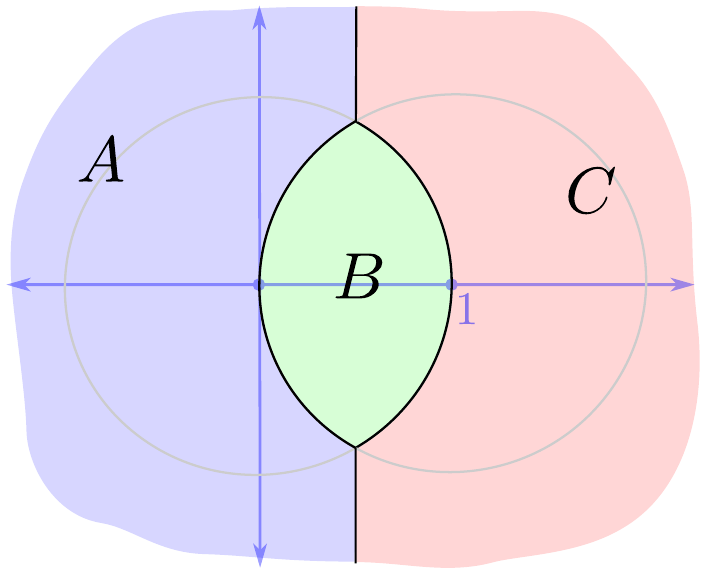}
        \caption{Fundamental domains for the action of $z\mapsto \frac{1}{1-z}$ on $\hat{\C}$. The fundamental domain $B$ parametrizes the shapes of ideal hyperbolic tetrahedra.}
        \label{fig:fund_doms}
\end{figure}

\section{Tetrahedra Shapes}
\label{sec:shapes}

Let $\Delta$ be a tetrahedron in the veering triangulation $\mathcal{T}$ of $M_{\varphi^\circ}$. Since edges of $\mathcal{T}$ are essential as edges on the fiber $\Sigma^\circ$, they are essential in $M_{\varphi^\circ}$, i.e., no edge of $\mathcal{T}$ can be homotoped into a cusp. This means that all the edges of $\mathcal{T}$ can be pulled tight to geodesics. So, although $\mathcal{T}$ may not be a geometric triangulation, we can pull it tight so that it is a geodesic triangulation, making $\Delta$ isometric---by an orientation preserving isometry---to a (possibly degenerate) ideal tetrahedron in $\h^3$. As noted at the beginning of \Cref{sec:bg}, $\Delta$ comes with a labelling of its vertices $\{v_0,v_1,v_2,v_3\}$, induced by the vertex labelling of the standard 3-simplex. The map $M_{(v_0,v_1,v_2)}(z)=\frac{(z-v_1)(v_0-v_2)}{(z-v_2)(v_0-v_1)}$ maps $v_0,v_1,v_2$  to $1,0,\infty$, respectively, and maps the final vertex $v_3$ to some $z\in\C$. The map $z\mapsto \frac{1}{1-z}$, which cyclically permutes $0,1$ and $\infty$, acts on $\hat{\C}$ with fundamental domains $A,B$ and $C$, as shown in \Cref{fig:fund_doms}. By composing with this map as necessary, we can guarantee that $M_{(v_0,v_1,v_2)}$ maps $v_3$ into $B$. Define the \textbf{shape parameter} $z_\Delta$ of $\Delta$ to be the complex number $z_\Delta := M_{(v_0,v_1,v_2)}(v_3)\in B$. Note that when $\Delta$ is negatively oriented in $M_{\varphi^\circ}$, we will have $\mr{Im}(z_\Delta)<0$, and when it is positively oriented, $\mr{Im}(z_\Delta)>0$. If $\mr{Im}(z_\Delta)=0$, then $\Delta$ is degenerate, i.e., it is flat.

\begin{figure}
        \centering
        \begin{subfigure}[b]{0.5\textwidth}
                \includegraphics[scale=.49]{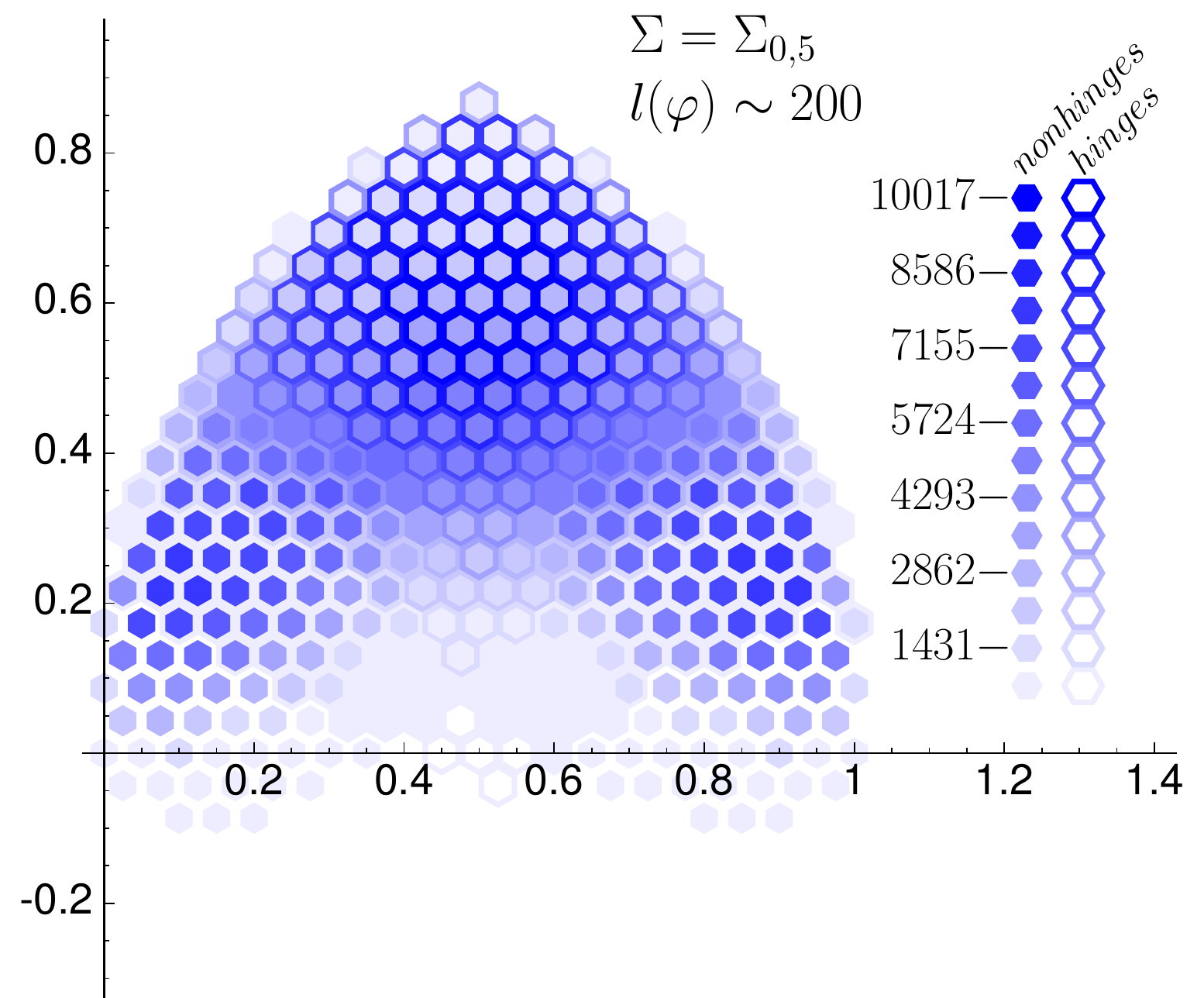}
        		\caption{}
        		\label{fig:S05_hist}
       \end{subfigure}
        \begin{subfigure}[b]{0.49\textwidth}
                \includegraphics[scale=.49]{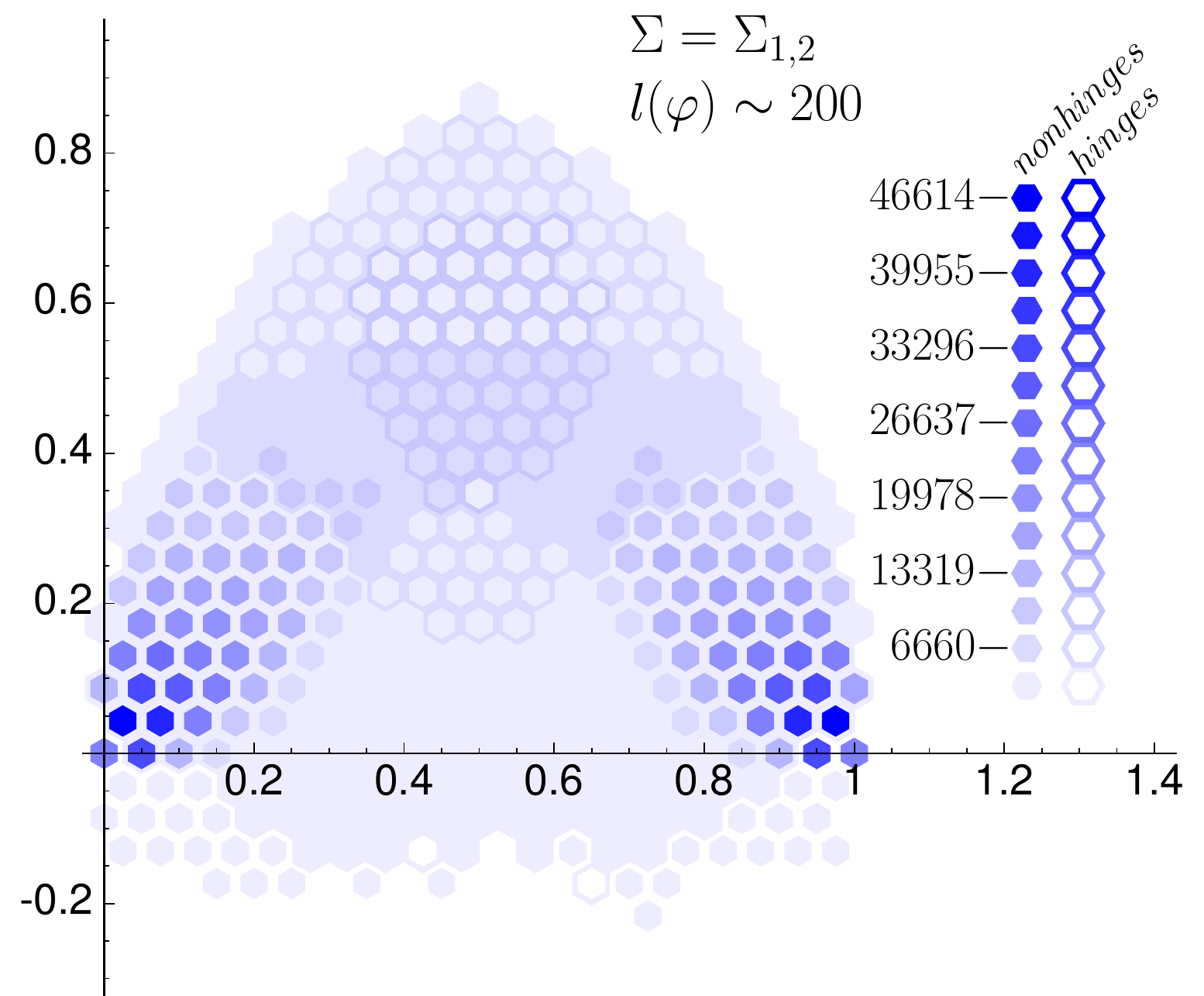}
                \caption{}
                \label{fig:S12_hist}
        \end{subfigure}       
       \caption{Tetrahedra shape histograms for veering triangulations with filled fiber $\Sigma_{0,5}$ (left) and $\Sigma_{1,2}$ (right). Each of these represents tetrahedra from 10000 triangulations, all having associated mapping classes of length $l(\varphi)\sim 200$. Each bin is separated into a central hexagon for non-hinges and a thick border for hinges, with opacity indicating the number of hinges/non-hinges, as indicated by the bar chart.}
        \label{fig:hists}
\end{figure}

\begin{figure}
       \begin{subfigure}[b]{0.244\textwidth}
        		\centering
                \includegraphics[scale=.275]{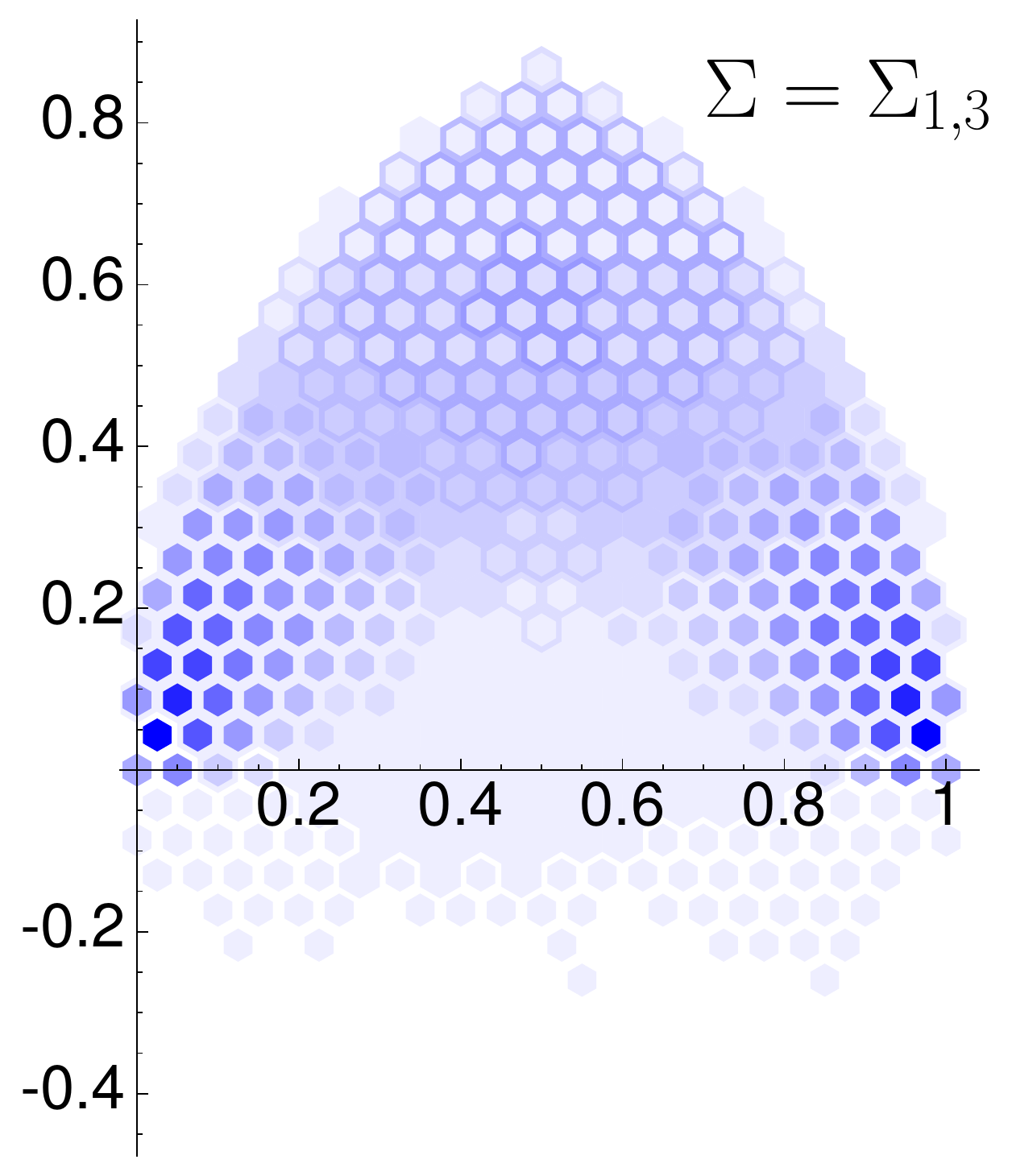}
        		\label{fig:ts_histo_13}
       \end{subfigure}      
       \begin{subfigure}[b]{0.244\textwidth}
        		\centering
                \includegraphics[scale=.275]{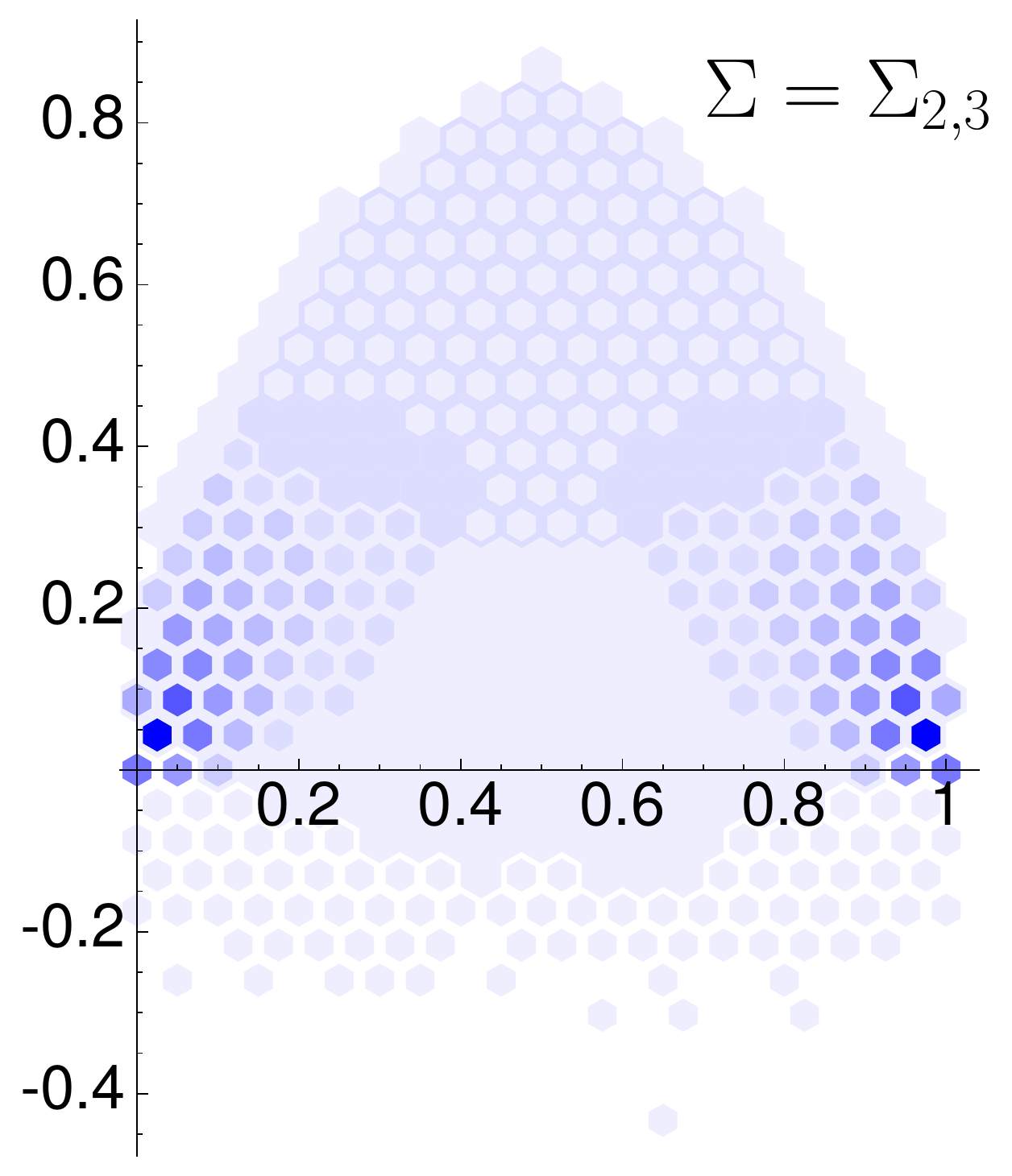}
        		\label{fig:ts_histo_23}
       \end{subfigure}       
       \begin{subfigure}[b]{0.244\textwidth}
        		\centering
                \includegraphics[scale=.275]{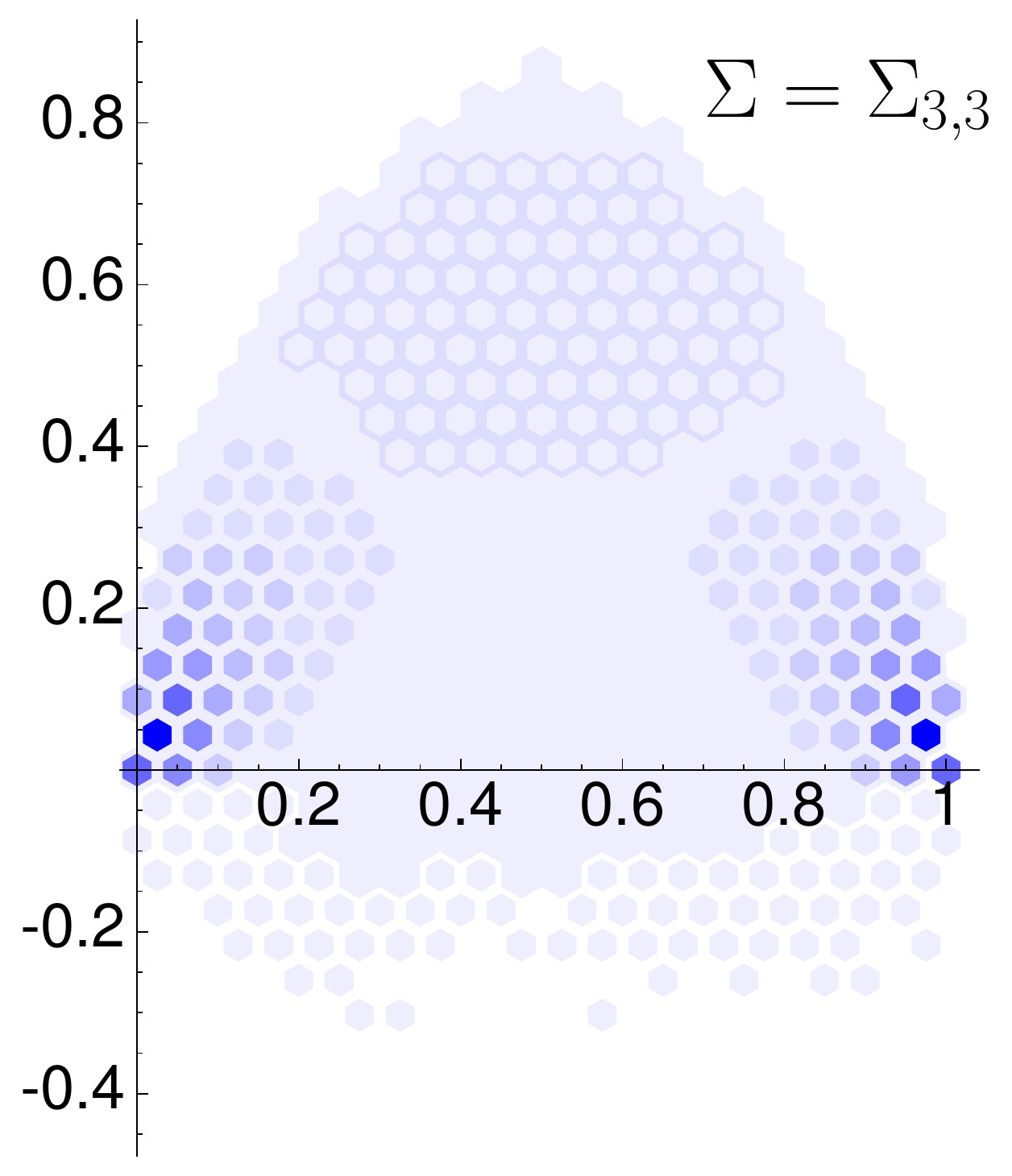}
        		\label{fig:ts_histo_33}
       \end{subfigure}
       \begin{subfigure}[b]{0.244\textwidth}
        		\centering
                \includegraphics[scale=.275]{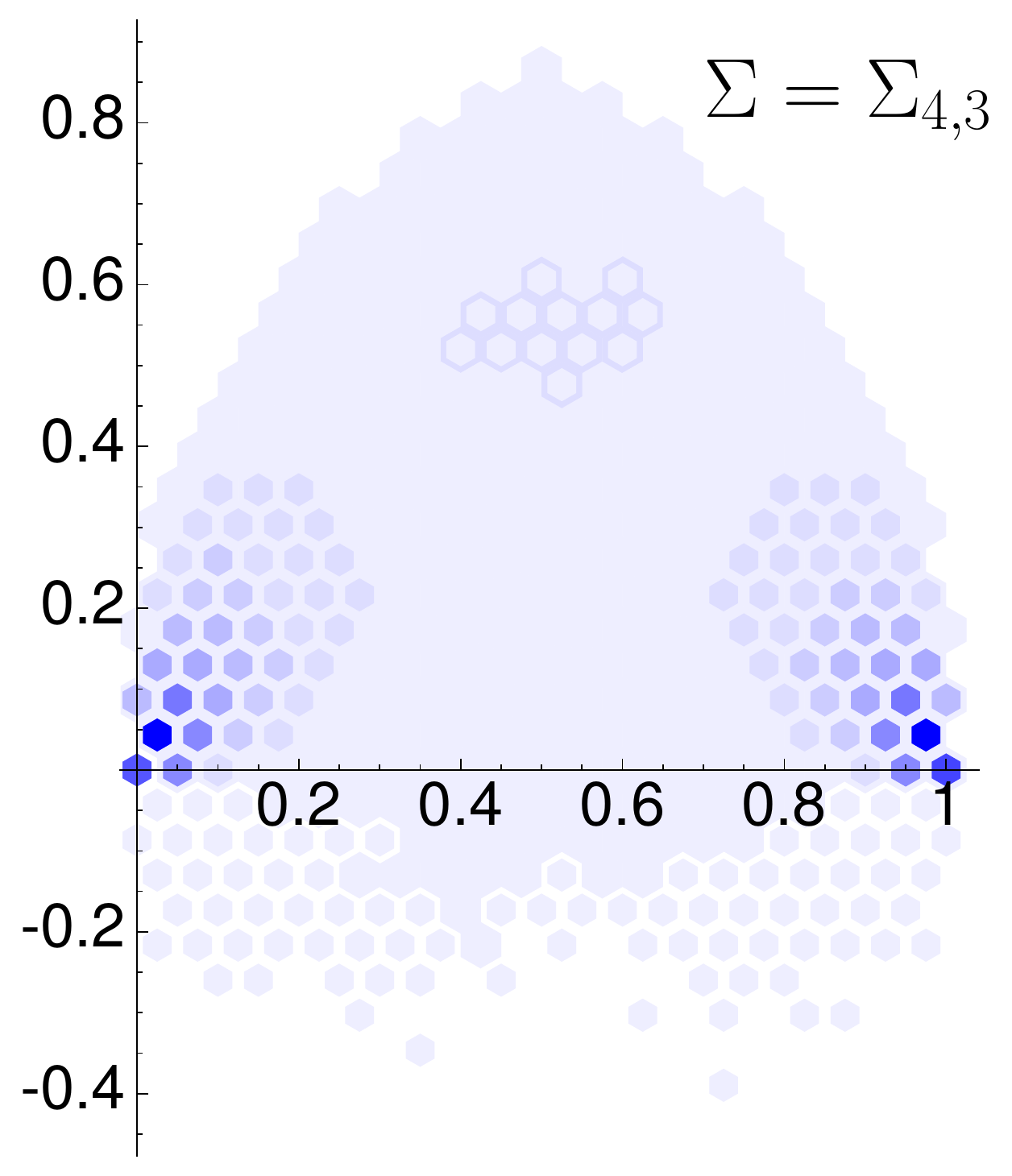}
        		\label{fig:ts_histo_43}
       \end{subfigure}
       \begin{subfigure}[b]{0.244\textwidth}
        		\centering
                \includegraphics[scale=.275]{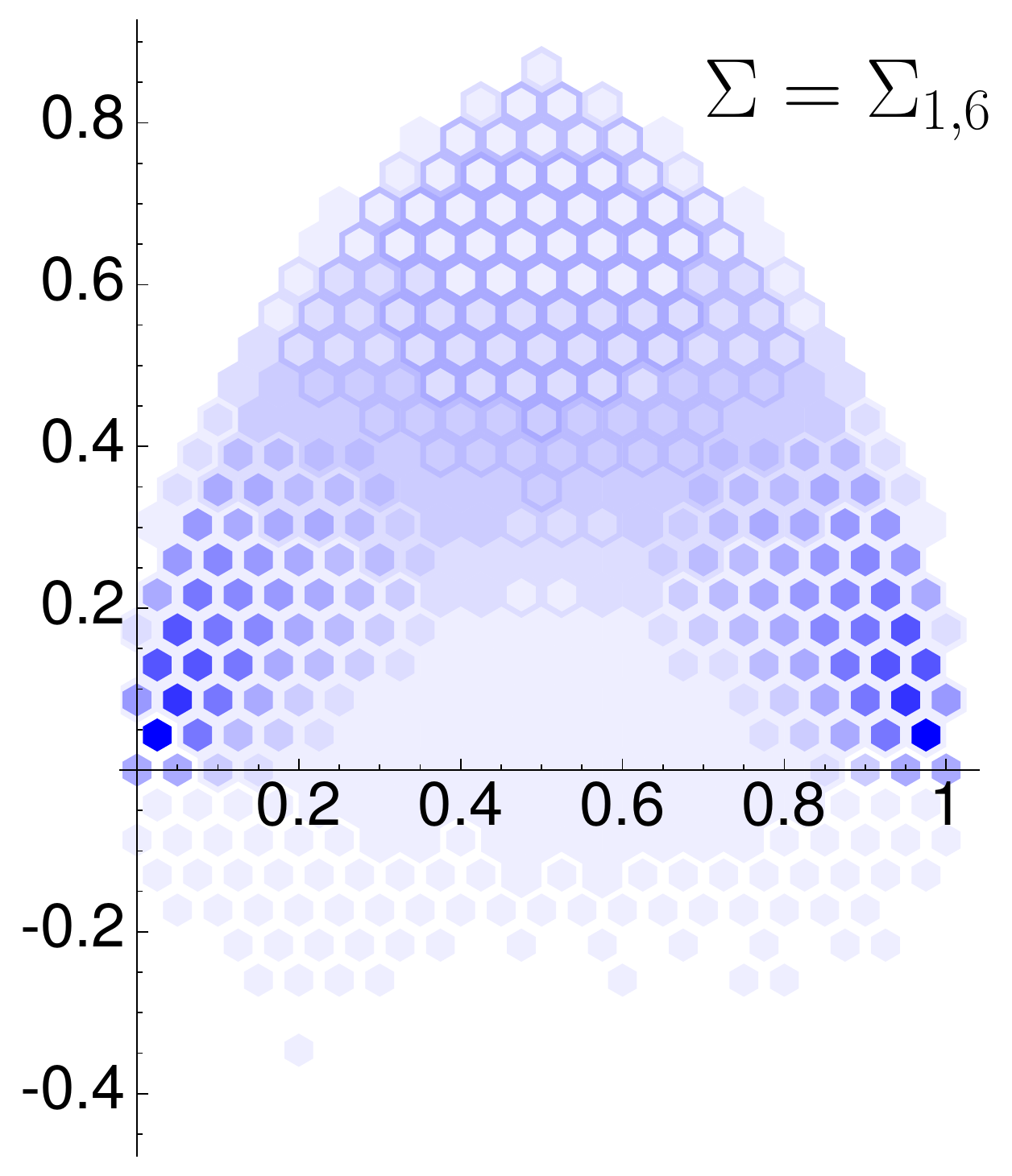}
        		\label{fig:ts_histo_16}
       \end{subfigure}
       \begin{subfigure}[b]{0.244\textwidth}
        		\centering
                \includegraphics[scale=.275]{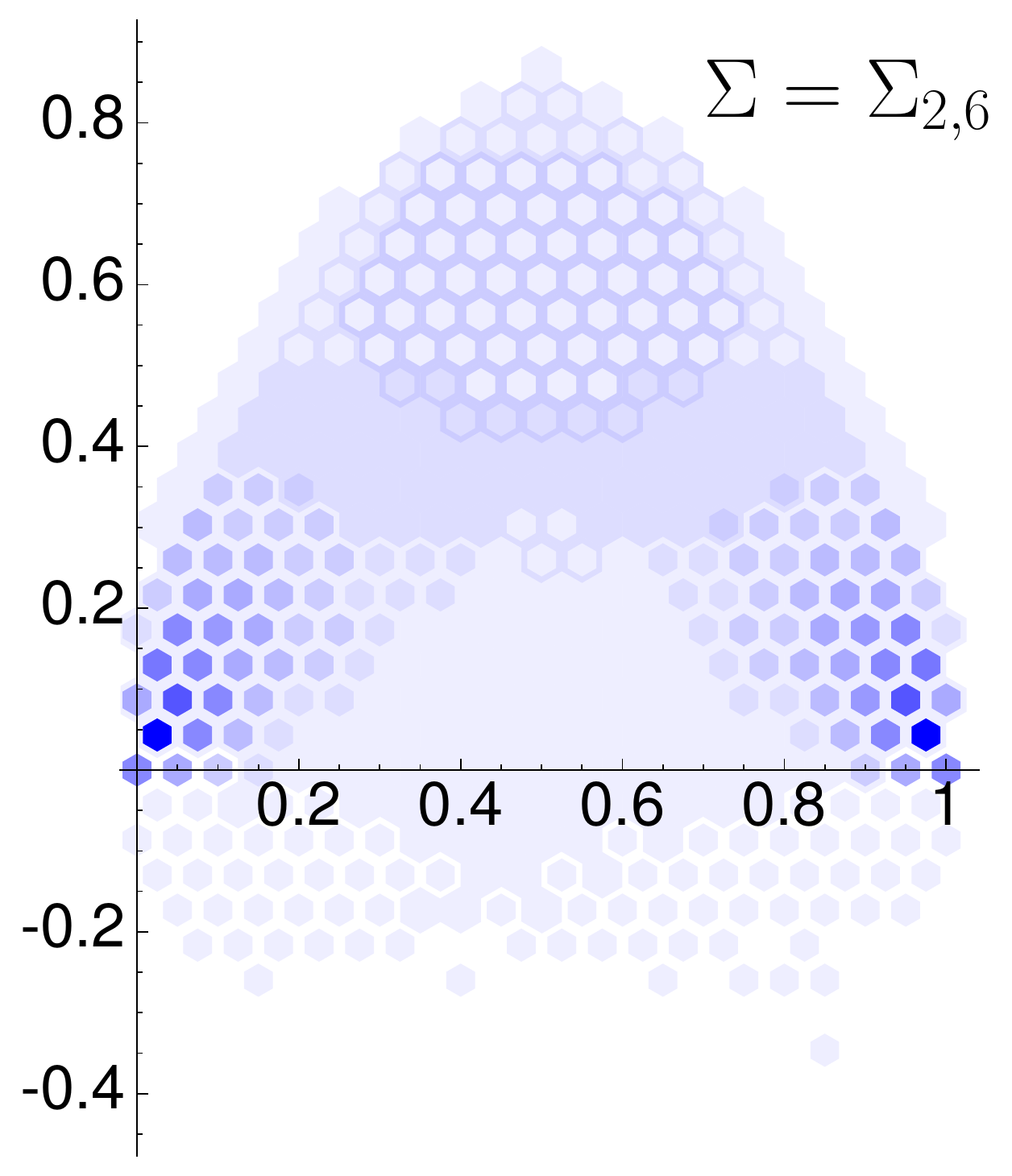}
        		\label{fig:ts_histo_26}
       \end{subfigure}
       \begin{subfigure}[b]{0.244\textwidth}
        		\centering
                \includegraphics[scale=.275]{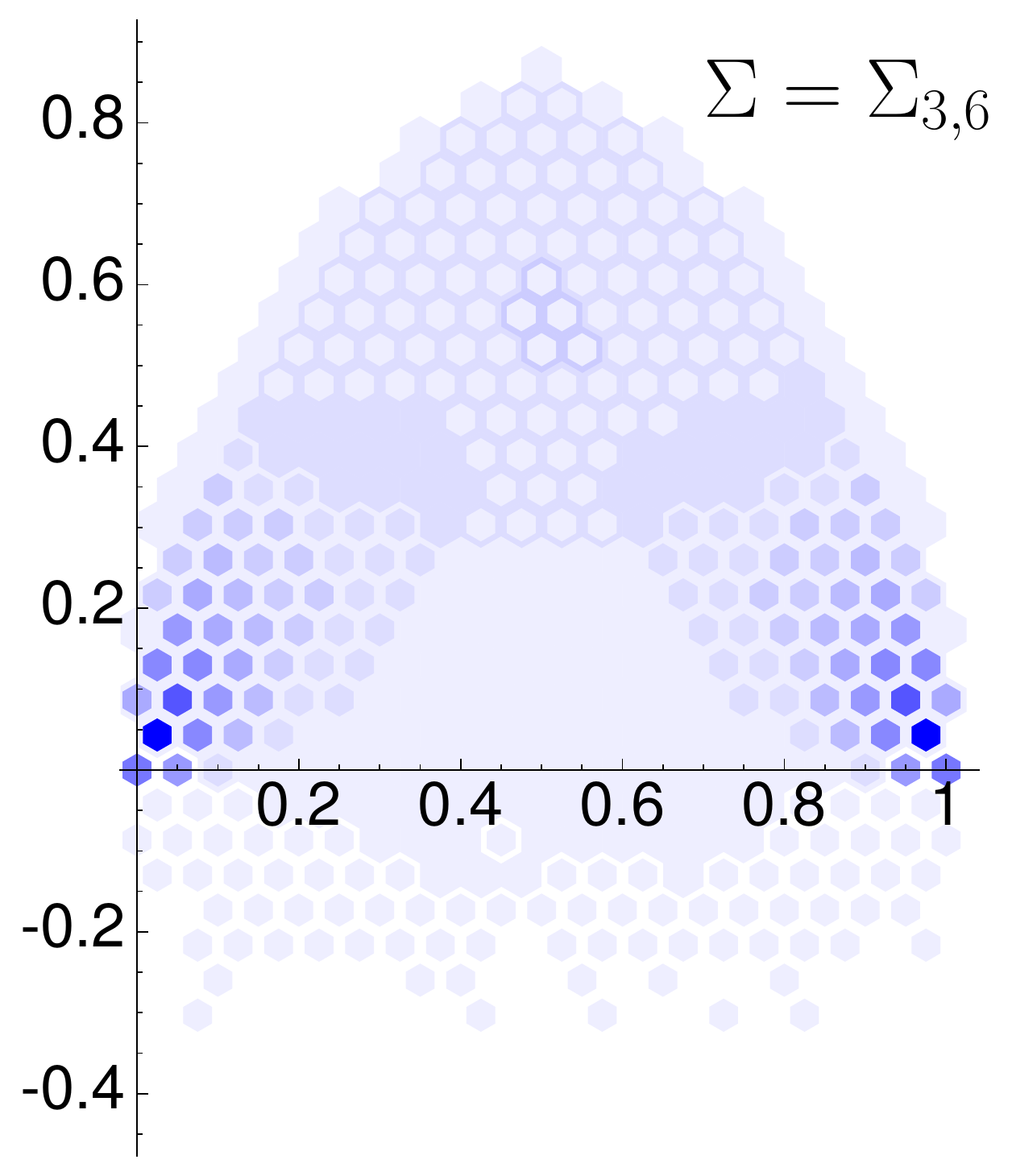}
        		\label{fig:ts_histo_36}
       \end{subfigure}
       \begin{subfigure}[b]{0.244\textwidth}
        		\centering
                \includegraphics[scale=.275]{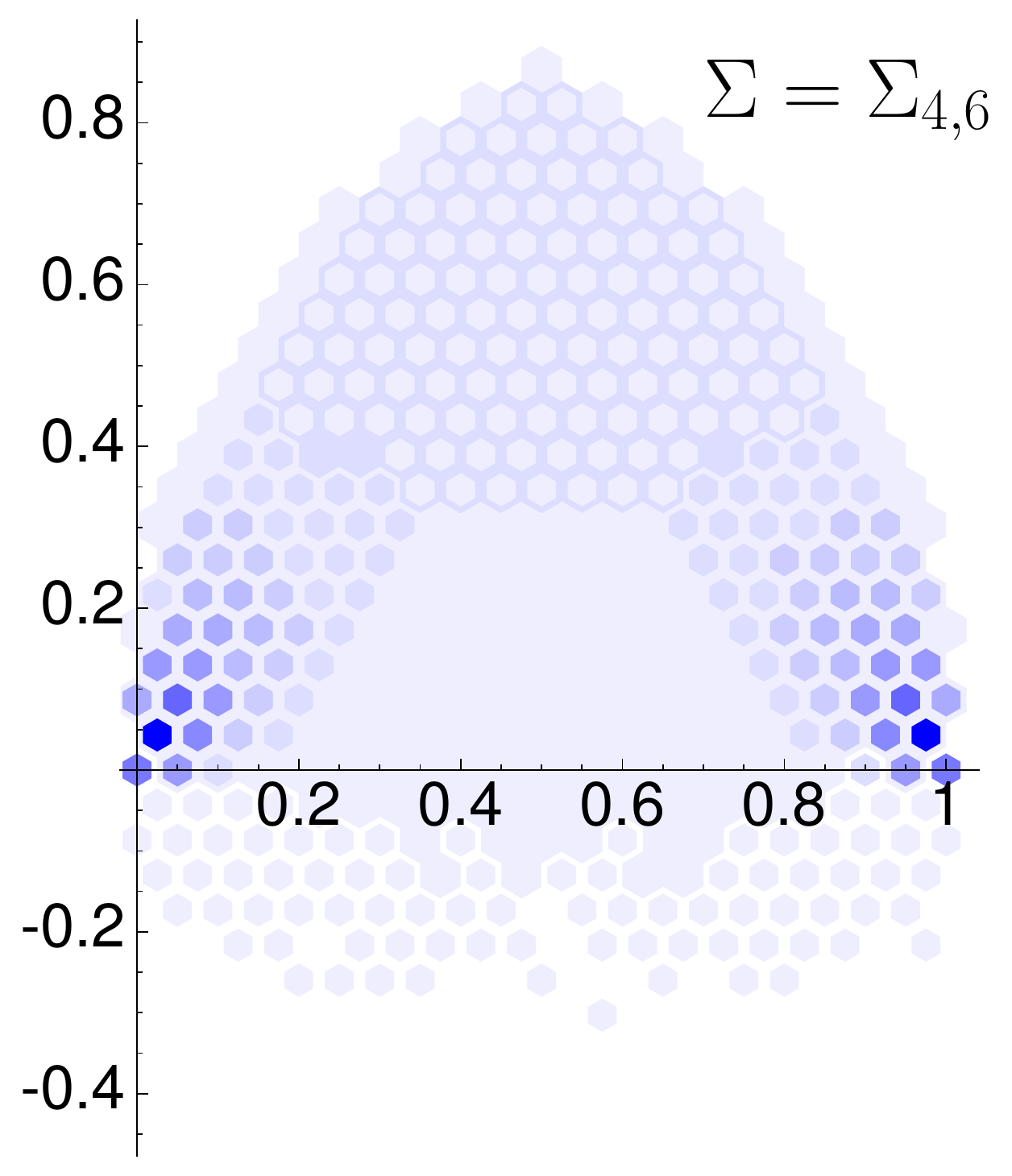}
        		\label{fig:ts_histo_46}
       \end{subfigure}
              \caption{Tetrahedra shape histograms with filled fiber as indicated. Each histogram is based on a sample of between $5000$ and $10000$ triangulations, all with associated mapping classes having word length $l(\varphi)\sim 200$.}
        \label{fig:hist_table}
\end{figure}

In this section we investigate the distribution in $B$ of tetrahedron shapes for veering triangulations. In particular, given a surface $\Sigma_{g,n}$ and a randomly generated pseudo-Anosov mapping class $\varphi$ of length $l(\varphi)$, generated as described above, we are interested in what shape parameters are realized by hinge and non-hinge tetrahedra in the veering triangulation of $M_{\varphi^\circ}$. 

Figures \Cref{fig:S05_hist} and \Cref{fig:S12_hist} shows a histogram of the tetrahedron shapes from a sample of 10000 veering triangulations with filled fiber $\Sigma_{0,5}$ and $\Sigma_{1,2}$, respectively, and word length $l(\varphi)\sim 200$. The opacity of a bin in this histogram reflects how full that bin is compared to the fullest bin, and hinges and non-hinges are shown separately, with the number of non-hinges in the bin indicated by the center hexagon, and the number of hinges by the hexagonal border. This allows one to easily see the overall fullness of the bin, as well as the relative abundance of hinges vs. non-hinges. Immediately, we see a sharp difference in the distribution of tetrahedra shapes for $\Sigma_{0,5}$ versus $\Sigma_{1,2}$. In particular, the $\Sigma_{0,5}$ sample has a much larger proportion of hinge tetrahedra, or \emph{hinge density}, than does the $\Sigma_{1,2}$ sample (about .513 versus .341), and the tetrahedra for $\Sigma_{1,2}$ are on average flatter (i.e., their shapes have small imaginary part). There are, however, notable similarities. For both, hinge shape parameters have, on average, greater imaginary part than non-hinges, and non-hinge shape parameters typically have real part away from $0.5$. We also note that there are relatively few negatively oriented tetrahedra (i.e., those whose shape parameters have negative imaginary part) in these histograms, compared to those that are positively oriented. On the other hand, \Cref{fig:fraction_geom} shows that for most filled fibers typical triangulations contain at least one negatively oriented tetrahedra.

As for shape parameter histograms for other filled fibers $\Sigma$, we find distributions very similar to those in \Cref{fig:hists}, with the exception of $\Sigma=\Sigma_{0,4}$ and $\Sigma=\Sigma_{1,1}$ (these two are special cases, as discussed in the introduction, and their shape parameter histograms are somewhat different). Some of these histograms are shown in \Cref{fig:hist_table}, with filled fibers as indicated. From each row of \Cref{fig:hist_table} we see that, if we fix the number of punctures of the filled fiber and increase the genus, non-hinges become more plentiful and on average flatter, and hinges become more scarce. The same is not generally true if we fix genus and increase the number of punctures, a fact which is reflected in the plot of \Cref{fig:cx_hd}, which shows the hinge density $\frac{|\mc{T}|_H}{|\mc{T}|}$ for various filled fiber surfaces, plotted against complexity. Note, though, that the plot does suggest that hinge density could be \emph{eventually} decreasing in number of punctures $n$, i.e., for $n>N$ with $N$ depending on the genus.

\begin{rem}
	From the histograms in Figures \Cref{fig:hists} and \Cref{fig:hist_table}, one might be tempted to conclude that hinges do not degenerate to 0 and 1 on the real axis. A closer look reveals that this conclusion is probably erroneous: in our data, we find hinges within a distance of $10^{-2}$ of 0 or 1, which is quite close (though it does not compare to the closest non-hinge, which is within about $10^{-6}$). This can be seen in the histogram for $\Sigma_{3,3}$, which shows at least one hinge in the bin centered at 0.
\end{rem}

\section{Volume}
\label{sec:volume}

\begin{figure}
        \centering
        \begin{subfigure}[b]{0.5\textwidth}
                \includegraphics[scale=.38]{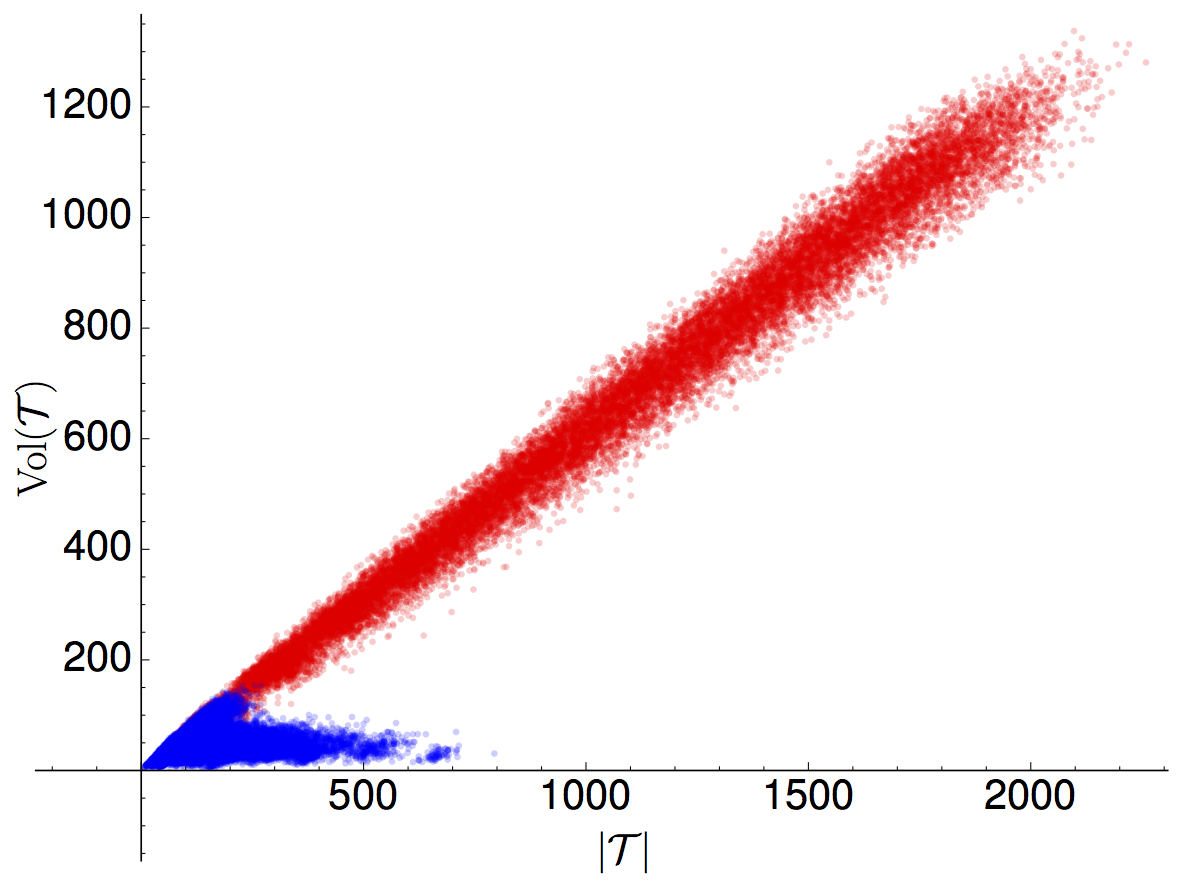}
        		\caption{}
        		\label{fig:linear_in_v1}
       \end{subfigure}
       \begin{subfigure}[b]{0.49\textwidth}
                \includegraphics[scale=.38]{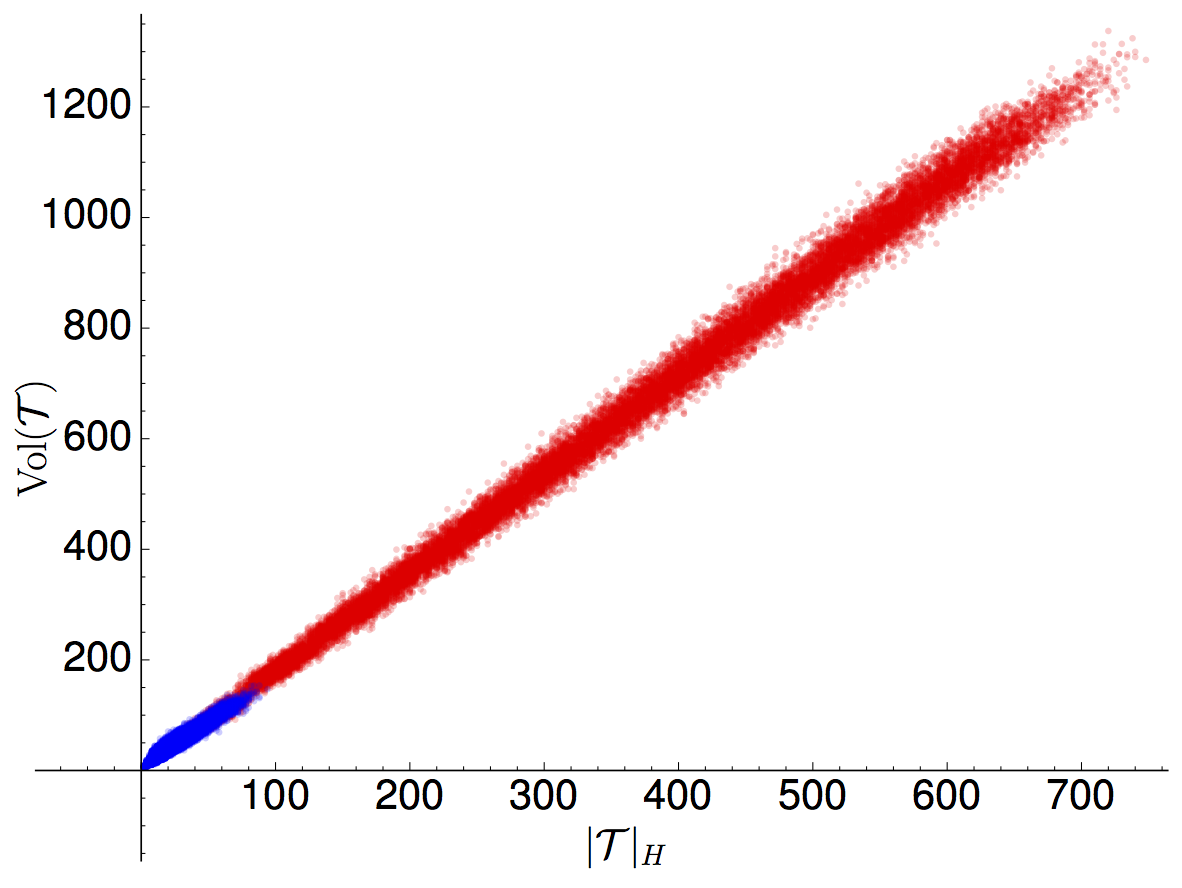}
        		\caption{}
        		\label{fig:linear_in_v2}
       \end{subfigure}       
       \caption{Left: Volume is roughly linear in the number of tetrahedra $|\mc{T}|$ in a veering triangulation. Red data points are triangulations coming from random walks on $\mr{Mod}(\Sigma)$; blue data points are triangulations coming from mapping classes having a subword consisting of a large power of a Dehn twist. Right: When volume is plotted against the number of hinge tetrahedra $|\mc{T}|_H$ the linear relation is considerably better, and even the non-random blue data points adhere to it. In both plots the sample is ~15000 random triangulations (red) and ~10000 non-random triangulations (blue).\ppp{9}}
        \label{fig:hingevol1}
\end{figure}

Combining work of Brock \cite{Br} and Maher \cite{Mah}, Rivin observed \cite[Theorem 5.3]{Ri} that for a random mapping class $\varphi_n$ of length $n$, one has $C_1n\le \mr{Vol}(M_{\varphi_n})\le C_2n$ with probability approaching $1$ as $n\to \infty$, for some constants $C_1,C_2$ depending on the fiber. Since the number of tetrahedra $|\mc{T}_n|$ in the veering triangulation associated to $\varphi_n$ is coarsely equal to $n$, it follows that the same statement holds with $|\mc{T}_n|$ replacing $n$ in the double inequality. Rivin also provides experimental results supporting his theorem. 
For the veering triangulations of mapping tori studied here we find, as predicted by Rivin's theorem,  a linear relation between $|\mc{T}|$ and volume---shown in the red scatter plot of \Cref{fig:linear_in_v1} for $\Sigma=\Sigma_{1,2}$. The relation is relatively coarse, though, and it is quite easy to construct arbitrarily bad outliers from the linear range. This is demonstrated by the blue data points in \Cref{fig:linear_in_v1}, which are triangulations coming from mapping classes having a subword consisting of a large power of a Dehn twist (n.b. these are not generic). On the other hand, if we instead look at volume as a function of the number of \emph{hinge} tetrahedra, which we will denote $|\mathcal{T}|_H$, the result is quite remarkable: not only are the random triangulations (red data points) more tightly grouped along the line, but the non-generic blue data points are no longer outliers (see \Cref{fig:linear_in_v2}).

\begin{figure}
        \centering
        \begin{subfigure}[b]{0.5\textwidth}
                \includegraphics[scale=.37]{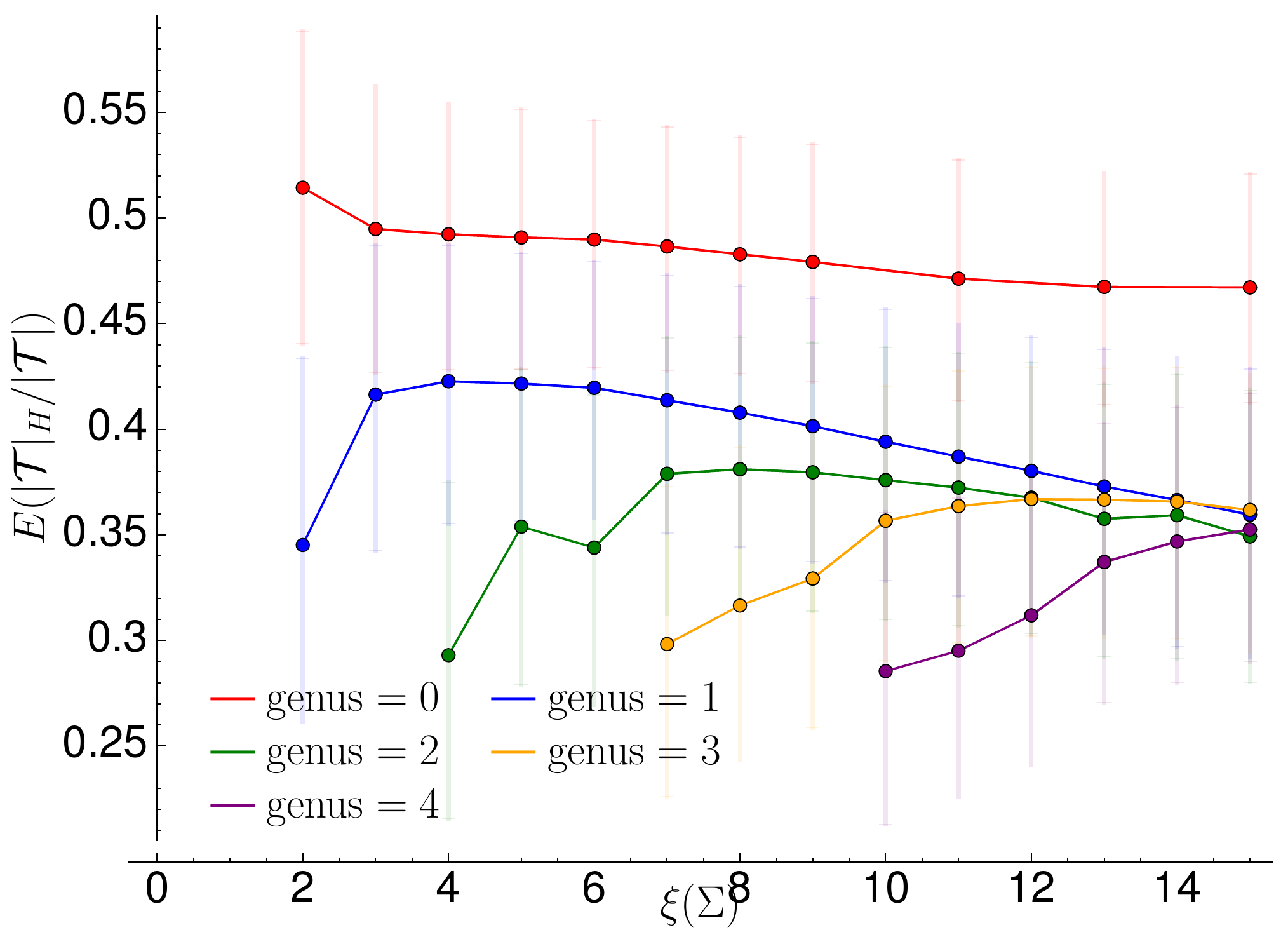}
        		\caption{}
        		\label{fig:cx_hd}
       \end{subfigure}
        \begin{subfigure}[b]{0.49\textwidth}
                \includegraphics[scale=.37]{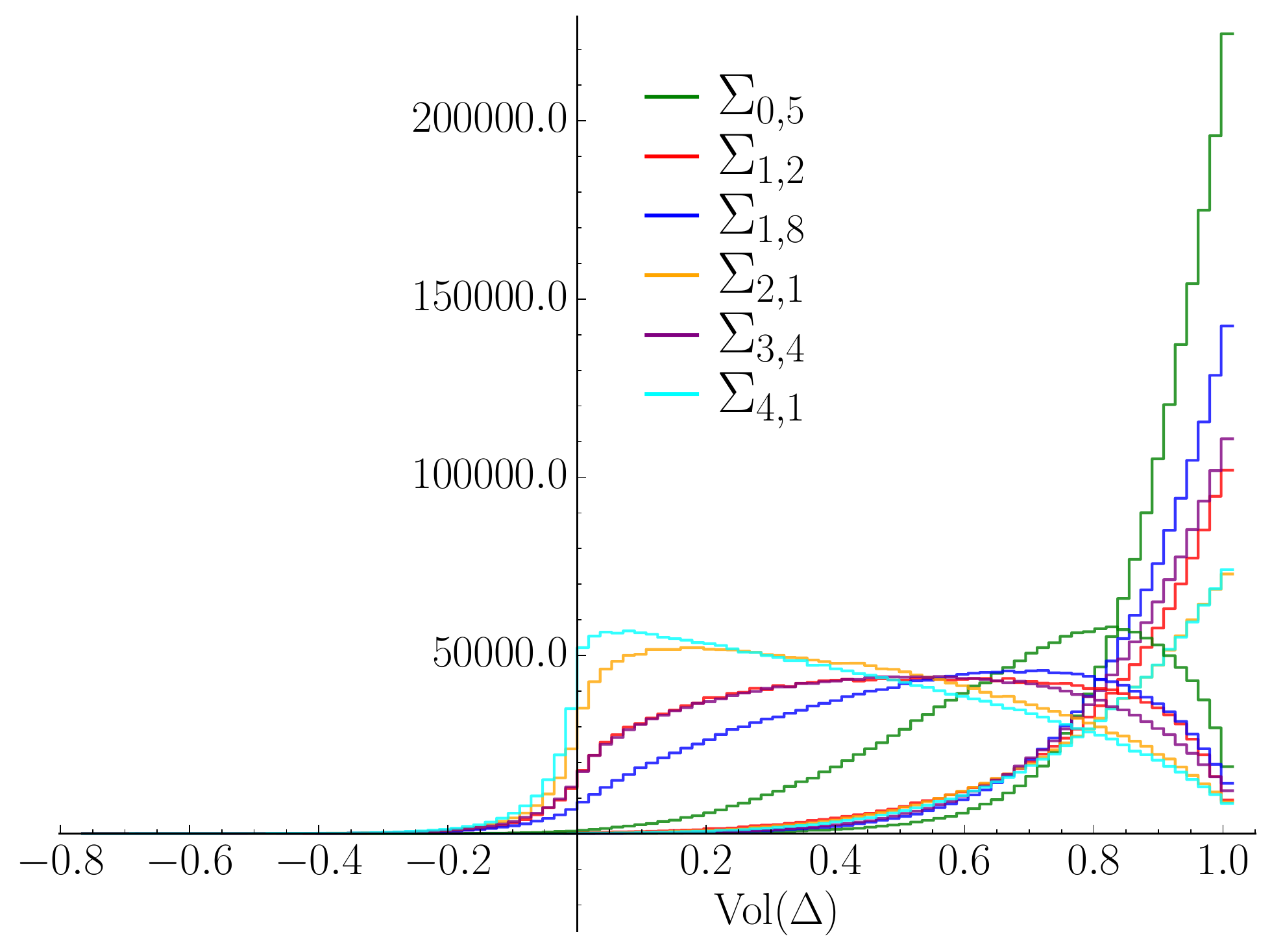}
                \caption{}
                \label{fig:tet_vol_hist}
        \end{subfigure}       
       \caption{Left: Average hinge density as a function of filled fiber complexity. Right: Histogram of tetrahedron volume for tetrahedra $\Delta$ coming from veering triangulations with associated surfaces as indicated. In order to show multiple histogram in one graphic, we only show the step graph determined by the height of the bins.}
        \label{fig:tv_hd}
\end{figure}

\placetextbox{0.86}{.89}{\rotatebox{90}{$\underbrace{\qquad\qquad\qquad\qquad\,\,\,}_{\text{\tiny hinges}}$}}
\placetextbox{0.86}{.76}{\rotatebox{90}{$\underbrace{\quad}_{\text{\tiny non-hinges}}$}}

\begin{figure}[h]
        \centering
       \begin{subfigure}[b]{0.49\textwidth}
                \includegraphics[scale=.38]{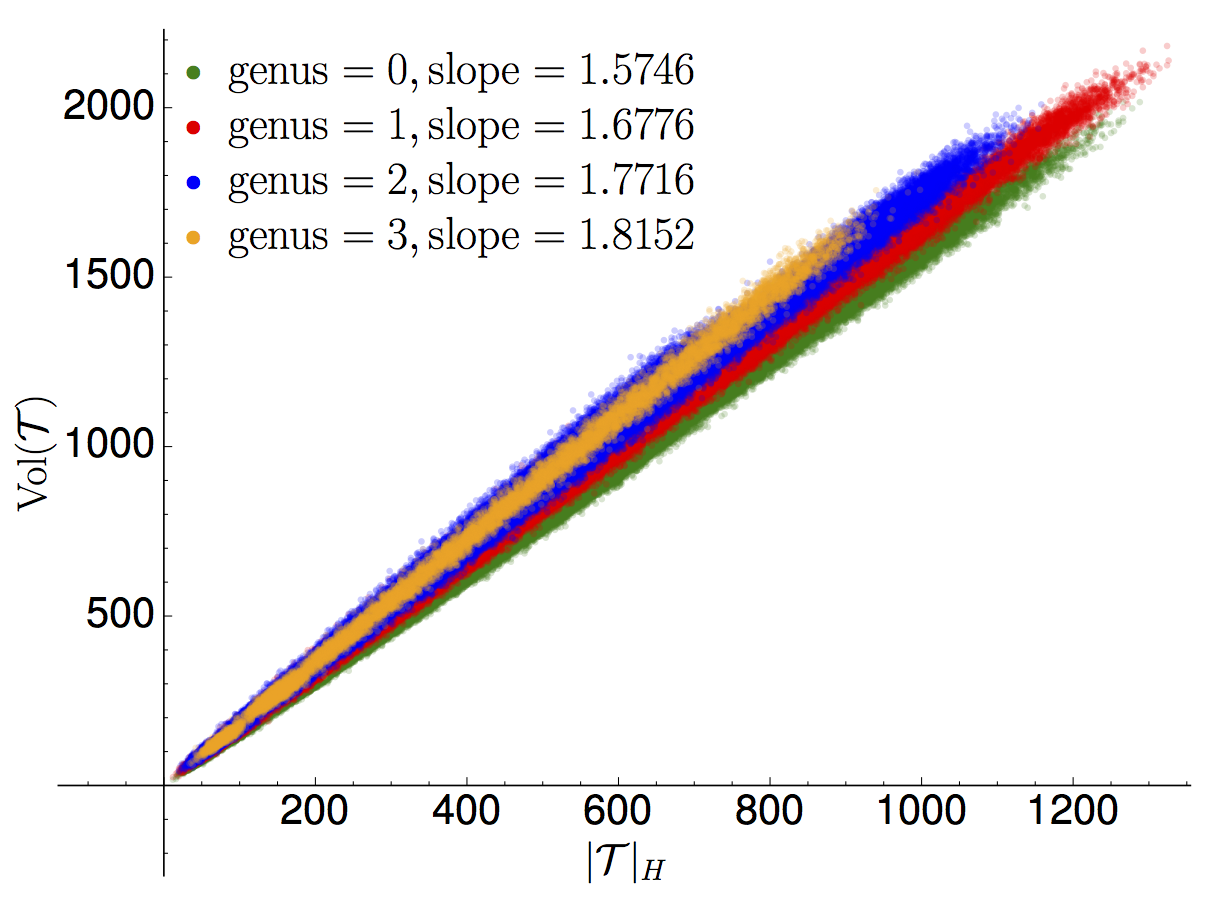}
                \caption{}
                \label{fig:hingevol_genus}
        \end{subfigure}
        \begin{subfigure}[b]{0.49\textwidth}
                \includegraphics[scale=.38]{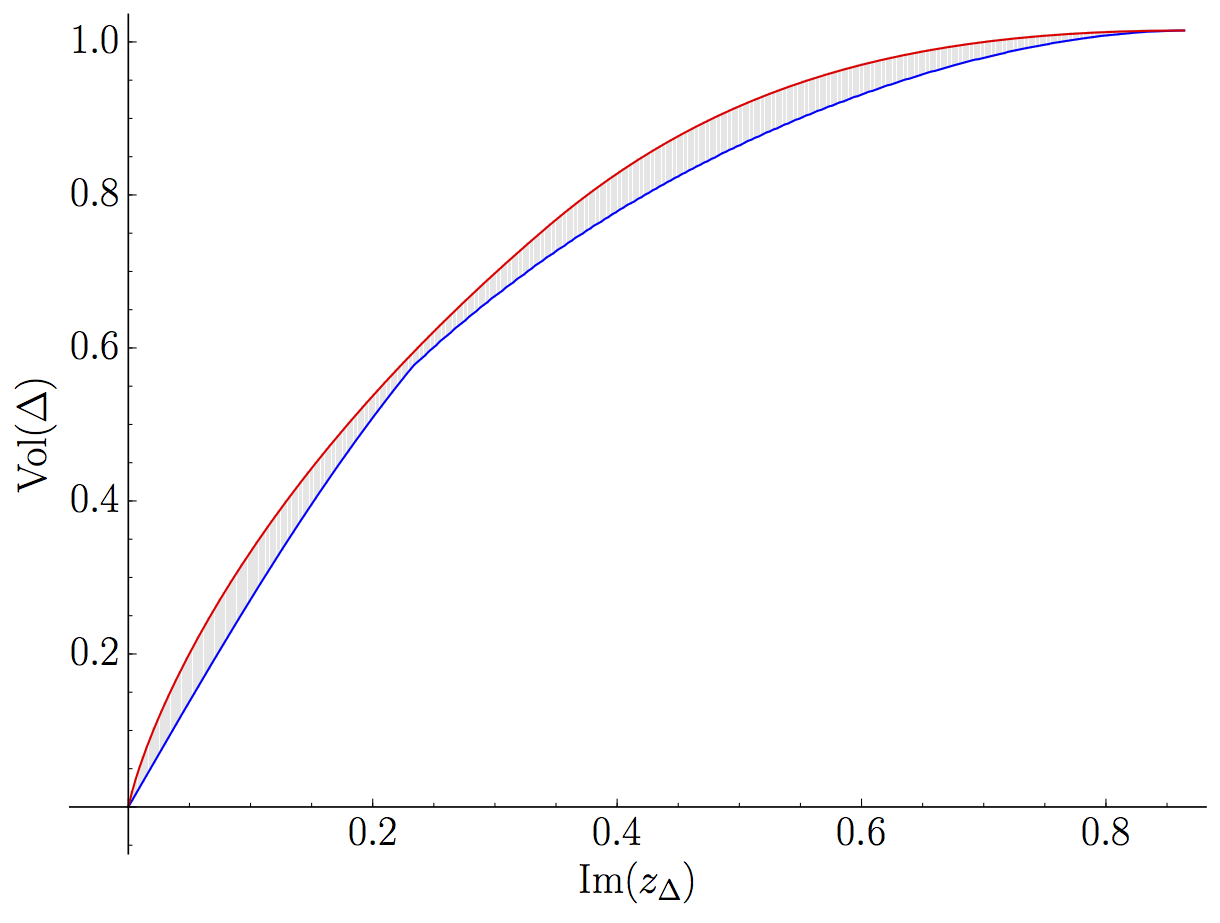}
                \caption{}
                \label{fig:loba2d}
        \end{subfigure}       
       \caption{Left: Scatter plots of $\mr{Vol}(\mc{T})$ against $|\mc{T}|_H$, for triangulations with filled fibers of various genera. Right: The red and blue plots show, respectively, the largest and smallest values that the $\mr{Vol}$ function takes for the given value of $\mr{Im}(z_\Delta)$.}
        \label{fig:hingevol_loba}
\end{figure}

To gain some insight into the above observations, consider first the results of the previous section, where we found that shape parameters for hinge and non-hinge tetrahedra are distributed quite differently in $B$. In particular, non-hinge tetrahedra are, on average, thinner than hinge tetrahedra (i.e., their shape parameters have smaller imaginary part). Perhaps unsurprisingly, these observations about tetrahedra shape reveal something about the expected volumes of hinge versus non-hinge tetrahedra. \Cref{fig:loba2d} shows a projection onto the $(\mr{Vol},\mr{Im}(z))$-plane of a plot of the volume function $\mr{Vol}:\{z\in B:\mr{Im}(z)\ge 0\}\to \R_{\ge 0}$. In particular, the red plot is $f_M(x)=\max_{\mr{Re}(z)}\{\mr{Vol}(z):\mr{Im}(z)=x\}$, and the blue plot is $f_m(x)=\min_{\mr{Re}(z)}\{\mr{Vol}(z):\mr{Im}(z)=x\}$. We can see from this that $\mr{Vol}(z_\Delta)$ is only subtly dependent on the real part $\mr{Re}(z_\Delta)$ of a tetrahedron's shape parameter, and is approximately determined by the imaginary part. Therefore, although there are in many cases far more non-hinge tetrahedra in a given triangulation, the majority of the volume appears to be concentrated in hinge tetrahedra. This is made somewhat clearer by plotting a histogram of tetrahedra volumes, with hinges and non-hinges plotted separately, as in \Cref{fig:tet_vol_hist}. Here again we see that hinges, on average, have very large volume compared to non-hinges. Thus it makes sense that the number of hinges $|\mathcal{T}|_H$ is a better predictor of volume than the total number of tetrahedra.

\begin{figure}
        \centering
        \begin{subfigure}[b]{0.5\textwidth}
                \includegraphics[scale=.38]{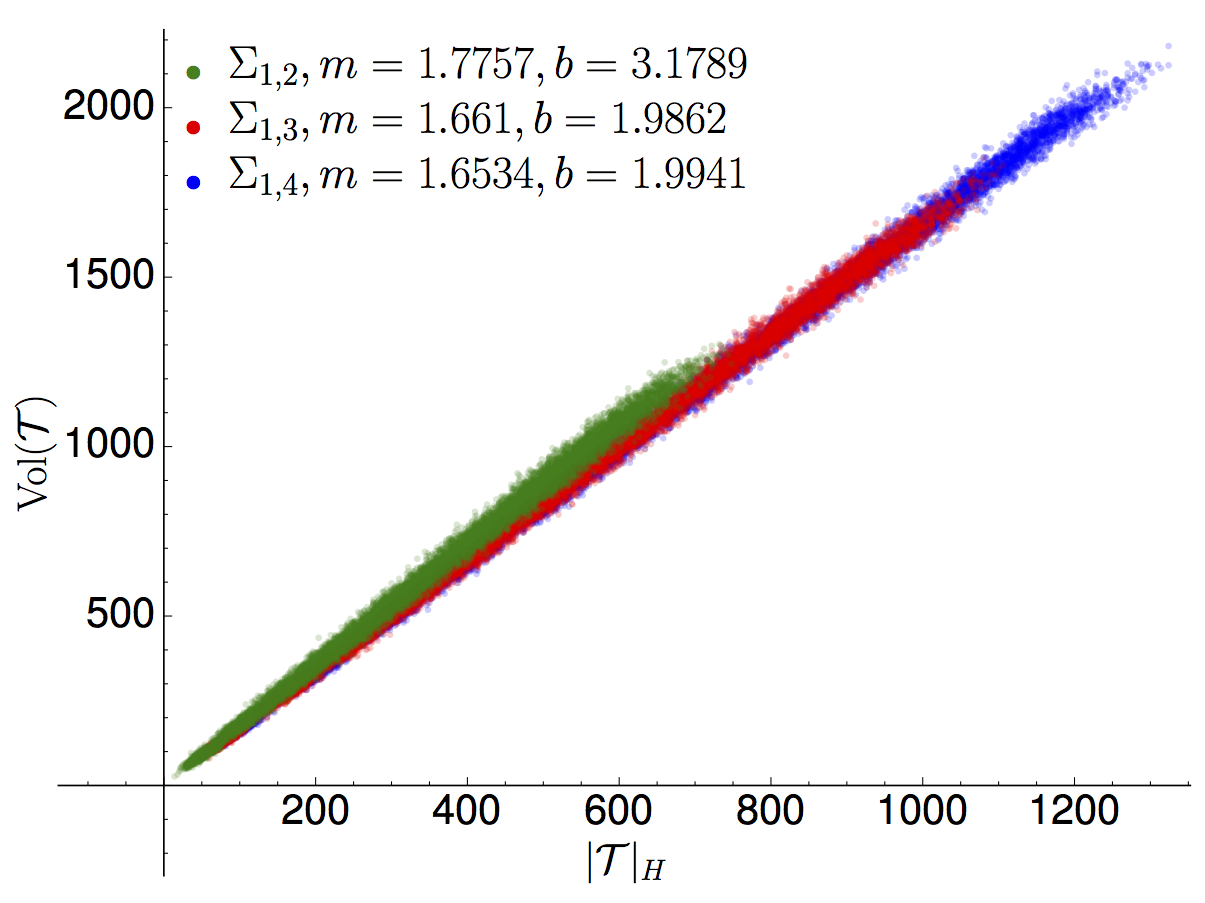}
        		\caption{}
        		\label{fig:hingevol_g1}
       \end{subfigure}
        \begin{subfigure}[b]{0.49\textwidth}
                \includegraphics[scale=.38]{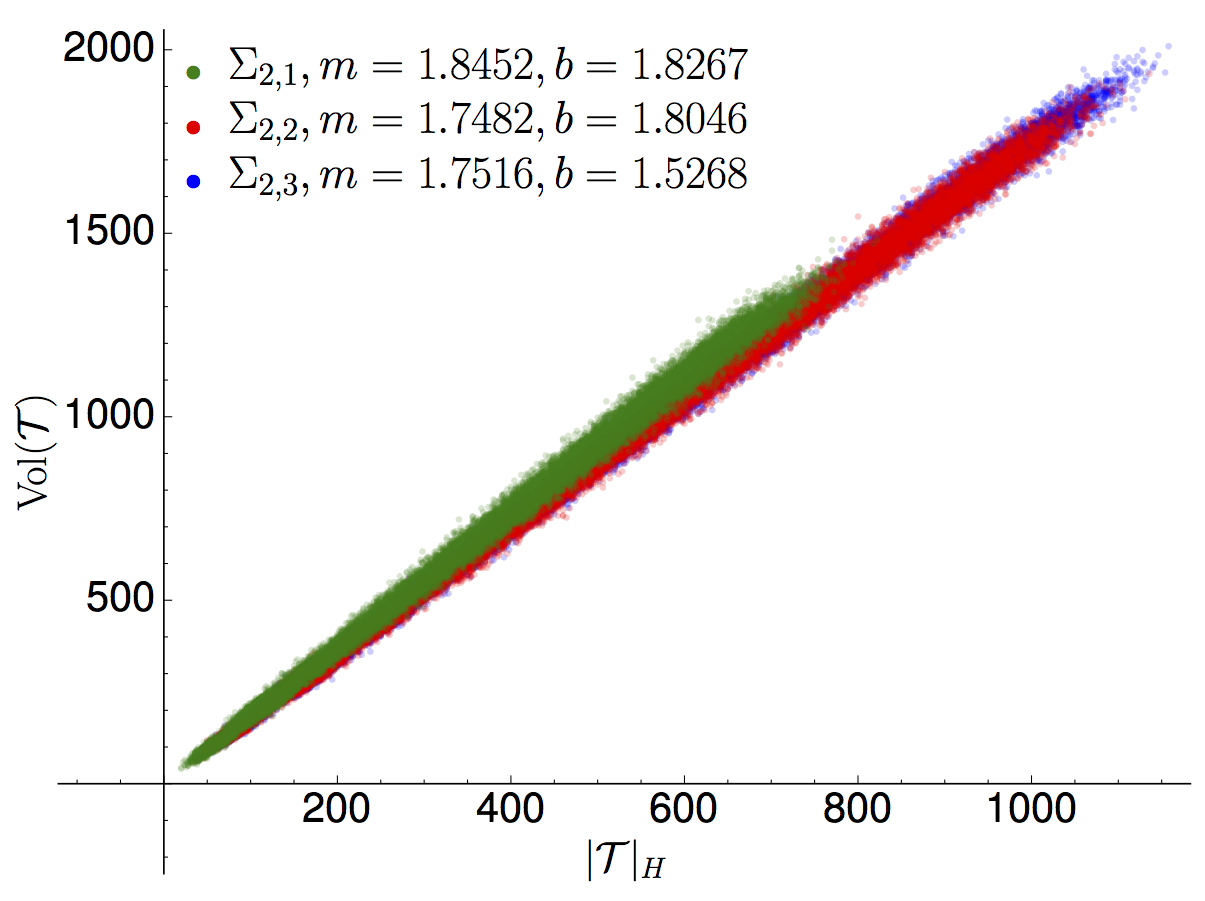}
                \caption{}
                \label{fig:hingevol_g2}
        \end{subfigure}       
       \caption{Here we break the scatter plots of \Cref{fig:hingevol_genus} up by number of punctures, for genus $1$ and $2$.}
        \label{fig:hingevol_g12}
\end{figure}

Returning to our experimental results, \Cref{fig:hingevol_genus} shows scatter plots of $|\mathcal{T}|_H$ against volume, for $\mr{genus}(\Sigma)$ ranging from $0$ to $3$. The surfaces that appear in this scatter plot are $\Sigma_{0,n}$ for $n\in\{5,8,11\}$, $\Sigma_{1,n}$ and $\Sigma_{2,n}$ for $n\in \{1,2,3\}$, and $\Sigma_{3,1}$, and there are approximately 12000 of each. In this graph we immediately see a striking linear relation, with a very narrow range of slopes  $\mr{Vol}(\mc{T})/|\mathcal{T}|_H$, even across multiple genera. At first glance it would probably be tempting to conjecture from this plot that the expected value of $\mr{Vol}(\mc{T})/|\mathcal{T}|_H$ increases with genus. There is reason for skepticism here, though. First, we only have low genus examples because of computational limits, so our sample size, in terms of genus, is very small. Second, if we consider, for example, the genus 1 examples, we find that slope does not appear to increase with complexity, and similarly for genus 2 (see Figures \Cref{fig:hingevol_g1} and \Cref{fig:hingevol_g2}). So it very well may be that higher complexity genus 3 surfaces would have slopes much less than $1.8133$, on average.

\begin{figure}
        \centering
        \includegraphics[scale=.7]{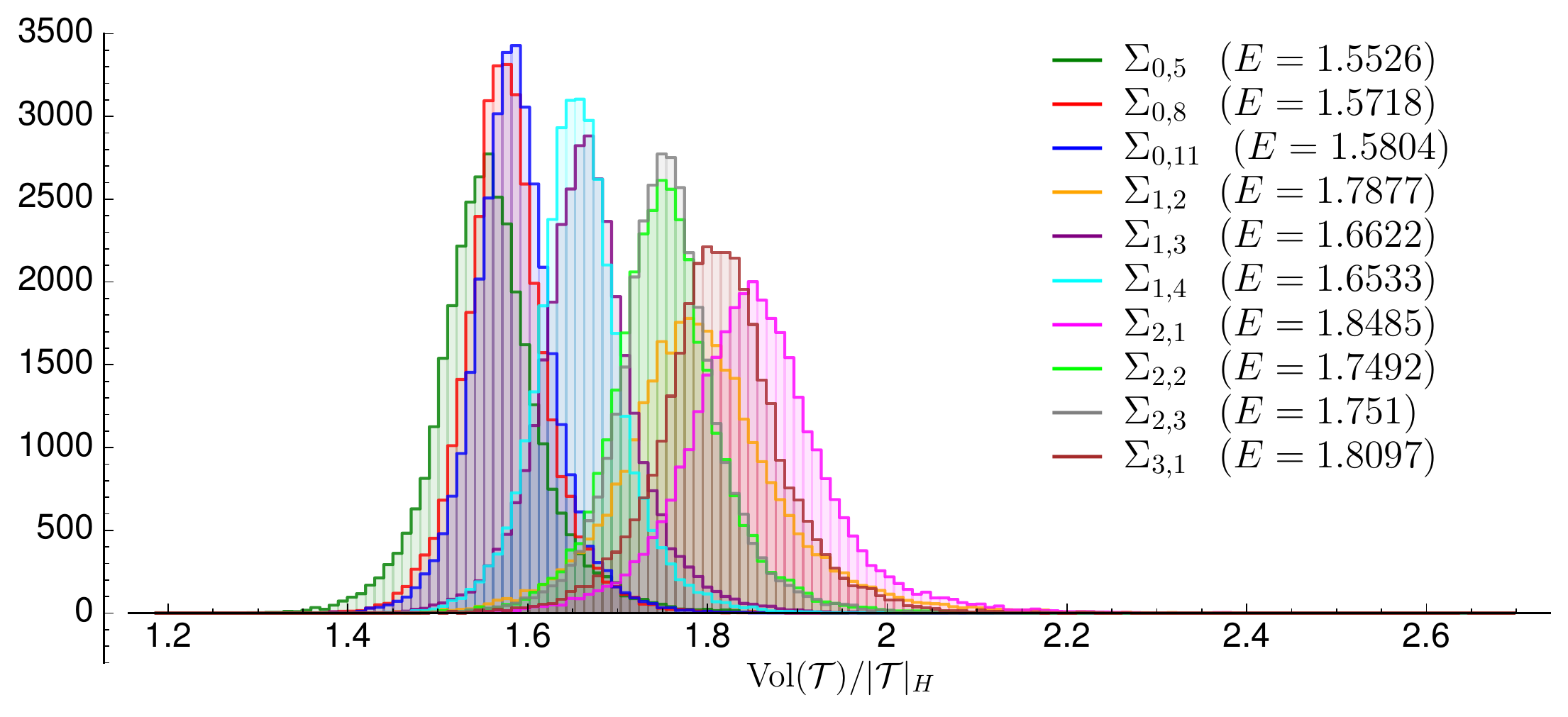}
       \caption{Left: histograms of $\mr{Vol}(\mc{T})/|\mc{T}|_H$, for veering triangulations with associated surfaces as indicated. $E$ is the expected value, or mean, of the histogram.}
        \label{fig:hingevol_hist_all}
\end{figure}

Another way to visualize the data in \Cref{fig:hingevol_genus} is to plot the histogram of the quotients $\mr{Vol}(\mc{T})/|\mathcal{T}|_H$ for each filled fiber $\Sigma$, as in \Cref{fig:hingevol_hist_all}. Note that for each of these surfaces the histogram \emph{appears} to approximate a normal distribution, however a QQ plot shows that in all cases there is a right skew. The value $E$ displayed for each of these histograms is the expected value of the histogram (as a discrete probability distribution), which will always be a close approximation of the mean of the underlying data. Unsurprisingly, these means are very close to the best fit slopes of the corresponding scatter plots, as one can see for some of the surfaces by comparing the values in Figures \Cref{fig:hingevol_g1} and \Cref{fig:hingevol_g2} to those in \Cref{fig:hingevol_hist_all}.

There is an upper bound on the quantity $\mr{Vol}(\mc{T})/|\mathcal{T}|_H$ in terms of the Euler characteristic of the fiber, due to Aougab--Brock--Futer--Minsky--Taylor. We state their result here, and give a proof which was communicated to the author by Dave Futer, for the sake of completeness. But first, recalling the construction of the veering triangulation as described in \Cref{sec:bg}, we make the following observation: For each maximal splitting $(\sigma,\mu)\rightharpoonup (\sigma',\mu')$ in the train track splitting sequence, one or more tetrahedra are layered onto $\Sigma^\circ$. Hence the tetrahedra in $\mc{T}$ are naturally partitioned into subsets, according to the layer that they belong to. If we order these subsets according to the order of the train track splitting, then order elements of each subset arbitrarily, we introduce an ordering on the tetrahedra in $\mc{T}$. This is a cyclic ordering, since the train track sequence is cyclic. In \Cref{sec:systole} it will be convenient to make it a total order, which we can do by choosing some arbitrary tetrahedra to be the first element. For what follows we will also need the following definition, due to Minsky--Taylor \cite{MiTa}:
\bg{defin}
A \textbf{pocket} in $\mc{T}$ is a union $P=\bigcup_i \Delta_i$ of tetrahedra $\Delta_i\in \mc{T}$, such that $\mr{int}(P)$ is connected; $\del P = \del P^-\cup\del P^+$, with $\del P^-\cong \del P^+$; and the induced triangulation on $\del P^+$ is obtained by performing edge flips on the triangulation of $\del P^-$, each of which corresponds to layering on a tetrahedron $\Delta_i$ (hence $\del P^-$ is below $P$ with respect to the transverse orientation, and $\del P^+$ is above).

\end{defin}

\bg{prop}[Aougab--Brock--Futer--Minsky--Taylor]\label{ABFMT} Let $M$ be a mapping torus with fiber $\Sigma^\circ$ and veering triangulation $\mathcal{T}$, and let $|\mathcal{T}|_H$ be the number of hinge tetrahedra in $\mathcal{T}$. Then
$$
\frac{\mr{Vol}(\mc{T})}{|\mathcal{T}|_H} \le v_3 (4|\chi(\Sigma^\circ)|+1),
$$
where $v_3$ is the volume of an ideal regular tetrahedron.
\end{prop}

\proof[Proof of proposition] The proof follows a general strategy introduced by Lackenby \cite{La}: drill out certain closed curves (this always increases the volume), retriangulate in a convenient way, and count the number of tetrahedra in the resulting triangulation. Since the volume of a hyperbolic manifold is always no greater than $v_3$ times the number of tetrahedra in any triangulation, this will give the desired result.

Consider a hinge tetrahedron $\Delta_h$, and let $\Delta_{1}, \dots, \Delta_{m}$ be the non-hinge tetrahedra appearing between $\Delta_h$ and the next hinge in the cyclic ordering. Let $\tilde{\Delta}_{1}, \dots, \tilde{\Delta}_{m}$ be lifts of the $\Delta_i$ to the infinite cyclic cover, all of which are contained in a single fundamental domain for the action of the monodromy. Let $\mathbf{T}=\bigcup_i \tilde{\Delta}_i$, and let each of $\mathbf{T}_1,\mathbf{T}_2,\dots, \mathbf{T}_k$ be a union of tetrahedra such that $\{\mr{int}(\mathbf{T}_i)\}_i$ are the connected components of $\mr{int}(\mathbf{T})$. Since a non-hinge with more red edges than blue (henceforth \textit{red} non-hinges) can never meet along a face of a non-hinge having more blue edges than red (\textit{blue} non-hinges), each component $\mathbf{T}_i$ is composed of either all red non-hinges or all blue non-hinges. Now, let $\mathbf{T}_i$ be a component which we may assume, without loss of generality, contains only red non-hinges. Since red non-hinges are formed by flipping red edges, no blue edges change as the tetrahedra in $\mathbf{T}_i$ are layered on, and hence all blue edges must be contained in $\del\mathbf{T}_i$. Since $\mathbf{T}_i$ is triangulated by red non-hinges, every tetrahedra face in $\mathbf{T}_i$ is a triangle with two red edges and one blue edge, and of these the triangles in $\del \mathbf{T}_i$ can be considered to be either below $\mathbf{T}_i$ or above $\mathbf{T}_i$ depending upon their coorientation. In particular, $\mathbf{T}_i$ is a pocket, sandwiched between a lower boundary $\del\mathbf{T}_i^-$ and an upper boundary $\del\mathbf{T}_i^+$ of triangular faces. 

\begin{figure}
        \centering
        \begin{subfigure}[b]{0.30\textwidth}
                \includegraphics[scale=.4]{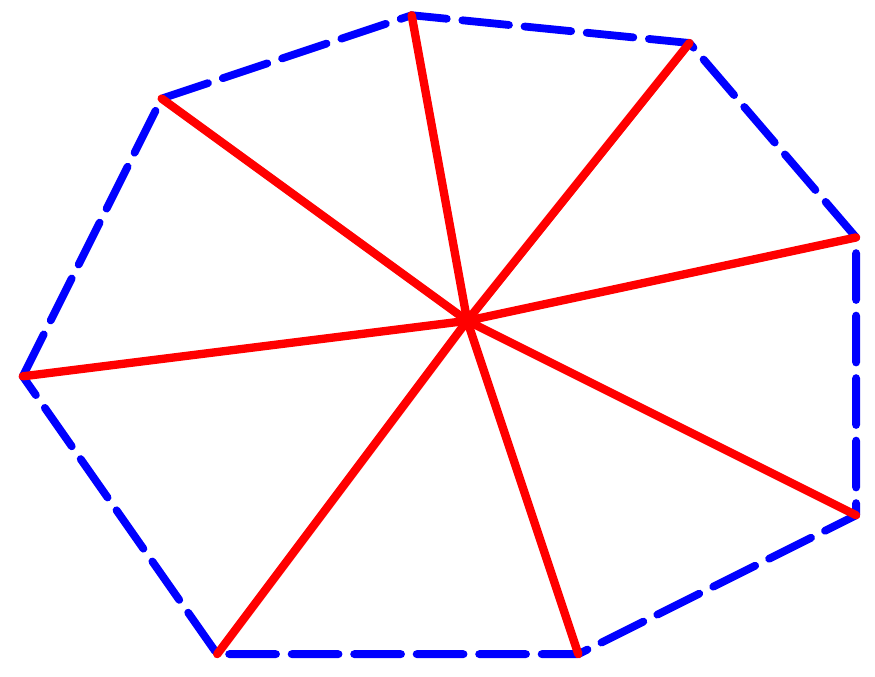}
                \caption{}
                \label{fig:punc_d}
        \end{subfigure}     
        \begin{subfigure}[b]{0.35\textwidth}
                \includegraphics[scale=.55]{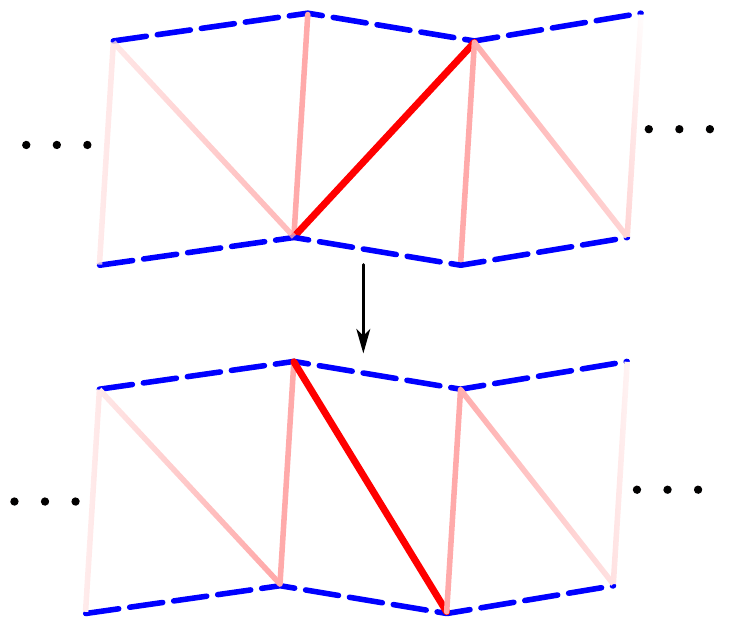}
        		\caption{}
        		\label{fig:flip}
       \end{subfigure}
        \begin{subfigure}[b]{0.30\textwidth}
                \includegraphics[scale=.55]{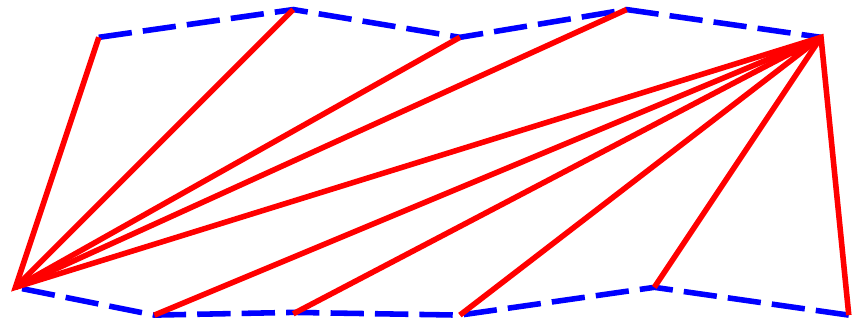}
        		\caption{}
        		\label{fig:bad_strip}
       \end{subfigure}
       \caption{Left: If $\del\mathbf{T}^-$ has a vertex meeting no blue edges, then none of its edges are flippable. Center: For an edge flip that results in a red non-hinge, the upper vertex shifts left, and the lower shifts right. Right: The worst case scenario for a strip of triangles that does not close up to form an annulus.}
        \label{fig:strips}
\end{figure}

To estimate the volume of $M$, we will need to retriangulate each pocket $\mathbf{T}_i$. First, we observe that every vertex of the induced triangulation of $\del\mathbf{T}_i^-$ must meet a blue edge, otherwise $\del\mathbf{T}_i^-$ would be as in \Cref{fig:punc_d}, and no edge flips would be possible. It follows that $\del\mathbf{T}_i^-$ is a strip of triangles as in Figures \Cref{fig:flip} and \Cref{fig:bad_strip}, which may or may not close up to form an annulus. If $\del\mathbf{T}_i^-$ does not close up, then it is a disk and $\mathbf{T}_i$ is a ball. In this case we can retriangulate $\mathbf{T}_i$ by picking a vertex $v$ and coning to it from each triangular face of $\del\mathbf{T}_i$ that does not have $v$ as a vertex. This new triangulation will have the same number of tetrahedra as there are triangles in $\del \mathbf{T}_i$.

If $\del\mathbf{T}_i^-$ does close up to form an annulus, then $\mathbf{T}_i$ is a solid torus. In this case, let $\alpha_i$ be the core curve of $\mathbf{T}_i$, and replace $\mathbf{T}_i$ by $\mathbf{T}_i\setminus \alpha_i$. We can triangulate our new $\mathbf{T}_i$ by coning the faces of $\del\mathbf{T}_i$ to the new vertex created at the removed core curve, so that $\mathbf{T}_i$ again has the same number of tetrahedra as $\del \mathbf{T}_i$.

With each $\mathbf{T}$ retriangulated as described above, let $\mc{T}'$ be the union of their projections $\pi(\mathbf{T})$, where $\pi$ is the covering map for the infinite cyclic cover. In other words, $\mc{T}'$ is the new triangulation obtained by carrying out the coning and drilling, as described above, in $\mathcal{T}$. Let $M'$ be the drilled manifold whose triangulation is $\mathcal{T}'$.

Note that we can consider $\del\mathbf{T}_i^-$ to be a subset of a fiber $\Sigma^\circ$ (in the infinite cyclic cover), and this is precisely the transverse projection of $\mathbf{T}_i$ onto $\Sigma^\circ$. Furthermore, all the projections $\del\mathbf{T}_j^-$ for $1\le j\le k$ are disjoint in $\Sigma^\circ$ (except possibly along edges). Otherwise, there would have to be some tetrahedra $\Delta' \not\subset \mathbf{T}$ layered on after the tetrahedra of some $\del\mathbf{T}_{j_1}^-$ and before those of some $\del\mathbf{T}_{j_2}^-$, which is impossible given how $\mathbf{T}$ was defined. 

Since any triangulation of $\Sigma^\circ$ has $2|\chi(\Sigma^\circ)|$ triangles, and since $\del \mathbf{T}^-=\bigcup_i \del\mathbf{T}_i^-$ is a subsurface (after projecting) of $\Sigma^\circ$, and similarly for $\del \mathbf{T}^+$, it follows that the number of tetrahedra in $\mathbf{T}$, which is equal to the number of triangles in $\del\mathbf{T}=\del\mathbf{T}^-\cup \del\mathbf{T}^+$, is at most $4|\chi(\Sigma^\circ)|$. Hence the triangulation $\mc{T}'$ has at most $|\mc{T}|_H(4|\chi(\Sigma^\circ)|+1)$ tetrahedra. Since the Gromov norm $\| \cdot\|$ of a manifold increases with drilling, and is no greater than the number of tetrahedra in a triangulation, we get $\mr{Vol}(\mc{T}) \le v_3 \|M\| \le v_3 \|M'\| \le v_3|\mc{T}|_H(4|\chi(\Sigma^\circ)|+1)$.
\qed

\begin{figure}[h]
        \centering
        \begin{subfigure}[b]{0.5\textwidth}
        	\includegraphics[scale=.37]{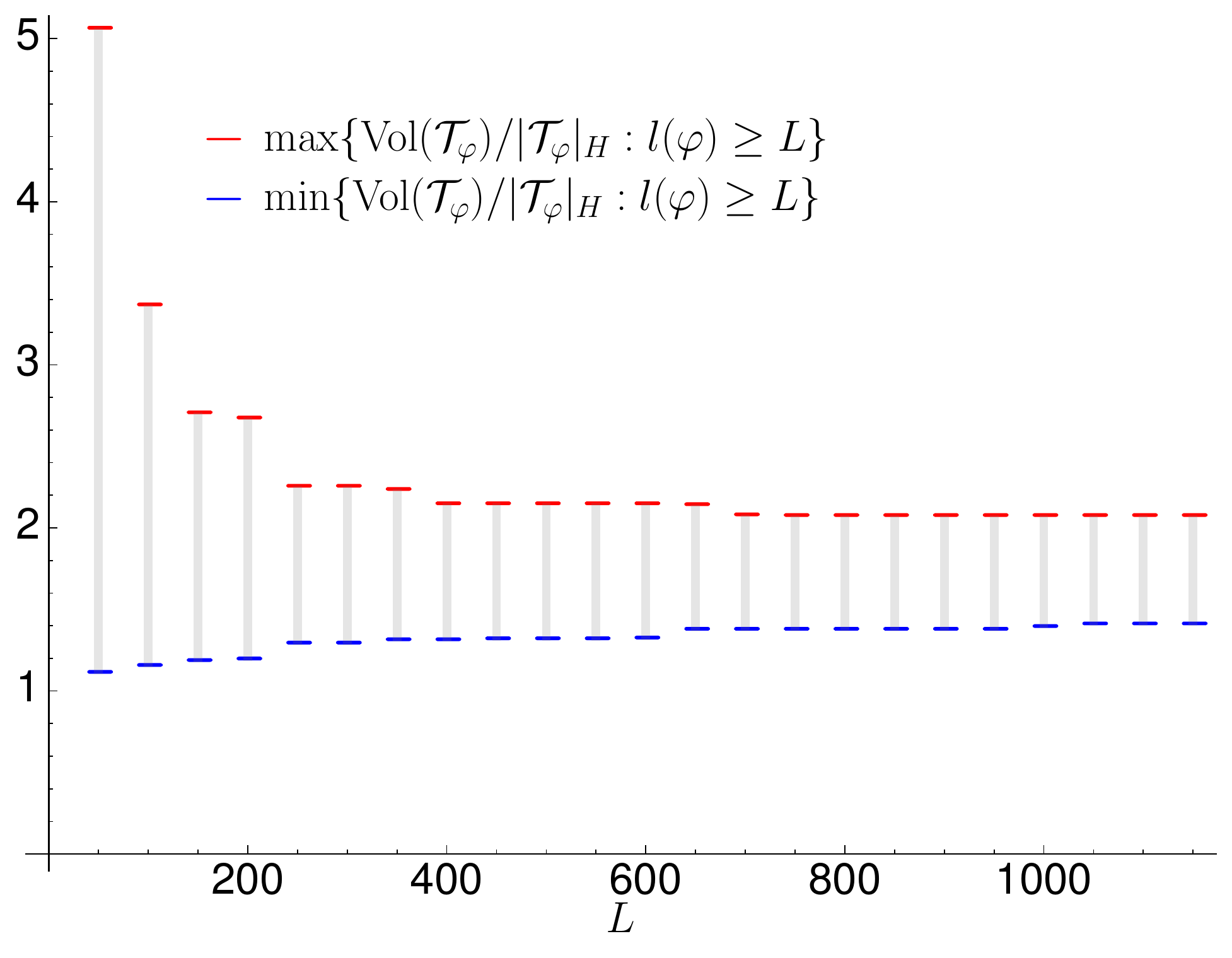}
        	\caption{}
        	\label{fig:slope_ranges}
       \end{subfigure}
       \begin{subfigure}[b]{0.49\textwidth}
                \includegraphics[scale=.37]{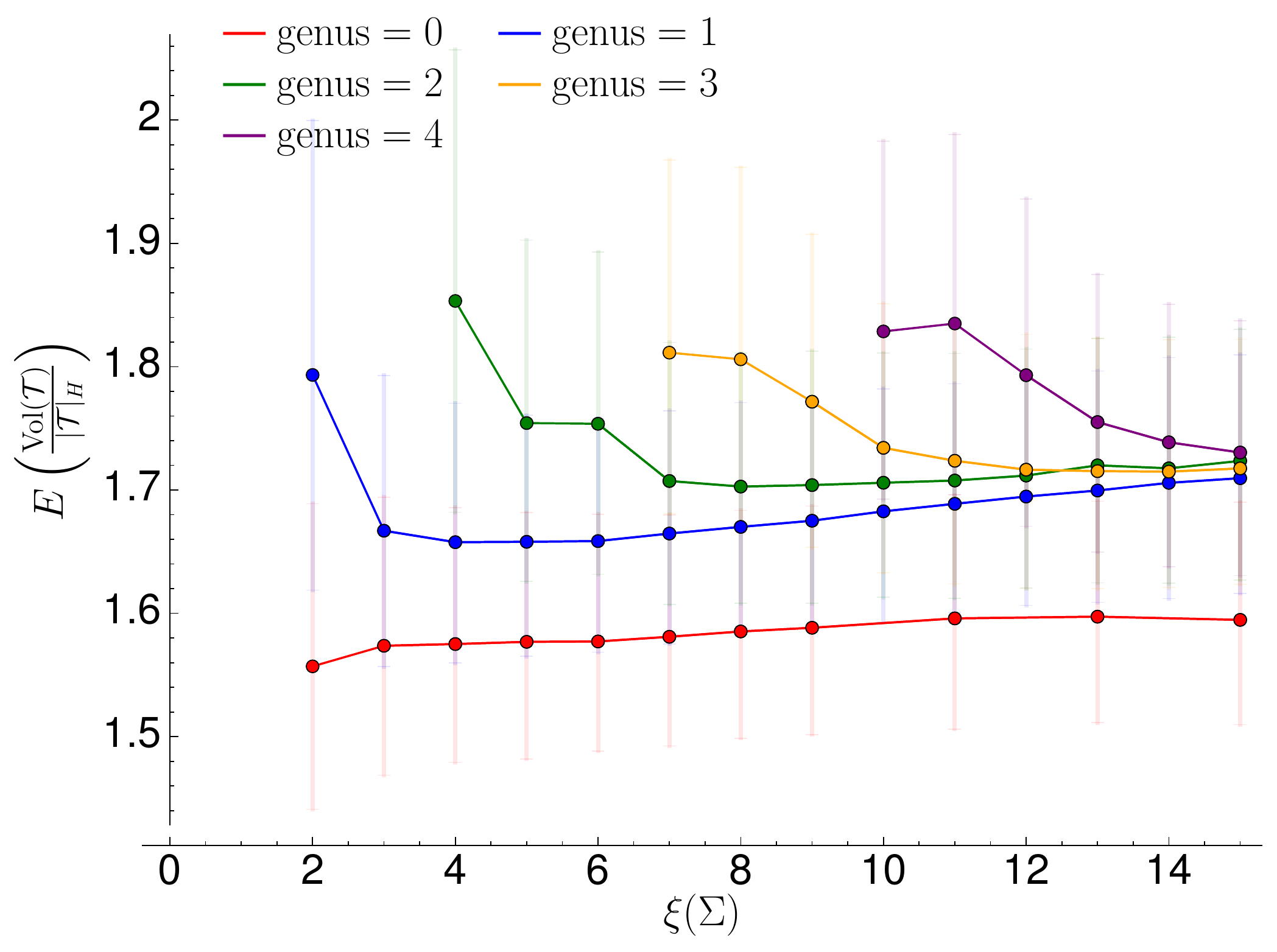}
                \caption{}
                \label{fig:Eval}
        \end{subfigure}
       \caption{Left: At each $L$ value, the vertical bar shows the range of values $\mr{Vol}(\mc{T})/|\mc{T}|_H$, for all triangulations with associated mapping classes $\varphi$ of length $l(\varphi) \ge L$. Right: surface complexity $\xi(\Sigma)$, plotted against the mean of the values $\mr{Vol}(\mc{T})/|\mc{T}|_H$ over all triangulations with that associated surface. Each dot is a surface, and dots for surface of the same genus are joined by line segments and color coded, as indicated. Vertical bars show the range of the middle 90\% of each distribution.}
        \label{fig:}
\end{figure}

The above proposition and the experimental results of this section naturally lead to the question of whether there exists a universal upper and lower bound for the slope $\mr{Vol}(\mc{T})/|\mathcal{T}|_H$, independent of the filled fiber $\Sigma$. For the $\sim 750000$ examples in our data set, \Cref{fig:slope_ranges} shows the range of the values $\frac{\mr{Vol}(\mc{T})}{|\mc{T}|_H}$ realized for all words of length greater than $x$. For the shortest words, of length $50\le l\le 150$, the slope can be quite large, and is likely much larger for words in the $0\le l\le 50 $ range (which are not part of our sample). As $L$ increases, the range of slopes for generic mapping classes becomes narrow, according to \Cref{fig:slope_ranges}, suggesting the least upper bound for large word length could be quite low. This of course cannot actually be the case, since if $\varphi$ is a short word with large slope, the finite cover obtained from the mapping class $\varphi^k$ will have the same slope for any $k$. Furthermore, our data is limited to low genus examples, and it is unclear thus far whether slope should be expected to increase as genus increases.

In order to get a better handle on the relation between slope and genus, we will need data for surfaces of higher complexity. Unfortunately, computing veering triangulations for large words becomes intractable as complexity increases. To get around this limitation, we observe from \Cref{tab:word_len} that the expected value of $\mr{Vol}(\mc{T})/|\mc{T}|_H$ for a sample of words of length $100\le l(\varphi)\le 1200$, is fairly reliably predicted by sampling a large number of mapping classes all of length $\sim200$. \Cref{fig:Eval} shows the expected value of $\mr{Vol}(\mc{T})/|\mc{T}|_H$ of each of $47$ surfaces $\Sigma$, plotted against the complexity $\xi(\Sigma)$. For each genus, the plot points for surfaces of that genus are joined by line segments and color coded per the key.

\begin{table}
\centering
\caption{The average value of $\mr{Vol}(\mc{T})/|\mc{T}|_H$ over 12000 veering triangulations with associated mapping class words of length $100\le l(\varphi)\le 1200$, is very close to the average we get from a sample of 10000 triangulations all with word length $l(\varphi)=200$.}

\begin{tabular}{|r|r|r|r|r|r|} \hline
 & $\Sigma_{0,5}$ & $\Sigma_{0,8}$ & $\Sigma_{0,11}$ & $\Sigma_{1,2}$ & $\Sigma_{1,3}$ \\ \hline \hline
 $l(\varphi)\le 1200$ & $1.5546$ & $1.5742$ & $1.5827$ & $1.7874$ & $1.6639$ \\ \hline
 $l(\varphi)=200$ & $1.557$ & $1.5769$ & $1.5853$ & $1.7933$ & $1.667$ \\ \hline
\end{tabular}

\hspace*{0.005cm}
\begin{tabular}{|r|r|r|r|r|r|} \hline
 & $\Sigma_{1,4}$ & $\Sigma_{2,1}$ & $\Sigma_{2,2}$ & $\Sigma_{2,3}$ & $\Sigma_{3,1}$ \\ \hline \hline
 $l(\varphi)\le 1200$ & $1.6567$ & $1.8484$ & $1.7514$ & $1.7537$ & $1.8112$ \\ \hline
 $l(\varphi)=200$ & $1.6577$ & $1.8533$ & $1.7543$ & $1.7537$ & $1.8114$ \\ \hline
\end{tabular}
\label{tab:word_len}
\end{table}

It seems quite plausible that all of the subgraphs in \Cref{fig:Eval} (except for genus $= 0$) could eventually converge, say to slightly over 1.7, but we would need a great deal more data points to form a hypothesis with any confidence. Unfortunately, we have pushed the complexity for each genus to the limit of what can be computed in a reasonable amount of time, as was demonstrated in \Cref{fig:time}. 

As to the question of whether there is a universal \emph{lower} bound on $\mr{Vol}(\mc{T})/|\mathcal{T}|_H$, our data points more clearly to a positive answer. In fact, the existence of a lower bound would follow from a conjecture of Casson \cite{FuGu:survey}, if it were proved. This conjecture asserts that the volume of any non-negative angle structure is less than or equal to the volume of the manifold. By a result of Futer--Gu\'{e}ritaud \cite{FuGu:veering}, any veering triangulation $\mc{T}$ has an angle structure in which all hinge tetrahedra are regular (i.e., all angles are $\frac{\pi}{3}$), and all non-hinge tetrahedra are flat. For this angle structure, the total volume of the tetrahedra (i.e., the volume of the angle structure) is $v_3\cdot |\mathcal{T}|_H$, where $v_3$ is the volume of a regular ideal tetrahedron. Combining this with Casson's conjecture, we get that $\mr{Vol}(\mc{T})/|\mathcal{T}|_H \ge v_3$. In fact, if this inequality holds then it is sharp, since equality is attained when $M_{\varphi^\circ}$ is the figure-8 knot complement (or any other manifold in its fibered commensurability class).

\bg{rem}
It is worth noting that our experimental data also provides evidence \textit{in favor} of the conjecture of Casson mentioned above (though in a restricted setting). In addition to the hundreds of thousands of generic mapping classes underlying the data in this section which satisfy the conjecture, \Cref{sec:other:geo} gives a family of non-generic examples for which it appears that $\mr{Vol}(\mc{T})/|\mathcal{T}|_H$ is converging to $v_3$ from above.
\end{rem}

\section{Systole Length and Chains of Non-Hinge Tetrahedra}
\label{sec:systole}

The \textbf{systole} of $M_{\varphi^\circ}$ is the shortest (closed, possible non-unique) geodesic in $M_{\varphi^\circ}$. Recall from the previous section that the tetrahedra in the veering triangulation $\mc{T}$ can be given a total ordering consistent with the partial order induced by the layering. Call a chain of consecutively ordered tetrahedra in $\mc{T}$ \textbf{monochromatic} if every tetrahedra in the chain has the same number of red edges (and hence also of blue edges). Let a \textbf{maximal (monochromatic) non-hinge chain} $M_\mr{ch}=\{\Delta_1,\dots, \Delta_m\}$ be a maximal length connected (along faces) sub-chain of a monochromatic chain of non-hinges. That is, we consider all monochromatic chains of non-hinges, and separate each of these into components that are connected via face gluings, then take the largest resulting component. In this section we demonstrate a relationship between the systole length and the length of the maximal non-hinge chain length $|M_{\mr{ch}}|$. These experiments are motivated by the Length Bound Theorem of Brock--Canary--Minsky, and in particular our results are consistent with a corollary to the Length Bound Theorem, given below. To state this corollary we will first need some definitions.

\begin{defin}
The \textbf{arc and curve complex} $\mc{A}(Y)$ of a compact surface $Y$ is the simplicial complex whose vertices are homotopy classes of essential simple closed curves and properly embedded arcs (if $Y$ is an annulus the homotopies must additionally fix the endpoints of arcs), and whose simplices correspond to tuples of disjoint arcs/curves.
\end{defin}

\begin{defin}
Given a subsurface $Y\subset \Sigma$ and stable/unstable laminations $\mc{L}^+,\mc{L}^-$, define the \textbf{subsurface distance} $d_Y(\mc{L}^+,\mc{L}^-)$ between $\mc{L}^+$ and $\mc{L}^-$ to be the distance in the arc and curve complex between the lifts of $\mc{L}^\pm$ to the (closure of the) cover of $\Sigma$ homeomorphic to $Y$.
\end{defin}

\begin{figure}
        \centering
        \begin{subfigure}[b]{0.46\textwidth}
                \includegraphics[scale=.8]{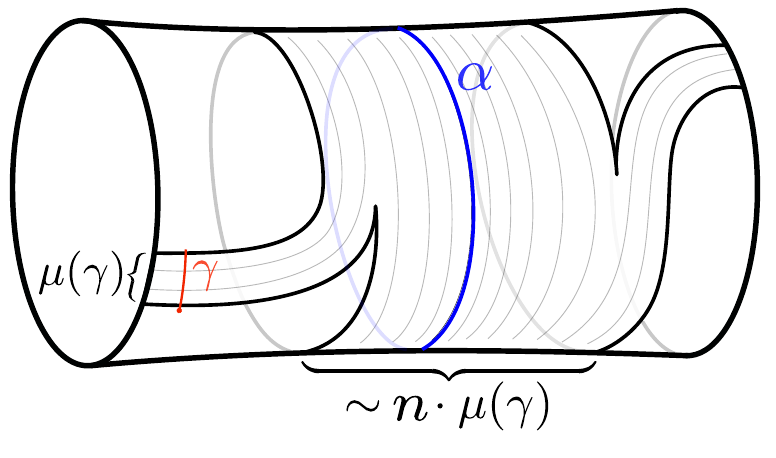}
                \caption{}
                \label{fig:annulus_lam}
        \end{subfigure}     
        \begin{subfigure}[b]{0.53\textwidth}
                \includegraphics[scale=1.2]{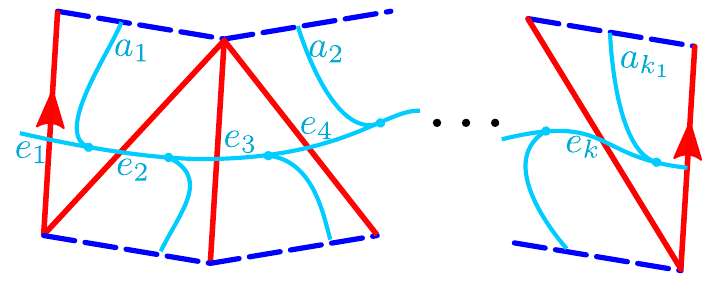}
        		\caption{}
        		\label{fig:strip}
       \end{subfigure}
       \caption{Left: An intersection of an annulus $A$ with the stable lamination $\mc{L}^+$: if $\mc{L}^+$ wraps $n$ times around the annulus, then $d_\alpha(\mc{L}^+,\mc{L}^-)=2n \pm 1$. Right: A possible train track for a more complicated intersection of $\mc{L}^+$ with an annulus. The weight of each $e_i$ branch will be roughly $n$ times the total weight of the $a_i$ branches, if the lamination wraps $n$ times around the annulus.}
        \label{fig:annulus}
\end{figure}

If $Y$ is an annulus with core curve $\delta$, then we will write $d_\delta(\mc{L}^+,\mc{L}^-)=d_Y(\mc{L}^+,\mc{L}^-)$. In this case, the subsurface distance just measures, within $\pm 1$, the number of times a component of $Y\cap \mc{L}^+$ wraps arround $Y$, relative to a component of $Y\cap \mc{L}^-$. For example, in \Cref{fig:annulus_lam}, if each component of $A\cap \mc{L}^+$ wraps $n$ times around $A$ in one direction, then each component of $A\cap\mc{L}^-$ (not shown) wraps $n$ times around $A$ in the other direction, so $d_\alpha(\mc{L}^+,\mc{L}^-)$ will be about $2n$. In this case the corresponding train track in the annulus will have branches parallel to the annulus with weight approximately $n$ times the total weight of edges entering the annulus from one side. 

For a surface $\Sigma$ and pseudo-Anosov mapping class $\varphi$, let $l_\varphi(\alpha)$ denote the hyperbolic length of the curve $\alpha$ in the mapping torus $M_\varphi$. The following is the Length Bound Theorem of Brock--Canary--Minsky, given in a form that appears as Corollary 7.3 in \cite{SiTa}. This version is somewhat weaker than the theorem in its original form, but it is well suited to our needs.

\bg{thm}[Length Bound Theorem \cite{BCM}]\label{LBT}There are $D>1,\epsilon> 0$, depending only on $\Sigma$, such that for any pseudo-Anosov $\varphi\in \mr{Mod}(\Sigma)$ and any curve $\alpha$ in $\Sigma$ with $l_\varphi(\alpha)\le \epsilon$:
$$
D^{-1}\cdot \frac{2\pi \cdot S_\varphi(\alpha)}{d_\alpha^2(\mc{L}^+,\mc{L}^-)+S_\varphi^2(\alpha)}\le l_\varphi(\alpha)\le D\cdot \frac{2\pi \cdot S_\varphi(\alpha)}{d_\alpha^2(\mc{L}^+,\mc{L}^-)+S_\varphi^2(\alpha)}
$$
where $S_\varphi(\alpha):=1+\sum_{Y\in \mc{Y}_\alpha} [d_Y(\mc{L}^+,\mc{L}^-)]_K$. Moreover, if $d_\alpha^2(\mc{L}^+,\mc{L}^-)\ge E$, for $E$ depending only on $\Sigma$, then $l_\varphi(\alpha)\le \epsilon$.
\end{thm}

In the above definition of $S_\varphi(\alpha)$, $\mc{Y}_\alpha$ is the set of all subsurfaces of $\Sigma$ which have $\alpha$ as a boundary component, and $[\cdot]_K$ is a cutoff function defined by $[x]_K=x$ if $x\le K$ and $[x]_K=0$ otherwise. The constant $K$ is shown to exist in Minsky \cite{Mi:kleinian}, but will not be important for what follows. From the above theorem it follows immediately that if $d_\alpha(\mc{L}^+,\mc{L}^-)$ is sufficiently large, then $l_\varphi(\alpha)$ will be small, and with a little more effort we get the following corollary:

\bg{cor}
Let $M_{\varphi^\circ}$ be a mapping torus with fiber $\Sigma_\circ$, monodromy $\varphi^\circ$, and veering triangulation $\mc{T}$. Let $|M_{\mr{ch}}|$ be the length of the maximal non-hinge chain in $\mc{T}$, and let $\sigma$ be the systole in $M_{\varphi^\circ}$. There exists $C,N\ge 0$ depending only on $\Sigma_\circ$ such that if $|M_{\mr{ch}}| \ge N$ then
$$
l_\varphi(\sigma)\le C\cdot \frac{1}{|M_{\mr{ch}}|}
$$

\end{cor}

We prove the corollary below, but for those that wish to skip the details, here is the rough idea: large maximal non-hinge chains must be layered around annuli, and hence the train track branches that are parallel to the core curve of such an annulus have large weight, so that the lamination must wrap many times around that annulus and hence has large sub-surface projection. We then apply the Length Bound Theorem, which implies that the core curve of the annulus is short, and therefore so is the systole.

\proof[Proof of Corollary]

Suppose we have a maximal non-hinge chain $\Delta_1,\Delta_2,\dots, \Delta_m$ of tetrahedra in a triangulation $\mc{T}$, with $m=|M_{\mr{ch}}|\ge \max\left(\frac{4}{3}(|\chi |^2+|\chi|),\frac{E|\chi|+|\chi|^2}{2}\right)=N$, where $\chi=\chi(\Sigma_\circ)$ and $E$ is as in the Length Bound Theorem. Let $\mathbf{T}=\bigcup_i \Delta_i$. By our definition of a maximal non-hinge chain, $\mr{int}(\mathbf{T})$ is connected, and we may also assume that all of these non-hinges are red (i.e., they have more red edges than blue). 

Just as in the proof of \Cref{ABFMT}, $\mathbf{T}$ is a pocket with upper and lower boundaries $\del\mathbf{T}^+$ and $\del\mathbf{T}^-$, each of which is a strip of red-red-blue triangles, which may or may not close up to form annuli.

In the first case, in which the strip of triangles does not close up, $m=|M_{\mr{ch}}|$ will be bounded by the square of the number of triangles in the strip. To see this, observe that when an edge flips, its upper vertex shifts left and its lower vertex shifts right, as shown in \Cref{fig:flip}, so that the worst case scenario is a diagonal edge going from the bottom left of the strip to the top right, as shown in \Cref{fig:bad_strip}. If there are $k_1$ blue edges along the top of the strip, and $k_2$ edges along the bottom, then one can check that only $k_1\cdot k_2$ flips are possible before no more edges are flip-able. Since there are $k_1+k_2$ triangles in the strip and exactly $2|\chi(\Sigma_\circ)|$ in any triangulation of the surface, and since the strip is a proper subsurface (otherwise it would have to close up into an annulus), we get $m\le k_1\cdot k_2\le \left(\frac{k_1+k_2}{2}\right)^2<|\chi|^2$. 

If $\del\mathbf{T}^-$ is an annulus, however, then $m$ can be arbitrarily large, and in this case $\mathbf{T}$ will be a solid torus---in this case we will refer to $\mathbf{T}$ as an \textbf{annular pocket}. Let $\alpha$ be the core curve of $\del\mathbf{T}^-$, let $k=k_1+k_2$ be the number of triangles in $\del\mathbf{T}^-$, and let $\tau$ be the train track prior to the spitting corresponding to the layering on of $\Delta_1$. Then $\del\mathbf{T}^- \cap \tau$ has $k$ branches $e_1,\dots, e_k$ that are dual to red edges (i.e., parallel to the core curve of the annulus), and $k_1$ edges $a_1,\dots, a_{k_1}$ that exit the annulus through the upper blue edges (see \Cref{fig:strip}). Since there are $m$ tetrahedra in the pocket, at least one of these $e_i$ branches must split at least $\frac{m}{k}$ times. When a branch splits, it cannot be split again until the two adjacent branches split, and from this it follows that, if $e_1$ is the branch that splits at least $\frac{m}{k}$ times, then the adjacent branches $e_2$ and $e_k$ split at least $\frac{m}{k}-1$ times, and $e_3,e_{k-2}$ split at least $\frac{m}{k}-2$ times, and so on until we have gone fully around the annulus. It follows immediately that every branch $e_i$ splits at least $\frac{m}{k}- \frac{k}{2}$ times, but in fact we can do better: a more careful analysis (left to the reader) shows that every branch must split at least $\frac{m}{k}-\frac{k}{4}=\frac{4m-k^2}{4k}$ times. It follows that $d_\alpha(\mc{L}^+,\mc{L}^-)\ge 2\cdot \left(\frac{4m-k^2}{4k}\right) -1$, since every $e_i$ splitting at least $\frac{4m-k^2}{4k}$ times means that $\mu(e_i)\ge \left(\frac{4m-k^2}{4k}\right)\cdot \sum_j \mu(a_j)$, i.e., the lamination is wrapping around the annulus at least $\frac{4m-k^2}{4k}$ times. Replacing $m=|M_{\mr{ch}}|$ by $\frac{E|\chi|+|\chi|^2}{2}$ we obtain $d_\alpha(\mc{L}^+,\mc{L}^-) \ge E$, so that $l_\varphi(\alpha)\le \epsilon$.  Hence we have 
$$
l_\varphi(\alpha)\le D\cdot \frac{2\pi \cdot S_\varphi(\alpha)}{d_\alpha^2(\mc{L}^+,\mc{L}^-)+S_\varphi^2(\alpha)} \le \frac{2\pi D}{d_\alpha(\mc{L}^+,\mc{L}^-)} \le \frac{2\pi D}{2\left(\frac{2|M_\mr{ch}|-k^2}{2k}\right)-1}\le \frac{\pi D\cdot |\chi|}{|M_\mr{ch}|-|\chi|^2-|\chi|}\le\frac{4\pi D \cdot |\chi|}{|M_{\mr{ch}}|}
$$
since $k \le 2|\chi|$ and $\frac{3}{4}|M_\mr{ch}| \ge (|\chi|^2+|\chi|)$. Since the systole is the shortest closed curve in $M_{\varphi^\circ}$, the conclusion follows.
\qed

\pp
In the above proof we establish that if $|M_{\mr{ch}}|$ is large, then so is the sub-surface projection distance $d_\alpha(\mc{L}^+,\mc{L}^-)$. In fact the converse of this is also true. By a result of Minsky--Taylor \cite{MiTa}, if $d_\alpha(\mc{L}^+,\mc{L}^-)\ge 10$, then $d_\alpha(V^+,V^-)\ge d_\alpha(\mc{L}^+,\mc{L}^-)-10$, where $V^+$ (resp. $V^-$) is the set of tetrahedra edges in $\del\mathbf{T}^+$ (resp. $\del\mathbf{T}^-$), regarded as simplices in the arc and curve complex. Since $d_\alpha(V^+,V^-)$ is proportional to the number of tetrahedra $|M_\mr{ch}|$ in the pocket $\mathbf{T}$, the converse follows. A full converse of the corollary does not hold, however. That is, if the systole is very short, it does not necessarily follow that $|M_\mr{ch}|$ will be large---this is because if $S_\varphi(\alpha)$ is large but $d_\alpha(\mc{L}^+,\mc{L}^-)$ is small then the systole will be short, but $|M_\mr{ch}|$ only sees $d_\alpha(\mc{L}^+,\mc{L}^-)$.

\begin{figure}[h]
        \centering
        \begin{subfigure}[b]{0.5\textwidth}
                \includegraphics[scale=.38]{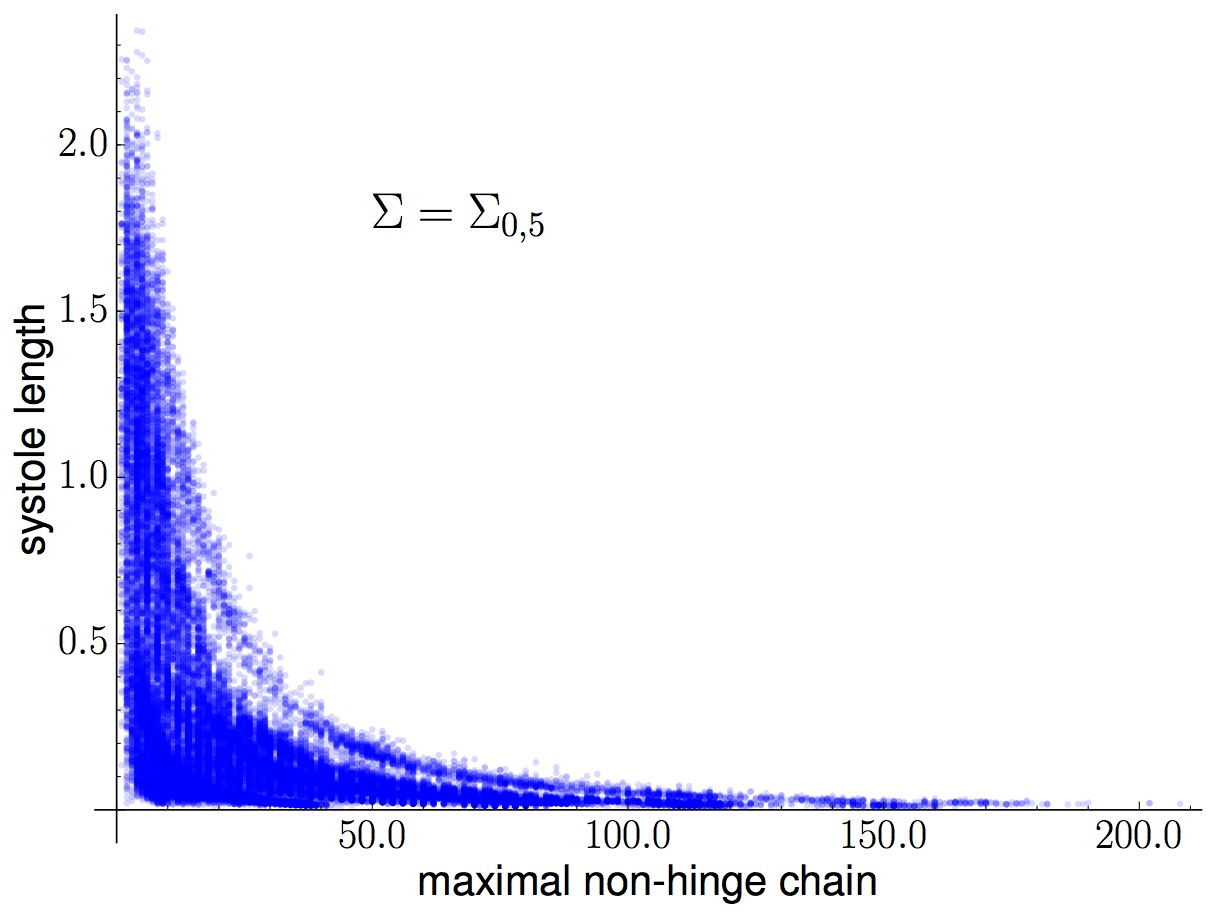}
        		\caption{}
        		\label{fig:chlen_sys_05_a}
       \end{subfigure}
        \begin{subfigure}[b]{0.49\textwidth}
                \includegraphics[scale=.38]{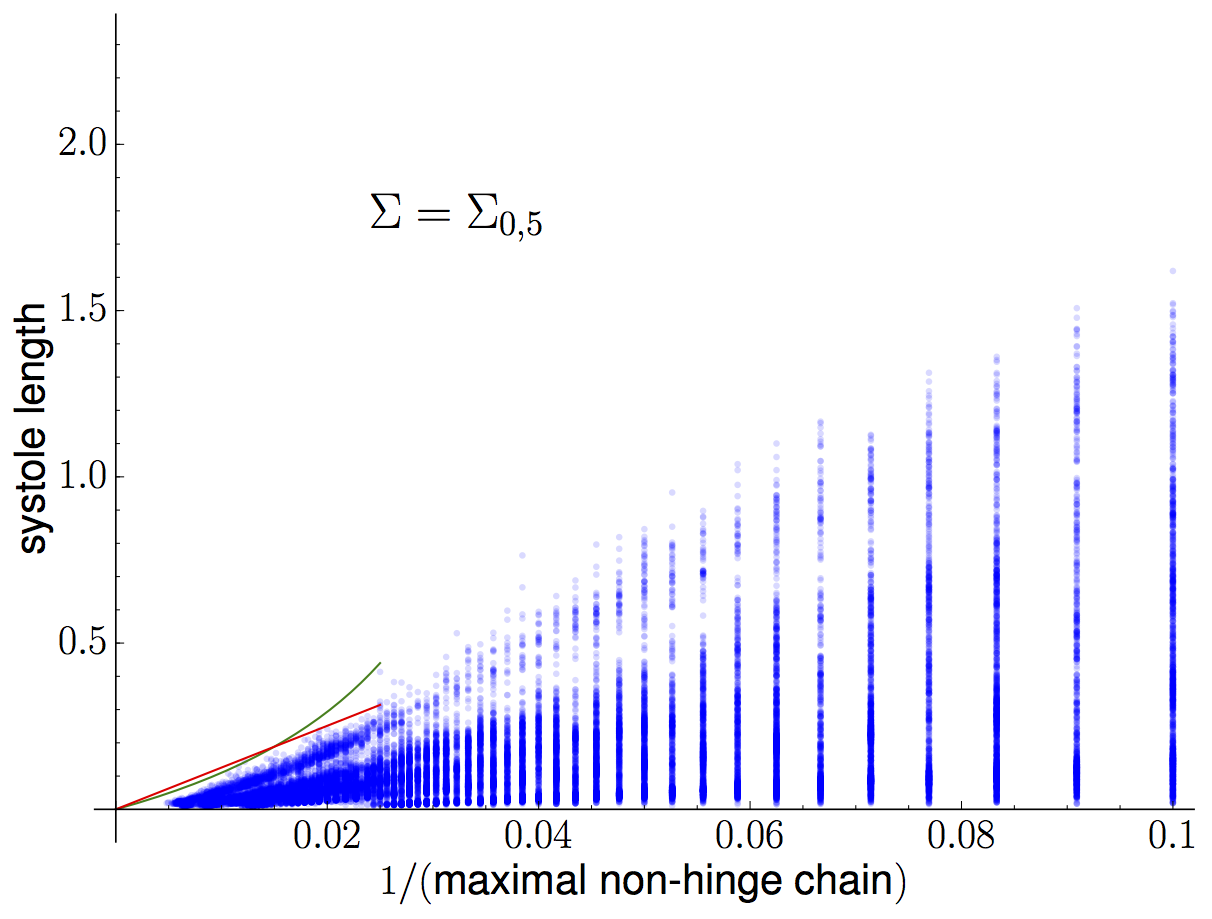}
                \caption{}
                \label{fig:chlen_sys_05_b}
        \end{subfigure}       
       \caption{The maximal chain of non-hinge tetrahedra (left) and its reciprocal (right) plotted against systole length, for $\sim 45000$ veering triangulations with associated surface $\Sigma_{0,5}$.\ppp{5}}
        \label{fig:chlen_sys_05}
\end{figure}

\begin{figure}
        \centering
        \begin{subfigure}[b]{0.5\textwidth}
                \includegraphics[scale=.38]{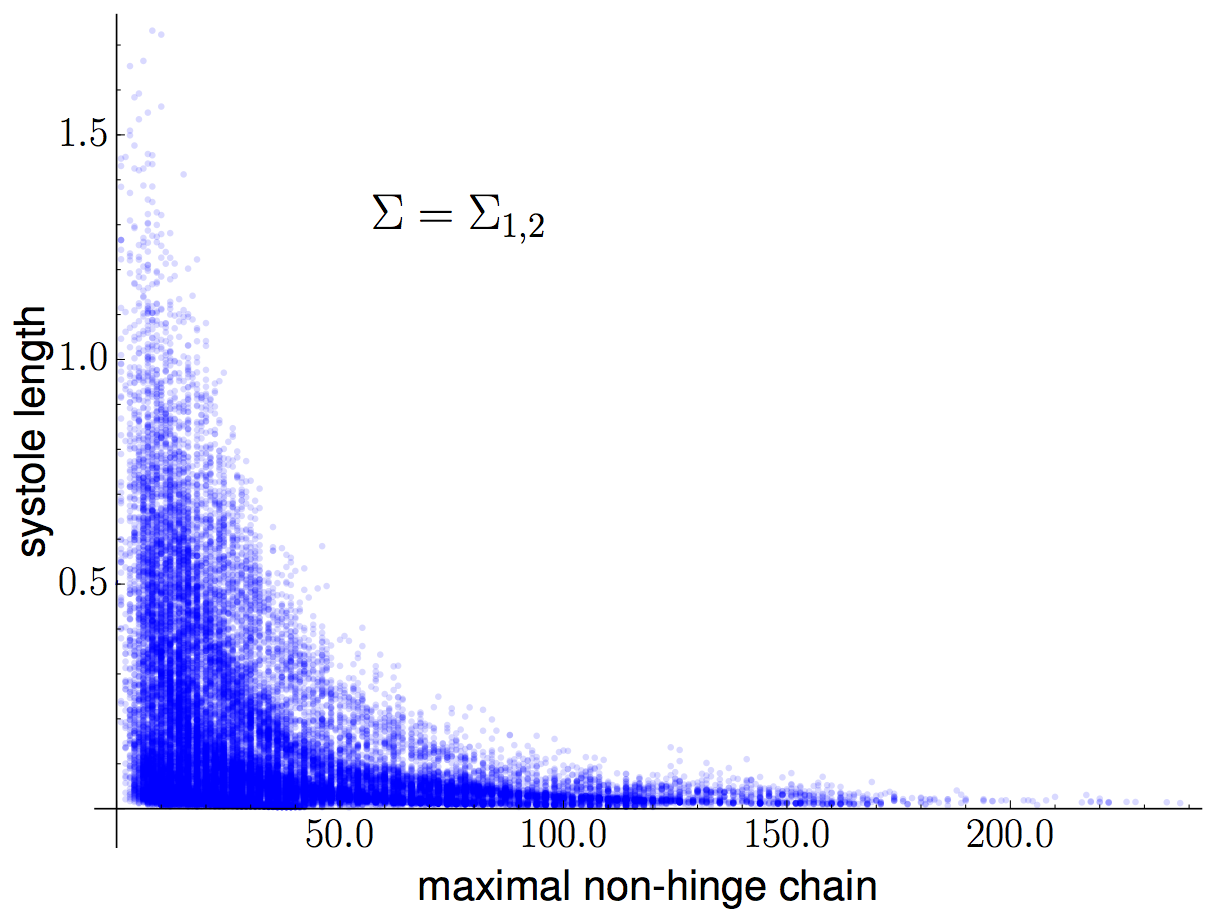}
        		\caption{}
        		\label{fig:chlen_sys_12_a}
       \end{subfigure}
        \begin{subfigure}[b]{0.49\textwidth}
                \includegraphics[scale=.38]{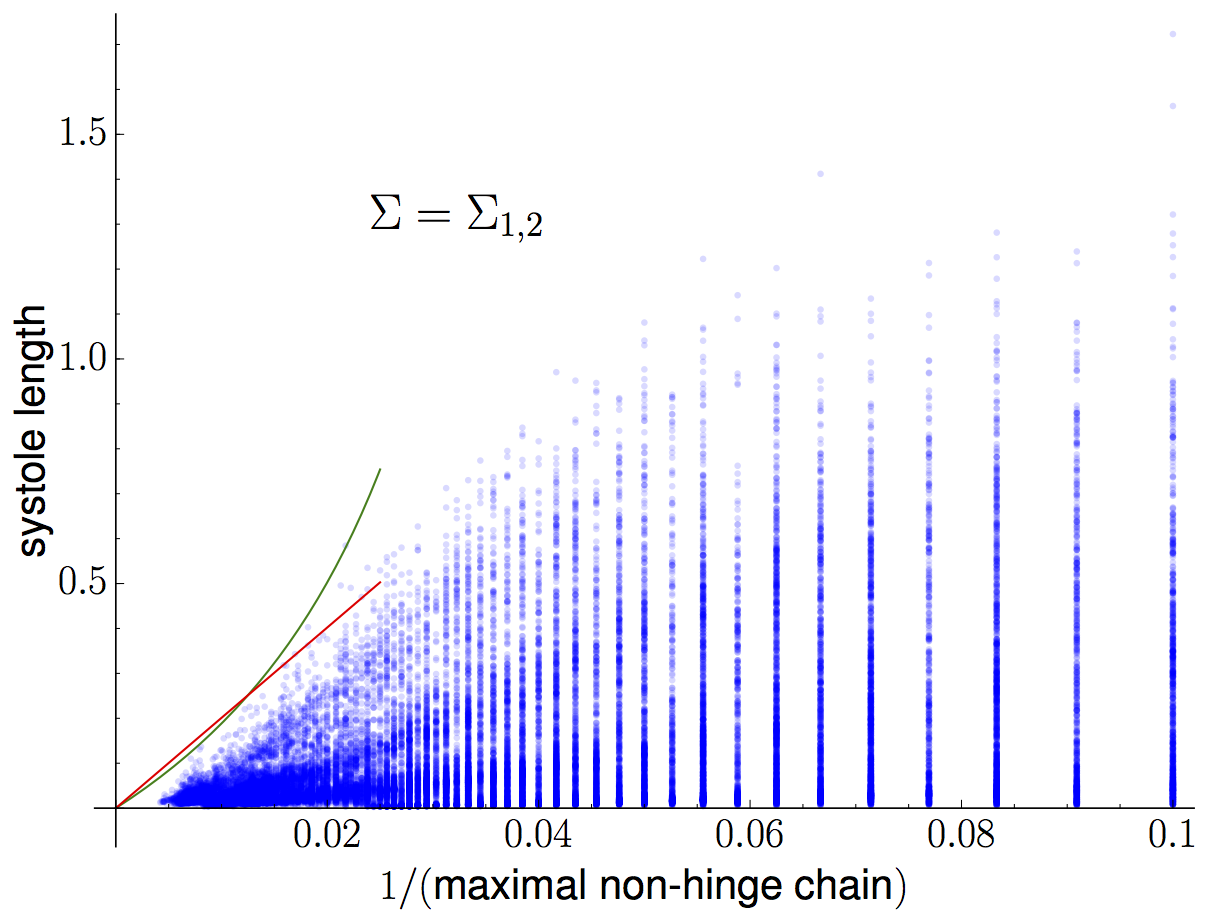}
                \caption{}
                \label{fig:chlen_sys_12_b}
        \end{subfigure}       
       \caption{The maximal chain of non-hinge tetrahedra (left) and its reciprocal (right) plotted against systole length, for $\sim 45000$ veering triangulations with associated surface $\Sigma_{0,5}$.}
        \label{fig:chlen_sys_12}
\end{figure}

Our experimental results, shown in Figures \Cref{fig:chlen_sys_05} and \Cref{fig:chlen_sys_12} for $\Sigma=\Sigma_{0,5}$ and $\Sigma=\Sigma_{1,2}$, respectively, demonstrate the upper bound on the systole length $l_\varphi(\sigma)$ in terms of the maximal non-hinge chain, as established in the corollary. We warn the reader that there is a caveat to these results: our sample is not truly random, due to computational limitations. See the note on methodology at the end of this section for a full explanation.

Figures \Cref{fig:chlen_sys_05_a} and \Cref{fig:chlen_sys_05_b} show our results for $\Sigma=\Sigma_{0,5}$ in two different formats. The first shows the maximal length of a chain of non-hinge tetrahedra on the $x$-axis, plotted against the systole length on the $y$-axis. In the second we have inverted the $x$-values, so that $1/|M_\mr{ch}|$ is plotted against systole length. In the latter plot, the $x$-axis has been cut off at $0.1$, in order to avoid compressing the left side of the plot too much. Figures \Cref{fig:chlen_sys_12_a} and \Cref{fig:chlen_sys_12_b} show the analogous scatter plots for $\Sigma=\Sigma_{1,2}$. In both cases the second plot suggests a coarse linear bound on systole length in terms of $\frac{1}{|M_\mr{ch}|}$, which is in agreement with the corollary. Consistent with the above observation about the absence of a converse to the corollary, we have many examples for which the systole is very small, but the maximal non-hinge chain length is not large. 

Note that in terms of the constant $D$ from the Length Bound Theorem, we can take $C=4\pi D\cdot |\chi|$. In principle, then, we should be able to get a lower bound on $D$ from our experimental data. The red line in \Cref{fig:chlen_sys_05_b} is the result of taking $D=0.25$. In fact if we use the (better) bound $l_\varphi(\sigma)\le \frac{\pi D\cdot |\chi|}{|M_\mr{ch}|-|\chi|^2-|\chi|}$ from the proof of the corollary, we get the green curve with $D=0.7$. Similarly, for the red and green curves in \Cref{fig:chlen_sys_12_b} we get $D=.4$ and $D=1.2$, respectively. The caveat to this observation is that, for it to be sound, we would need to know the size of the constant $N$ (and therefore $E$), since the corollary is only established for $|M_\mr{ch}|\ge N$.

Our definition of $M_\mr{ch}$ was chosen in part because as defined it is computable, but this definition has a couple of shortcomings. First, the total ordering on the tetrahedra is not canonical unless the splitting sequence for the veering triangulation has no multiple splits. In practice, though, this is not much of a problem. For $\Sigma_{0,5}$, of the $\sim 43000$ triangulations represented in \Cref{fig:chlen_sys_05}, only about $3\%$ have multiple splits in their train track splitting sequences, so for the vast majority the ordering on the tetrahedra is canonical. For $\Sigma_{1,2}$, there are always multiple splits according to our data, but the most branches we see splitting simultaneously is $4$. This means that in a different total ordering, $|M_\mr{ch}|$ could change by at most $4$, which is not very significant. 

The second drawback of our definition of $M_\mr{ch}$ is that it doesn't necessarily see full annular pockets. That is, we could have some number of non-hinges layer onto an annulus, then as the weight of the branches in this annulus decrease, some other branch far away could become maximal weight, resulting in a hinge being layered on before the remaining non-hinges are layered onto our original annulus. $M_\mr{ch}$ would then see this annular pocket as two separate non-hinge chains, while the full pocket will dictate the length of the core curve of the annulus. Again, though, for our results this is unlikely to play a large role. In the case of $\Sigma_{0,5}$, the fiber is most frequently $\Sigma_{0,6}$, for which all triangulations have $8$ triangles. Because the fiber is so small, it should be rare for a hinge to layer on disjointly from an annulus. The situation is  same for $S_{1,2}$, for which the fiber is typically $S_{1,4}$, so that triangulations again have $8$ triangles. 

\subsection*{A note on methodology} In order to compute the systole of a manifold, \texttt{SnapPy} must build the Dirichlet domain, which is computationally hard. For this reason, in this section we are forced to limit ourselves to short words, of length less than 100. On the other hand, if we sample $\mr{Mod}(\Sigma)$ by random walks, as we did in previous sections, then we do not get triangulations with long maximal chains. This is because long chains correspond, roughly, to large powers of Dehn twists in the mapping class word, since Dehn twisting many times in an annulus results in the lamination wrapping many times around the annulus. In order to get around this limitation, we include words of the form $w=w_0^kw_1w_2\dots w_m$ for various $k$ between $5$ and $40$, and words of the form $w=w_0^{k_0}w_1^{k_1}w_2^{k_2}\dots w_m^{k_m}$ for random positive integers $k_i$ selected from some interval. In both cases the $w_i$ are randomly selected from the set of generators discussed in \Cref{sec:method}.

\section{Other Results}
\label{sec:others}

\subsection{Geometric Veering Triangulations}
\label{sec:other:geo}

Despite the conjecture in \Cref{sec:nongeo}, it is still possible that there exist infinite families of mapping classes which yield geometric veering triangulations for surfaces other than $\Sigma_{1,1}$ and $\Sigma_{0,4}$. One approach to finding such a family, which was suggested to the author by Samuel Taylor, is to choose mapping classes $\varphi_k$ which have subwords $\psi^k$ which are supported on a subsurface $\Sigma_{1,1}\subset \Sigma$, and so that $\psi\in \mr{Mod}(\Sigma_{1,1})$ is not too complicated. The rough idea here is that when $k$ is large, the portion of the mapping torus coming from $\psi^k$ acting on $\Sigma_{1,1}$ should be (locally) geometrically similar to the triangulation of the mapping torus $M_{\psi^k}$ with fiber $\Sigma_{1,1}$. Since layered veering triangulations of once-punctured torus bundles are always geometric, we then have reason to hope that the layered veering triangulation of $M_{\varphi^\circ}$ will also be geometric.

\begin{figure}
        \centering
        \begin{subfigure}[b]{0.5\textwidth}
                \includegraphics[scale=.37]{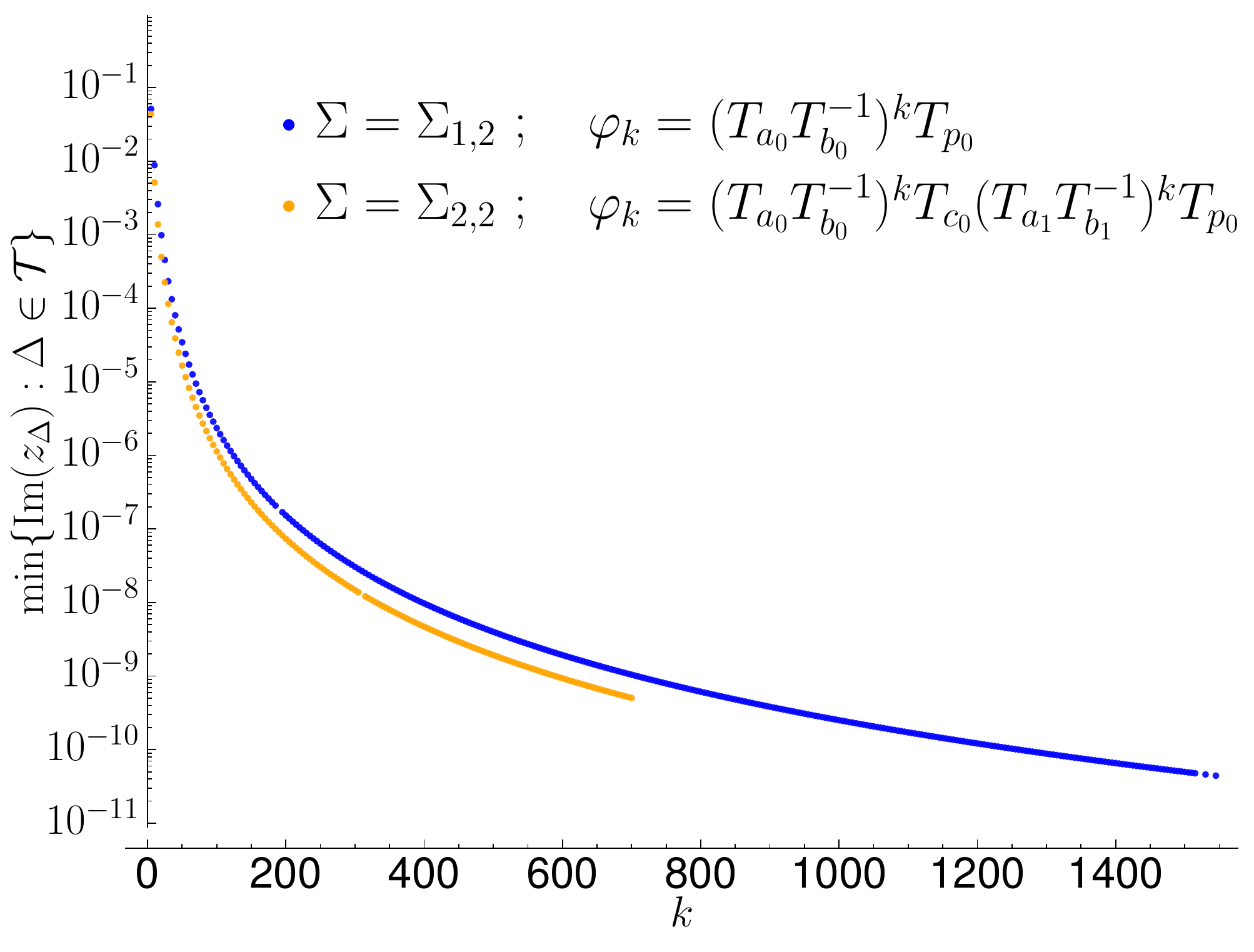}
        		\caption{}
        		\label{fig:min_tet_shape}
       \end{subfigure} 
       \begin{subfigure}[b]{0.49\textwidth}
                \includegraphics[scale=.37]{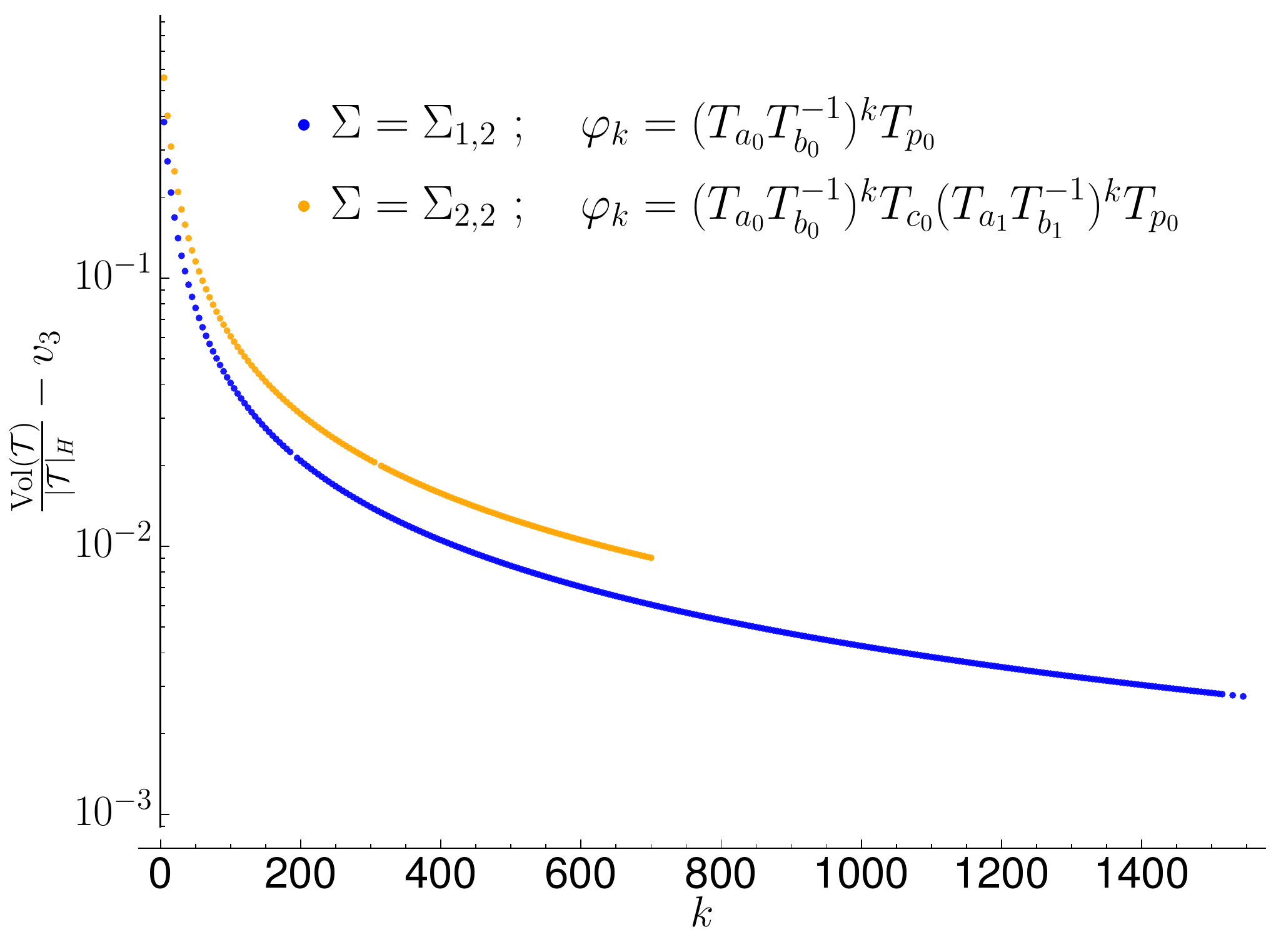}
        		\caption{}
        		\label{fig:casson_conv}
       \end{subfigure}  
       \caption{Left: This plot shows that as $k$ grows, the tetrahedron in $\mc{T}$ with smallest imaginary part becomes flatter in a very controlled way, but remained positive for $k$ up to about 1500. Right: The difference between $\frac{\mr{Vol}(\mc{T})}{|\mc{T}|_H}$ and the volume $v_3$ of a regular ideal tetrahedron becomes small as $k$ becomes large.}
        \label{fig:casson}
\end{figure}

We tested the two simplest such families, for $k$ as large as computationally possible, with promising results. For $\Sigma=\Sigma_{1,2}$ and $\varphi=(T_{a_0}T_{b_0}^{-1})^kT_{p_0}$ we were able to compute veering triangulations for $k$ up to $~1500$, and \texttt{SnapPy} reports that every one of these is geometric. As we can see in \Cref{fig:min_tet_shape}, which plots $k$ against the minimum of the imaginary parts of the tetrahedra shapes, as $k$ gets large some of the tetrahedra become very thin. For manifolds of this size, \texttt{SnapPy} is not able to rigorously verify the solutions to the tetrahedra gluing equations, so we cannot be sure that these tetrahedra shapes are accurate. However, we computed these examples with 212 bits of precision (64 decimal places), and the accuracy that \texttt{SnapPy} estimates for the tetrahedra shape is in all cases greater than $52$ decimal places. On the other hand, as we can see in \Cref{fig:min_tet_shape}, the flattest tetrahedron has imaginary part on the order of $10^{-11}$, which is quite large compared to $10^{-52}$.

Besides their apparent geometricity, these example have some other interesting properties. Recall from \Cref{sec:volume} that the conjectured lower bound for $\mr{Vol}(\mc{T})/|\mc{T}|_H$ is $v_3=1.0149416064...$, the volume of a regular ideal tetrahedron. \Cref{fig:casson_conv} shows $(\mr{Vol}(\mc{T})/|\mc{T}|_H-v_3)$  plotted against $k$, with a logarithmic $y$-axis, suggesting that as $k\to\infty$, this quantity becomes very small, possibly converging to $0$. In other words, $\mr{Vol}(\mc{T})/|\mc{T}|_H$ appears to be converging to the lower bound, or something very close to it.

Also somewhat interesting is the distribution of tetrahedron shapes for these triangulations. For the examples with $\Sigma=\Sigma_{1,2}$, $49.97\%$ of the tetrahedra in the sample are non-hinges, and of these $98.8\%$ have shapes that are within a distance of $\sim 0.023$ of either the origin or  $1+i0$---about half for each. Additionally, about $39.1\%$ of the hinges are within $\sim.023$ of $z_{\mr{reg}}$, the shape parameter of a regular ideal tetrahedron. Most of the remaining hinges also have shapes with imaginary parts at least $\frac{1}{2}$, so this is consistent with $\mr{Vol}(\mc{T})/|\mc{T}|_H$ being very close to $v_3$. For the example with $\Sigma=\Sigma_{2,2}$, the distribution of shapes is similar: $63.5\%$ are non-hinges, $97.1\%$ of which are concentrated near $0+i0$ and $1+i0$, and most of the hinges are close to $z_{\mr{reg}}$. This distribution of tetrahedra shapes is again consistent with $\mr{Vol}(\mc{T})/|\mc{T}|_H$ being very close to $v_3$.

Given the strength of the experimental results in this section, it is natural to ask whether these families of triangulations are in fact geometric. One way forward toward proving geometricity of these examples would be to try to parametrize the angle structure space of the triangulation corresponding to $\varphi_k$, in hopes that a maximum of the volume functional could be shown to exist in the interior of the space. Our results also raise the question of whether there are other families of geometric veering triangulations. Actually, it is likely that the construction used for the genus 1 and 2 examples given here could be generalized to higher genus, though we have not pursued this. Somewhat more ambitiously, we could ask: are there purely pseudo-Anosov subgroups of the mapping class group for which all associated veering triangulations are geometric? Or, in a different direction, are there families of veering triangulations that are canonical (apart from the once-punctured torus bundles and four-punctured sphere bundles)?

\subsection{Exponential Growth of Dilatation in Volume}

\begin{figure}
        \centering
        \begin{subfigure}[b]{0.5\textwidth}
                \includegraphics[scale=.37]{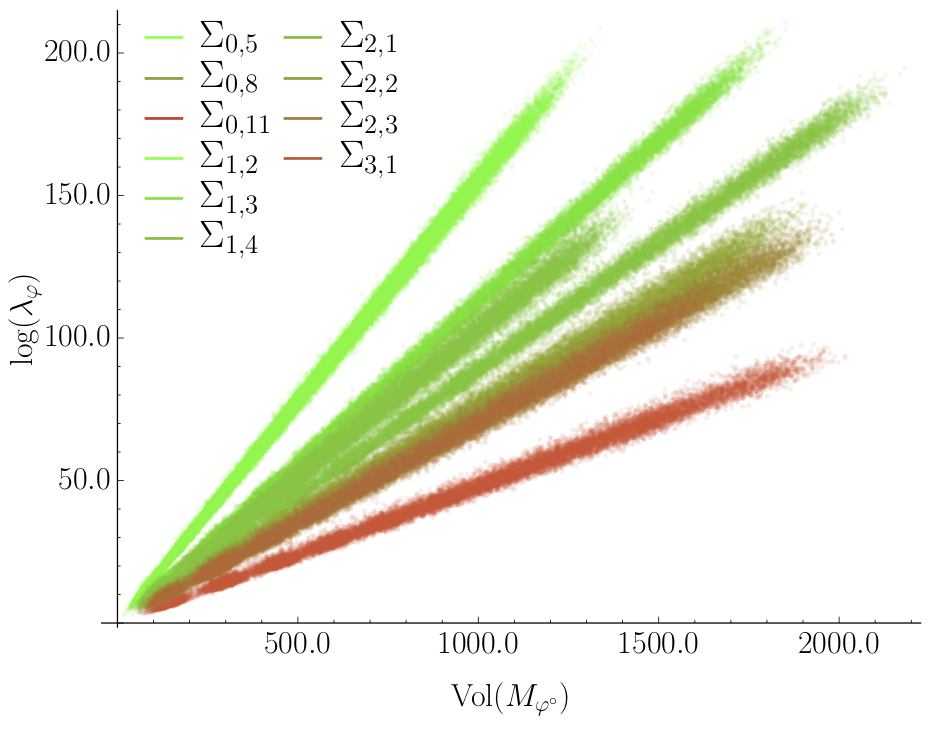}
        		\caption{}
        		\label{fig:dilatation}
       \end{subfigure} 
       \begin{subfigure}[b]{0.49\textwidth}
                \includegraphics[scale=.5]{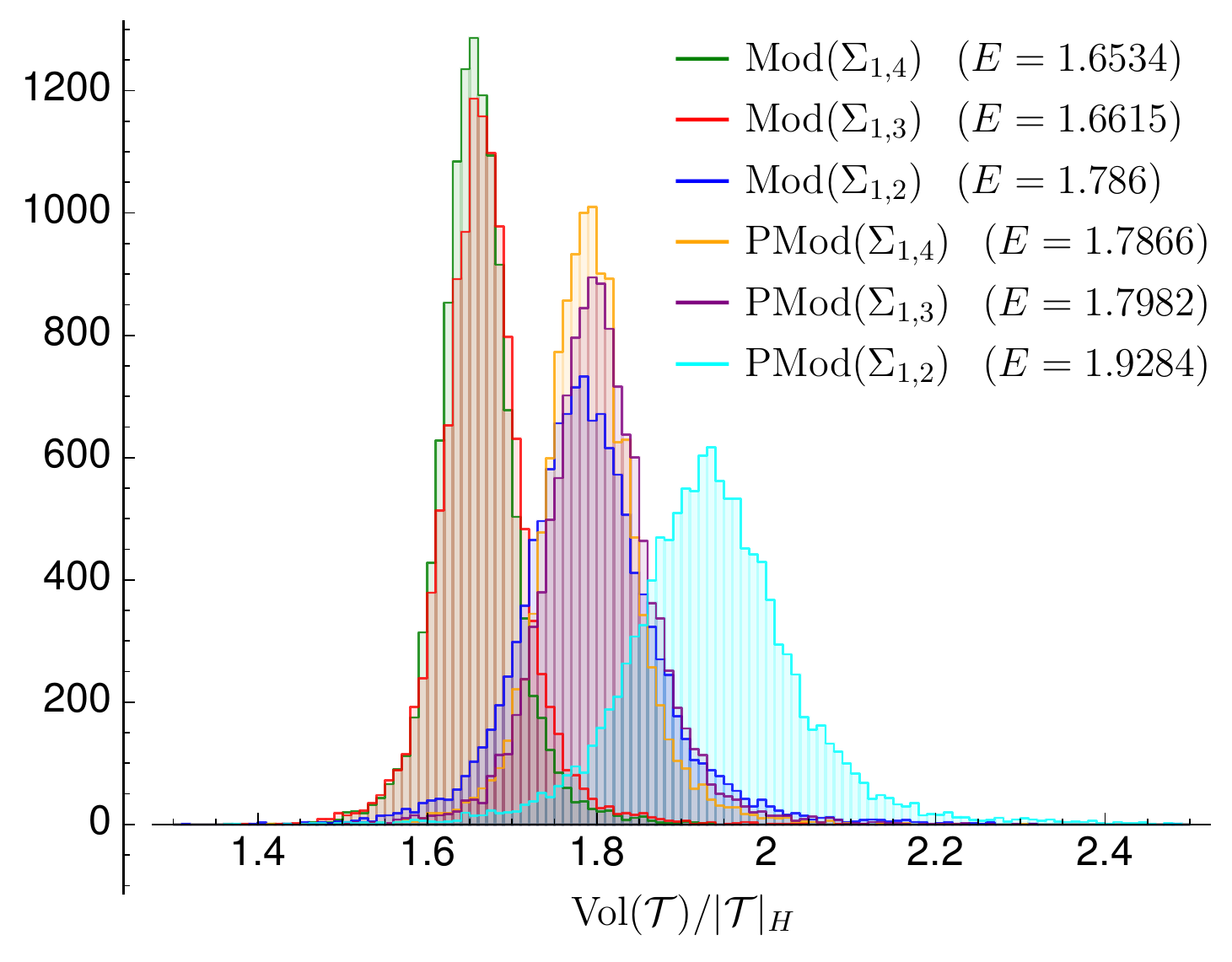}
                \caption{}
                \label{fig:hingevol_hist_g1}
        \end{subfigure}      
       \caption{Left: $\log$ of dilatation plotted against volume. The colors move toward the red end of the spectrum as surface complexity increases. Right: If we sample from the pure mapping class group $\mr{PMod}(\Sigma)$, the mean of $\mr{Vol}(\mc{T})/hinge(\mc{T})$ increases significantly.}
        \label{fig:dil_&_pmod}
\end{figure}

It is a result of Kojima and McShane that the dilatation $\lambda_\varphi$ of a mapping class $\varphi:\Sigma\to \Sigma$ satisfies $\log(\lambda_\varphi) \ge \frac{\mr{Vol}(M_\varphi)}{3\pi|\chi(\Sigma)|}$, and Kin--Kojima--Takasawa give a lower bound $\log(\lambda_\varphi)\le C_{\epsilon,\Sigma}\mr{Vol}(M_\varphi)$, where $\epsilon$ is the injectivity radius of $M_\varphi$.  For generic mapping classes, Rivin shows that the ratio $\frac{\log(\lambda_\varphi)}{\mr{Vol}(M_\varphi)}$ lies, with probability approaching $1$ as the word length $l(\varphi)\to \infty$, in some interval $[C_1,C_2]$ depending on $\Sigma$, where $C_1>0$. In addition, Rivin conjectures that $C_1=C_2$. Past experiments showing the linear relation between volume and log dilatation, in particular those of Kojima and Takasawa, have demonstrated well the bounds on $\frac{\log(\lambda_\varphi)}{\mr{Vol}(M_\varphi)}$ for generic examples. These experiments focused on manifolds of small volume, however, and because of this limitation do not offer clear evidence for or against Rivin's conjecture. Since our dataset consists of manifolds with volume as large as $\sim 2500$, we are able to see a much clearer linear relation, and also begin to see support for Rivin's conjecture. In \Cref{fig:dilatation}, each scatter plot point is a mapping torus $M_{\varphi_\circ}$ with fiber $\Sigma^\circ$, and is colored according to the complexity of the filled fiber $\Sigma$. In particular, those plot points with the highest complexity filled fiber are colored red, and the coloring transitions to green as complexity decreases. This coloring of the scatter plots reveals an apparent reciprocal relation between slope and complexity, which is perhaps not too surprising given the result of Kojima--McShane.

\subsection{Sampling from $\mr{PMod}(\Sigma)$ vs. $\mr{Mod}(\Sigma)$}

If $\Sigma$ has genus $g\ge 1$ and at least $n\ge2$ punctures, then the pure mapping class group $\mr{PMod}(\Sigma)$ will be a proper index $n$ subgroup of $\mr{Mod}(\Sigma)$, and we should expect that roughly one $n^{th}$ of our randomly sampled elements of $\mr{Mod}(\Sigma)$ will be pure mapping classes. In fact this is the case for our dataset, as we have checked, and if we plot separate tetrahedron shape histograms for elements of $\mr{PMod}(\Sigma)$ and $\mr{Mod}(\Sigma)\setminus\mr{PMod}(\Sigma)$, we find that the two histograms are roughly the same. 
On the other hand, if we sample directly from $\mr{PMod}(\Sigma)$, using as generators the Dehn twists about the $a_i, b_i, c_i$, and $p_i$ curves ((i.e., we omit twists about the $q_i$ curves in \Cref{fig:surf_gens} from our generating set), we find that the distribution of tetrahedra shapes is quite dramatically different from the shape distribution for $\mr{Mod}(\Sigma)$.

Similarly, \Cref{fig:hingevol_hist_g1} shows histograms of $\mr{Vol}(\mc{T})/|\mathcal{T}|_H$ for three different surfaces, and for random walks sampled from $\mr{PMod}(\Sigma)$, as well as the corresponding histograms for random walks sampled from $\mr{Mod}(\Sigma)$ (these latter appear in \Cref{fig:hingevol_hist_all} as well). Notice that the expected values for $\mr{PMod}(\Sigma)$ is greater than that of $\mr{Mod}(\Sigma)$ for each surface, by about $.14$. As with tetrahedra shape histograms above, this difference of expected value has more to do with difference in generating sets than the fact that we are comparing $\mr{PMod}(\Sigma)$ elements to $\mr{Mod}(\Sigma)$ elements. For example, if we sample $\mr{Mod}(\Sigma_{1,2})$ using the full generating set, half the elements will be pure mapping classes, but the histogram of $\mr{Vol}(\mc{T})/|\mathcal{T}|_H$ for these pure mapping classes will be virtually identical to that of the non-pure elements.
       
\section{Supplemental Materials}
\label{sec:supp}

The results of this paper are based on analysis of approximately $766,000$ veering triangulations. Our full data set, along with all Python scripts used in the process of generating and analyzing the data, plus a number of additional figures which do not appear in this paper, will be available for download from the author's web page, currently at \url{www.wtworden.org/research/esvt/}. Each of these triangulations is stored as part of a larger Python object, which contains a wealth of other computed information, including all invariants computed by \texttt{flipper} when the veering triangulation was computed, and several which were computed by \texttt{SnapPy}, along with a method to call the full \texttt{SnapPy} manifold. Since these objects are rather large, we also store a Python dictionary object which contains all of the more compact information of interest, without the full \texttt{flipper} triangulation object. While the full dataset including \texttt{flipper} triangulations is about 522 gb uncompressed, the dictionaries alone (plus log files) are only about 41 gb uncompressed. For this reason only the dictionaries are downloadable from our webpage at this time. Using a module called \texttt{vt\_tools}, all or a subset of these dictionaries can be loaded into a Python session, and can be queried with respect to any of the dictionary keys (for example, you could ask for all triangulations with $\mr{Vol}(\mc{T})<100$ and associated surface $\Sigma=\Sigma_{2,3}$). This \texttt{vt\_tools} module also contains all functions used to generate the figures in this paper (most of these require \texttt{Sage}, though). More detailed instructions for navigating the data set and querying tools will be provided on our website.

\bibliographystyle{alpha}
\bibliography{biblio}

\newcommand{\etalchar}[1]{$^{#1}$}
\begin{thebibliography}{CDGW17}

\bibitem[Ago11]{Ag}
Ian Agol.
\newblock Ideal triangulations of pseudo-{A}nosov mapping tori.
\newblock In {\em Topology and geometry in dimension three:}, volume 560 of
  {\em Contemp. Math.}, pages 1--17. Amer. Math. Soc., Providence, RI, 2011.

\bibitem[Aki99]{Ak}
Hirotaka Akiyoshi.
\newblock On the {F}ord domains of once-punctured torus groups.
\newblock {\em S{\=u}rikaisekikenky{\=u}sho K{\=o}ky{\=u}roku}, 1104:109--121,
  1999.

\bibitem[BCM12]{BCM}
Jeffrey~F. Brock, Richard~D. Canary, and Yair~N. Minsky.
\newblock The classification of {K}leinian surface groups, {I}{I}: The ending
  lamination conjecture.
\newblock {\em Ann. of Math.}, 176(1):1--149, 2012.

\bibitem[Bel16]{Be}
Mark Bell.
\newblock flipper (computer software).
\newblock Available at {\url{pypi.python.org/pypi/flipper}}, 2013--2016.
\newblock Version 0.9.8.

\bibitem[Bro03]{Br}
Jeffrey Brock.
\newblock Weil-{P}etersson translation distance and volumes of mapping tori.
\newblock {\em Communications in Analysis and Geometry}, 11(5):987--999, 2003.

\bibitem[CD06]{CaDi}
James~W. Cannon and Warren Dicks.
\newblock On hyperbolic once-punctured-torus bundles {I}{I}: fractal
  tessellations of the plane.
\newblock {\em Geometriae Dedicata}, 123(1):11--63, 2006.

\bibitem[CDGW17]{SnapPy}
Marc Culler, Nathan~M. Dunfield, Matthias Goerner, and Jeffrey~R. Weeks.
\newblock Snap{P}y, a computer program for studying the geometry and topology
  of $3$-manifolds.
\newblock Available at \url{http://snappy.computop.org} (Version 2.3.2),
  2009--2017.

\bibitem[DHO{\etalchar{+}}14]{Du:knots}
Nathan~M. Dunfield, A.~Hirani, M.~Obeidin, A.~Ehrenberg, S.~Battacharyya,
  D.~Lei, et~al.
\newblock Random knots: A preliminary report.
\newblock Available at
  \url{http://www.math.uiuc.edu/~nmd/slides/random_knots.pdf}, 2014.

\bibitem[DS10]{DiSa}
Warren Dicks and Makoto Sakuma.
\newblock On hyperbolic once-punctured-torus bundles {I}{I}{I}: Comparing two
  tessellations of the complex plane.
\newblock {\em Topology and its Applications}, 157(12):1873--1899, 2010.

\bibitem[DT03]{DuTh:haken}
Nathan~M. Dunfield and William~P. Thurston.
\newblock The virtual {H}aken conjecture: Experiments and examples.
\newblock {\em Geom. Topol.}, 7:399--441, 2003.

\bibitem[DT06]{DuTh:tunnel}
Nathan~M. Dunfield and Dylan~P. Thurston.
\newblock A random tunnel number one $3$-manifold does not fiber over the
  circle.
\newblock {\em Geom. Topol.}, 10:2431--2499, 2006.

\bibitem[FG11]{FuGu:survey}
David Futer and Fran{\c{c}}ois Gu{\'e}ritaud.
\newblock From angled triangulations to hyperbolic structures.
\newblock {\em Contemporary Mathematics}, 541:159--182, 2011.

\bibitem[FG13]{FuGu:veering}
David Futer and Fran{\c{c}}ois Gu{\'e}ritaud.
\newblock Explicit angle structures for veering triangulations.
\newblock {\em Algebraic \& Geometric Topology}, 13(1):205--235, 2013.

\bibitem[FH82]{FlHa}
William Floyd and Allen Hatcher.
\newblock Incompressible surfaces in punctured-torus bundles.
\newblock {\em Topology and its Applications}, 13(3):263--282, 1982.

\bibitem[GM17]{GaMa}
Vaibhav Gadre and Joseph Maher.
\newblock The stratum of random mapping classes.
\newblock {\em Ergodic Theory and Dynamical Systems}, pages 1--17, 2017.

\bibitem[Gu{\'e}06]{Gu:thesis}
Fran{\c{c}}ois Gu{\'e}ritaud.
\newblock G{\'e}om{\'e}trie hyperbolique effective et triangulations
  id{\'e}ales canoniques en dimension trois.
\newblock {\em PhD Thesis}, 2006.

\bibitem[Gu{\'e}16]{Gu:cannon}
Fran{\c{c}}ois Gu{\'e}ritaud.
\newblock Veering triangulations and {C}annon-{T}hurston maps.
\newblock {\em J. Topology}, 9(3):957--983, 2016.

\bibitem[HIK{\etalchar{+}}16]{hikmot}
Neil Hoffman, Kazuhiro Ichihara, Masahide Kashiwagi, Hidetoshi Masai,
  Shin’ichi Oishi, and Akitoshi Takayasu.
\newblock Verified computations for hyperbolic 3-manifolds.
\newblock {\em Experimental Mathematics}, 25(1):66--78, 2016.

\bibitem[HIS16]{HoIsSe}
Craig~D Hodgson, Ahmad Issa, and Henry Segerman.
\newblock Non-geometric veering triangulations.
\newblock {\em Experimental Mathematics}, 25(1):17--45, 2016.

\bibitem[HRST11]{HoRuSeTi:veering}
Craig~D Hodgson, J~Hyam Rubinstein, Henry Segerman, and Stephan Tillmann.
\newblock Veering triangulations admit strict angle structures.
\newblock {\em Geom. Topol.}, 15(4):2073--2089, 2011.

\bibitem[Iss12]{Is}
Ahmad Issa.
\newblock {\em Construction of non-geometric veering triangulations of fibered
  hyperbolic $3$-manifolds}.
\newblock PhD thesis, Master’s thesis, University of Melbourne, 2012.

\bibitem[J{\o}r03]{Jo}
Troels J{\o}rgensen.
\newblock On pairs of once-punctured tori.
\newblock In {\em Kleinian Groups and Hyperbolic 3-Manifolds}, volume 299 of
  {\em London Math. Soc. Lec. Notes}, pages 183--208. Cambridge University
  Press, 2003.

\bibitem[Lac04]{La}
Marc Lackenby.
\newblock The volume of hyperbolic alternating link complements.
\newblock {\em Proc. London Math. Soc. (3)}, 88(1):204--224, 2004.
\newblock With an appendix by Ian Agol and Dylan Thurston.

\bibitem[Mah10]{Mah}
Joseph Maher.
\newblock Linear progress in the complex of curves.
\newblock {\em Transactions of the American Mathematical Society},
  362(6):2963--2991, 2010.

\bibitem[Min10]{Mi:kleinian}
Yair Minsky.
\newblock The classification of {K}leinian surface groups, {I}: models and
  bounds.
\newblock {\em Ann. Math.}, 171(1):1--107, 2010.

\bibitem[MT17]{MiTa}
Yair~N Minsky and Samuel~J Taylor.
\newblock Fibered faces, veering triangulations, and the arc complex.
\newblock {\em Geom. Funct. Anal.}, 27:1450--1496, 2017.

\bibitem[Riv14]{Ri}
Igor Rivin.
\newblock Statistics of random $3$-manifolds occasionally fibering over the
  circle.
\newblock {\em arXiv preprint arXiv:1401.5736}, 2014.

\bibitem[ST17]{SiTa}
Alessandro~N Sisto and Samuel~J Taylor.
\newblock Largest projections for random walks and shortest curves in random
  mapping tori.
\newblock {\em Math. Res. Lett.}, 2017.
\newblock to appear.

\bibitem[Thu78]{Th}
William Thurston.
\newblock Geometry and topology of $3$-manifolds, lecture notes.
\newblock {\em Princeton University}, 1978.

\bibitem[Wee05]{We}
Jeffrey Weeks.
\newblock Computation of hyperbolic structures in knot theory.
\newblock In William Menasco and Morwen Thistlethwaite, editors, {\em Handbook
  of knot theory}, pages 461--480. Elsevier Science, 2005.

\end{thebibliography}

\end{document}